\documentclass[12equipped]{article}
\usepackage{amssymb,amsmath}

\input epsf.tex

\title{Airy structures and symplectic geometry of topological recursion}

\author {Maxim Kontsevich, Yan Soibelman}

\begin{document}

\maketitle

\newtheorem{thm}{Theorem}[subsection]
\newtheorem{defn}[thm]{Definition}
\newtheorem{lmm}[thm]{Lemma}
\newtheorem{rmk}[thm]{Remark}
\newtheorem{prp}[thm]{Proposition}
\newtheorem{conj}[thm]{Conjecture}
\newtheorem{exa}[thm]{Example}
\newtheorem{cor}[thm]{Corollary}
\newtheorem{que}[thm]{Question}
\newtheorem{ack}{Acknowledgements}
\newtheorem{exe}{Exercise}

\newcommand{\C}{{\bf C}}
\newcommand{\K}{{\bf k}}
\newcommand{\R}{{\bf R}}
\newcommand{\N}{{\bf N}}
\newcommand{\Z}{{\bf Z}}
\newcommand{\Q}{{\bf Q}}
\newcommand{\GG}{\Gamma}
\newcommand{\A}{A_{\infty}}
\newcommand{\g}{{\mathfrak{g}}}
\newcommand{\CC}{{\mathcal C}}
\newcommand{\MM}{{\mathcal M}}
\newcommand{\LL}{{\mathcal L}}
\newcommand{\FF}{{\mathcal F}}
\renewcommand{\H}{{\mathcal{H}}}
\newcommand{\BB}{{\mathcal{B}}}
\newcommand{\epi}{\twoheadrightarrow}
\newcommand{\mono}{\hookrightarrow}
\newcommand\ra{\rightarrow}
\newcommand\uhom{{\underline{Hom}}}
\newcommand\OO{{\mathcal O}}
\newcommand{\WS}{\widehat{S}}
\newcommand{\x}{\widehat{x}}
\newcommand{\y}{\widehat{y}}
\newcommand{\WW}{\mathcal{W}}
\newcommand{\D}{\mathbb{D}}

\newcommand{\epp}{\varepsilon}

\newcommand{\M}{{\mathsf{M}}}

\newcommand{\G}{{\mathsf{G}}}

\newcommand{\Gr}{{\mathsf{Gr}}}
\newcommand{\T}{{\mathsf{T}}}

\tableofcontents

\section{Introduction}

\subsection{Linear algebra data of the topological recursion}

Conventional approach to topological recursion (TR for short, see [EynOr1]) produces
a collection of meromorphic polydifferentials $\omega_{g,n}$ on
a spectral curve $\Sigma$ starting with the ``initial data"
$\omega_{0,1},\omega_{0,2}$. 
From $\omega_{0,1},\omega_{0,2}$ one derives a certain integral kernel $K(z,w)$, which behaves like a differential form in 
one variable and vector field in the other one.
Forms $\omega_{g,n}$ except of $\omega_{0,1},\omega_{0,2}$ 
belong to the symmetric powers $Sym^n(V)$ of a certain infinite-dimensional vector space $V$ of meromorphic $1$-forms. The form $\omega_{0,1}$ is auxiliary and is used only via so-called recursion kernel. The form $\omega_{0,2}$ is not an honest element of $Sym^2(V)$. 

 In this paper we suggest to change the set of initial data, by replacing the above forms by four tensors assigned to $V$, which we denote by  $A,B,C,\epsilon$. These tensors can be derived from
$\omega_{0,1},\omega_{0,2}$. 
The tensor $A$ is essentially $\omega_{0,3}$, and the tensor $\epsilon$ is  $\omega_{1,1}$. There are explicit formulas for the tensors $B$ and $C$ which we will discuss in the paper.
The vector space $V$ is the space of such meromorphic $1$-forms on the spectral curve which have trivial residues at the ramification points.

From this new point of view the topological recursion becomes a special case of the deformation quantization. More precisely it gives a cyclic vector
(wave function) of a deformation quantization module ($DQ$-module for short) corresponding to a certain quadratic\footnote{Here we mean that  the subvariety under consideration  is the set of common zeroes of a collection of polynomials of degree at most 2.}
Lagrangian subvariety
in the infinite-dimensional symplectic vector space
$V\oplus V^{\ast}$.  The Lagrangian subvariety can be described explicitly in terms of the tensors $A,B,C,\epsilon$.
The construction is purely algebraic and it does not depend on the notion of  spectral curve. Furthermore, it makes $V$ into an algebra over a certain PROP defined by the above four tensors.

From this perspective spectral curves are not fundamental objects in the formalism of topological recursion,
contrary to the conventional point of view.

In particular, the vector space $V$ can be arbitrary (e.g. finite-dimensional or infinite-dimensional).\footnote{In [ABCO] the authors study the dynamics of Young diagrams via an ``abstract''  version of topological recursion which seems to 
give an Airy structure in our sense.}

\subsection{Lie algebra of the topological recursion}

The notion of {\it Airy structure} underlying our approach to TR is named after the archetypical example of 
the ``Airy spectral curve'' $x=y^2$.
In fact there are two notions: {\it classical Airy structure} and {\it quantum Airy structure}.

Starting with the tensors $A,B,C$ we define a {\it classical} Airy structure. It  is equivalent to the above-mentioned quadratic Lagrangian subvariety. Skew-symetrization of the tensor $C$ endows $V$ with a Lie coalgebra structure, hence a Lie algebra structure on $\g:=V^{\ast}$.

The tensor $\epsilon$  is responsible for the quantization of the classical Airy structure. It gives a trivialization (i.e. a way to represent the cocycle as a coboundary)
of a certain $2$-cocycle on the Lie algebra $\g$. 
The cocycle itself is a restriction to the above Lie algebra of the famous ``Japanese cocycle'' for the Lie algebra $gl(V\oplus V^{\ast})$.
In the case when $dim\,V<\infty$ the restricted cocycle is naturally a coboundary. In the infinite-dimensional case (which we
are mostly interested in) a choice of trivializing $1$-cochain $\epsilon$ is an additional piece of information.
In the case of  the Airy spectral curve one finds that $\g\simeq \C[[x]]\,\partial_x\ltimes (\C[[x]]/\C)$.

\subsection{Spectral curves and the corresponding Airy structures}

Our notion of spectral curve is more general than the one in the conventional TR. 
Namely, we consider a smooth complex Poisson surface $P$ (which can be non-compact),
such that the divisor $D$ of zeros of the Poisson bivector field is smooth. We assume that the open symplectic leaf $S:=P-D$
is endowed with a Lagrangian foliation $\FF$, which extends by continuity  to $P$ and the extension is tangent to $D$. Then an $(\FF,D)$-transversal spectral curve is 
a smooth compact curve $\Sigma\subset P$  which intersects transversally $D$ at finitely many points $p_\beta$ with multiplicities $k_\beta\ge 1$, and which is tangent with order $1$ to leaves of the foliation 
$\FF$ at finitely many points $r_\alpha$.

These data determine a germ of complex analytic space ${\mathcal B}_1$, which is the moduli space of deformations of $\Sigma$ in the class of spectral curves. It is a subspace of a bigger moduli space $\BB_0$ of spectral curves defined in the same way, but without the tangency condition.

Consider the family of  symplectic vector spaces 
\[{\mathcal H}_\Sigma:={\mathbb H}^1(\Sigma, {\Omega}^0_\Sigma(-\sum_\beta k_\beta p_\beta)\to {\Omega}^1_\Sigma(\sum_\beta k_\beta p_\beta))\,.\] 
We prove that the family ${\mathcal H}_\Sigma$ gives rise to a symplectic vector bundle on ${\mathcal B}_0$, endowed with an {\it affine} flat connection.
It follows that (locally) ${\mathcal B}_0$ can be embedded into any affine symplectic vector space ${\mathcal H}_\Sigma$  as a Lagrangian subvariety. In other words we have a (germ of) Lagrangian submanifold in an affine symplectic space.

Then the general story of deformation quantization  ensures that there is  a cyclic $DQ$-module which quantizes the Lagrangian image 
of ${\mathcal B}_1$. This cyclic $DQ$-module can be realized as a representation of the (completed) Weyl algebra of ${\mathcal H}_\Sigma$. Locally (along the Lagrangian germ) it is unique up to isomorphism. Nevertheless there is no preferable choice of a cyclic vector.

It turns out that our quantum Airy structure gives a canonical cyclic vector $\psi_{\BB_1}$ (wave function). In fact it exists even for a bigger space $\BB_0$.
More precisely, the finite-dimensional simplectic vector space ${\mathcal H}_\Sigma$ can be represented as the quotient $G/G^\perp$, where $G$ is an infinite-dimensional  coisotropic subspace of a symplectic vector space $W$ of meromorphic $1$-forms with trivial residue on the disjoint union of formal discs about $r_\alpha$ (here $G^\perp$ is the skew-orthogonal to $G$). Then $\BB_1$ can be obtained as a symplectic (a.k.a Hamiltonian) reduction
of an infinite-dimensional quadratic Lagrangian submanifold ${ L}\subset W$, i.e. $\BB_1$ is isomorphic to the image of ${L}\cap G$ under the natural projection $G\to G/G^\perp={\mathcal H}_\Sigma$. The Lagrangian submanifold ${L}$ gives rise to a classical Airy structure on any
Lagrangian complement to a  tangent space to ${L}$. 

 Any Lagrangian complement ${\mathcal V}_\Sigma$ to the finite-dimensional Lagrangian subspace $T_\Sigma \BB_0\subset \H_\Sigma$ gives (after the pullback by epimorphism  $G\epi \H_\Sigma$) an infinite-dimensional Lagrangian complement $V_\Sigma$ to the tangent space to $L$. 
 This gives us a classical Airy structure on ${V}_\Sigma$, which can be quantized. The corresponding wave function $\psi$ is given by the TR, and it corresponds to $\psi_{\BB_1}$. Notice that the choice of $V_\Sigma$ corresponds to a choice of the Bergman kernel in the language of conventional TR.

\subsection{Topological recursions}

Another interesting observation is that there are different versions of TR.
In particular, we define two sets of tensors: the set $(A_{TR}, B_{TR}, C_{TR}, \epsilon_{TR})$ corresponding to the conventional (Eynard-Orantin) version of the TR, and the set $(A,B,C,\epsilon)$ mentioned above, which is based on certain residue constraints arising in the deformation theory of the spectral curve $\Sigma$. These two sets of tensors are different. Furthermore, the set $(A_{TR}, B_{TR}, C_{TR}, \epsilon_{TR})$ does not satisfy the axiomatics of Airy structures. Nevertheless the comparison of the two approaches can be made. For that we  introduce a certain subspace $V_{odd}\subset V$ respected by  both sets of tensors. Furthermore they define on it the same Airy structure (which is a natural Airy substructure of the one given by $(A,B,C,\epsilon)$). This leads to a new recursion formulas, which are different from those of Eynard-Orantin, but produce the same collection of polydifferentials $\omega_{g,n}, n\ge 1, 2-2g-n<0$.

Moreover, our modified recursion is   in a sense simpler than  the original Eynard-Orantin TR. Our formulas do not involve local involutions (compare Sections \ref{section: TR} and \ref{sec: modified}).

\subsection{Contents of the paper}
Classical and quantum Airy structures are introduced in Section 2.
We also prove there the result that connects TR with deformation quantization. The proof is purely algebraic, so it works for Airy structures 
over an arbitrary field of characteristic zero. An important corollary of our formalism is the fact that the tensors arising in this ``abstract'' topological recursion are symmetric.

In Section 3 we consider the local version of TR.
First we recall the basics on the conventional TR and derive the corresponding tensors $(A_{TR}, B_{TR}, C_{TR}, \epsilon_{TR})$.
As we have already pointed out they do not satisfy the axioms of the Airy structure, contrary to the alternative approach via the residue constraints.
We explain in Section 3 that two approaches agree if we restrict both tensors to the smaller subspace $V_{odd}$.

Finally we compare two approaches in general. The comparison reveals the role of Bergman kernel and the notion of primitive Airy structure which we introduced in Section 2.

In Section 4 we discuss the notion of spectral curve as well the deformation theory of spectral curves. It turns out that the natural framework for that is the one of a complex Poisson 
surface which has an open symplectic leaf endowed with an analytic foliation.
This allows us to consider e.g. non-compact spectral curves, like those which are associated with irregular Hitchin systems (see e.g. [KoSo]).

In  Section 5 we discuss the notion of affine symplectic connection which allows us to embed the germ of the moduli space of spectral curves into an affine symplectic space as a germ of Lagrangian submanifold.

In Section 6 we discuss an infinite-dimensional version of the material of Section 5 in the case of spaces of embedded formal discs. There we use the ideas of formal differential geometry of Gelfand and Kazhdan (see [GelKaz]). 

In Section 7 we explain  the formal version of the embedded Lagrangian germ. It is described in terms of the Hamiltonian reduction of a certain infinite-dimensional Lagrangian submanifold of a Tate vector space with respect to a coisotropic subspace. In that construction only the coistoropic space depends on the global geometry of the spectral curve. Rest of the data are local: they are derived from the formal neighborhoods of points, where the spectral curve is tangent to the foliation. An important corollary is that $TR$ solves the so-called Holomorphic Anomaly Equation, i.e. quantizes the Lagrangian image of $\BB_0$ in a covariantly constant way.

In Section 8 we discuss a possible extension of the quantization of the moduli space of smooth spectral curves to the case when there is no foliation at all. In the case when the foliation exists, the moduli space of smooth spectral curves can be embedded locally as a Lagrangian submanifold in an affine symplectic space (e.g. in examples it can be an open dense subvariety of the base of Hitchin integrable system). In the case when there is no foliation, the moduli space still can be embedded locally as a Lagrangian submanifold,
but the ambient symplectic manifold does not have a global affine structure. The ambient space has a weaker structure which we call a {\it semi-affine symplectic structure}.

In Section 9 we propose two speculations. One is about a possible relation of the quantization of the moduli space of smooth spectral curves (e.g. the above-mentioned open subset of the Hitchin base)  with the quantization of an individual spectral curve. Another speculation is about differential-graded Airy structures derived from the moduli spaces of {\it compact} 3-dimensional Calabi-Yau varieties.

Few final remarks.

We warn the reader that this paper is not about applications of the topological recursion. It is about the origin of the recursion relations.
For that reason we make only a basic comparison of the formalism of Airy structures with the one of Eynard-Orantin. We leave for the future work many interesting aspects of the story, which (as we hope) admit a natural explanation in the framework of Airy structures. Those include e.g. the relation
of the TR to WKB  quantization of an {\it individual} spectral curve (see [DuMu]),
symplectic invariance in TR (see [EynOr2]). We do not discuss the passage from the multi-differentials $\omega_{g,n}$ to the corresponding functions $F_{g,n}$. There are many other questions which arise naturally in the framework of our formalism, but which we prefer to leave for the future work.
This includes e.g. the generalization of  TR to non-smooth spectral curves and also to so-called quantum curves.

{\it Acknowledgments.} We thank to J. Andersen, G. Borot, L. Chekhov, B. Eynard, M. Mulase for discussions on  various aspects of topological recursion. Y.S. thanks to IHES for excellent research conditions. His work was partially supported by an NSF grant.

\section{Airy structures}

\subsection{Classical Airy structures}

Let $\K$ be a field of characteristic zero, $V$ an $n$-dimensional $\K$-vector space. 
Let us choose an ordered basis on $V$ and denote by $y_1,...,y_n$ the corresponding coordinates. Denote by $x^1,...,x^n$ the dual coordinates on $V^{\ast}=Hom_{\K}(V,\K)$.

On the vector space $W:=T^{\ast}(V^{\ast})=V\oplus V^{\ast}$  we have the standard symplectic structure. The corresponding Poisson bracket in the coordinates are given by the formulas $\{y_i,x^j\}=\delta_{ij},\,\,\{x^i,x^j\}=\{y_i,y_j\}=0$. 

We denote by $Sym_{\le 2}(W)$ the Lie algebra $Sym^0(W)\oplus Sym^1(W)\oplus Sym^2(W)$ of polynomial  functions
on $W^{\ast}$ of degree less or equal than $2$, endowed with the natural Lie algebra structure induced 
by the above Poisson bracket.

\begin{defn} Classical Airy structure on $V$ is given by a collections of at most  quadratic polynomials (i.e. elements of $Sym_{\le 2}(W)$) of the form
$$H_i:=-y_i+\sum_{jk}a_{ijk}x^jx^k+2\sum_{jk}b_{ij}^kx^jy_k+\sum_{jk}c_i^{jk}y_jy_k, 1\le i\le n,$$
such that the vector space $\oplus_{1\le i\le n}\,\K\cdot H_i$ is closed under the Poisson bracket.
\end{defn}

To save the language we will often skip the word ``at most'' and call such polynomials {\it quadratic}. The same terminology will apply to
the variety defined by equations $H_i=0, 1\le i\le n$.

In a coordinate-free form a
classical Airy structure on $V$ 
is defined by three tensors:
\begin{itemize}
\item[a)] tensor $A=(a_{ijk})\in V^{\otimes 3}$,
\item[b)] tensor $B=(b_{ij}^k): V\to V\otimes V$,
\item[c)] tensor $C=(c_i^{jk}): V\otimes V\to V$, 
\end{itemize}
such that 
$A$ and $C$  are  symmetric with respect to the permutation of the indices $j,k$.

Let us analyze the constraints on $A, B, C$.

By definition, any classical Airy structure gives rise to a Lie subalgebra
$$\g:=V^{\ast}\simeq \oplus_{1\le i\le n}\K\cdot H_i\subset Sym_{\le 2}( W.)$$
Let us compute its structure constants. We have
\begin{align*}\{H_{i_1},H_{i_2}\} &=\{-y_{i_1}+2\sum_{j,k}b_{i_1j}^kx^jy_k+...,-y_{i_2}+2\sum_{j,k}
b_{i_2j}^kx^jy_k+...\}=\\
&=2\sum_k(b_{i_1i_2}^k-b_{i_2i_1}^k)y_k+\dots\,\,\,.\end{align*}

Since in the RHS we must have a linear combination of $H_1,...,H_n$ we conclude that
$$\{H_{i_1},H_{i_2}\}=\sum_k g_{i_1i_2}^kH_k,$$
where $g_{i_1i_2}^k=2(b_{i_2 i_1}^k-b_{i_1 i_2}^k)$ are the structure constants of $\g$.
 
Since there are no linear terms in $x^i, 1\le i\le n$ in the expression for the generator $H_k$, we can see (by directly computing the coefficient in front of $x^i$ in the LHS of the above formula) that
the tensor $a_{ijk}$ is symmetric with respect to arbitrary permutations of the indices $i,j,k$.

 In a similar way comparing in the above formula the coefficients in front of the quadratic monomials in $x^i$ and $y_j$   we can reformulate the notion of a classical Airy structure in the following way: it is a collection of tensors
 
 \[\begin{array}{l}A\in (V^{\otimes 3})^{S_3}\subset V^{\otimes 3},\\
 B\in Hom(V, V\otimes V),\\
 C\in Hom((V\otimes V)_{S_2}, V)=Hom(V\otimes V,V)^{S_2}\subset Hom(V\otimes V,V),\end{array}\]
 where the notation like $E^{S_m}$ (resp $E_{S_m}$) refers to invariants (resp. coinvariants)
of the action of the symmetric group $S_m$ on an $S_m$-vector space $E$.    
The above three tensors are required to satisfy the following three relations: for any $1\le i_1,i_2,i_3,i_4\le n$

\begin{equation}\begin{array}{l}
\mbox{Coeff}_{x^{i_3} x^{i_4}} \Delta_{i_1,i_2}=0,\,\,\mbox{Coeff}_{x^{i_3} y_{i_4}} \Delta_{i_1,i_2}=0,\,\,\mbox{Coeff}_{y_{i_3} y_{i_4}} \Delta_{i_1,i_2}=0\\
\mbox{ where } \Delta_{i_1,i_2}:=\{H_{i_1},H_{i_2}\}
- 2\sum_k (b^k_{i_2 i_1}-b^k_{i_1 i_2}) H_k.
\end{array}
\label{classical Airy relations}\end{equation}

These relations can be considered in a coordinate-free way as certain cumbersome purely quadratic  in $A,B,C$ expressions with values in respectively
$$\wedge^2 V \otimes Sym^2 V,\,\,\,Hom(V,\wedge^2 V \otimes V),\,\,\,Hom(Sym^2 V, \wedge^2 V).$$

Furthermore, we can reformulate the notion of classical Airy structure geometrically.
Let $M$ be a $\K$-affine symplectic space of dimension $2n$ and $\g$ is a Lie algebra
over $\K$ of dimension $n$ which acts on $M$ by Hamiltonian transformations. Assume that the corresponding moment map $\mu: M\to \g^{\ast}$ is at most quadratic. 
This means that the moment map is given explicitly by $\mu(p)=(H_1(p),...,H_n(p))$, where each $H_i$ is a polynomial of degree less or equal than two on $M\simeq \K^{2n}$.

Let $p\in \mu^{-1}(0)$ be a smooth point. 
The germ $L$ at $p$ of the algebraic variety $\mu^{-1}(0)$ is a smooth coisotropic subvariety of $M$ of dimension $n$, hence it is Lagrangian.
We will call such Lagrangian submanifolds (or germs of such submanifolds)
{\it quadratic Lagrangian submanifolds}.

The tangent space $T_p L\subset T_pM$ is a Lagrangian subspace. Let us choose a Lagrangian 
vector subspace $V\subset T_pM$ which is transversal to $T_pL$ (we call it {\it Lagrangian complement}). 

\begin{prp} The vector space $V$ carries a natural Airy structure. Conversely any classical Airy structure is isomorphic in the natural sense to the one given by this geometric construction.

\end{prp}
{\it Proof.} Indeed the above choices identify symplectically $M$ with $V\oplus V^{\ast}=T^{\ast}(V^{\ast})$, as long as we identified $V^{\ast}$ with $T_pL$. The quadratic Hamiltonians $H_i$ has the form
$$H_i:=-y_i+\sum_{jk}a_{ijk}x^jx^k+2\sum_{jk}b_{ij}^kx^jy_k+\sum_{jk}c_i^{jk}y_jy_k, 1\le i\le n.$$
Hence we have a classical Airy structure on $V$.

Conversely, if we have a classical Airy structure on $V$, we observe that $\g=V^{\ast}$ acts in the Hamiltonian way on $M=V\oplus V^{\ast}$, the moment map $\mu: V\oplus V^{\ast}\to \g^{\ast}\simeq V$ is given by
the formula $\mu(p)=(H_1(p),...,H_n(p))$ with the Hamiltonians as above, and that $p=0$ is a smooth point of $\mu^{-1}(0)$. This completes the proof. $\blacksquare$

\begin{rmk}
The space $\OO(M)_{\le 2}$ of polynomial functions on $M$ of degree less or equal than $2$ is a Lie subalgebra of the 
algebra $\OO(M)$ of polynomial functions on $M$. Then the above considerations mean that
a classical  Airy structure  is given by a choice of an $n$-dimensional Lie algebra $\g$ over $\K$ together with a monomorphism
of Lie algebras $\phi: \g\to \OO(M)_{\le 2}$, a choice of a point $p\in \mu^{-1}(0)$ where $\mu^{-1}(0)$ is a non-empty smooth variety of dimension $n$ and a choice of a transversal Lagrangian complement to the tangent space to $\mu^{-1}(0)$ at $p$ . 
\end{rmk}

Suppose that the field $\K=\C$ or $\R$. Then we have the Lie group $G=exp(\g)$ (say, connected and simply-connected). This group acts on $\mu^{-1}(0)$ by affine symplectic isomorphisms. On a small neighborhood of the smooth point $p$ the action is locally transitive with the discrete stabilizer. This means that for any two points $p_1,p_2$ which are sufficiently close to $p$  there is a unique  element $g=g_{p_1,p_2}\in G$ close to $id_G$  which maps $p_1$ to $p_2$ and hence inducing an identifications of affine spaces of Lagrangian complements to $T_{p_i}\mu^{-1}(0)\subset T_{p_i} M$,  and an equivalence of the corresponding Airy structures. It is not clear a priori (although very likely) that the obtained equivalence depends {\it algebraically}  (and not only analytically) on the Lie subalgebra  $\mathfrak{g}\subset \mathcal{O}(M)_{\le 2}$, and a pair of smooth points $p_1,p_2\in \mu^{-1}(0)$.

Finally, let us write the formula for the change of tensors $A,B,C$ with respect to the change of the transversal Lagrangian subspace at a given smooth point $p$. Any two Lagrangian complements $V_p^{(1)}, V_p^{(2)}$ to $T_p L\subset T_p M$ are both canonically
 identified with $T_p^* L=T_p M/T_p L$, hence identified with each other. One of them can be seen as the graph of linear map $V\to V^*$ which comes from a symmetric bilinear pairing $Sym^2(V)\to \K$. In local coordinates the calculation reduces to the application of a linear symplectic change of variables
 $$y_i\to y_i, \,\,\,x^i\to x^i+\sum_j s^{ij} y_j,$$
 where $s=(s^{ij})_{1\le i,j\le n}$ is a symmetric matrix.
 A direct calculation gives us the following {\bf gauge transformaton}:
 \begin{equation}\begin{array}{lll}
 a_{ijk} &\to &a_{ijk},\\
 & &\\
 b_{ij}^k&\to &b_{ij}^k+\sum_l a_{ijl} s^{lk},\\
 & &\\
 c_i^{jk}&\to &c_i^{jk}+\sum_{l_1,l_2} a_{il_1 l_2} s^{l_1 j} s^{l_2 k}+\sum_l b_{il}^k s^{lj}+\sum_l b_{il}^j s^{lk}.
 \end{array}\label{classical gauge change}\end{equation}

\subsection{Quantum Airy structures}

We keep the notation of the previous subsection. 
The vector space of at most quadratic functions $Sym_{\le 2}(W), W=V\oplus V^{\ast}$ 
is spanned by 
\[\{1,x^i, y_i, x^i x^j, x^i y_j, y_i y_j\}, \quad1\le i,j\le n\,.
\]
We define an associative $\K$-algebra $Sym^{quant}(W)$ generated by 
$1, \hbar, x^i, \hbar \partial_i$, with $1$ and $\hbar$ being central elements, 
subject to the relations $(\hbar \partial_i)x^j-x^j(\hbar\partial_i)=\hbar \delta_{ij}$ 
(i.e. $\partial_i=\partial/\partial x^i$ if we allow the division by $\hbar$). We define $Sym_{\le 2}^{quant}(W)$ as the vector subspace
of $Sym^{quant}(W)$ spanned by 
\[\{1, \hbar,  x^i, \hbar \partial_i,x^i x^j, \hbar(x^i\partial_j+\partial_j x^i), \hbar^2\partial_i \partial_j\},  \quad1\le i,j\le n\,.\] 
It is a Lie algebra with the bracket equal to $(ab-ba)/\hbar$. Clearly modulo $\hbar$ it is isomorphic to the Lie algebra 
$Sym_{\le 2}(W)$. 

\begin{defn} A {\bf quantization} of a classical Airy structure  is given by a monomorphism
$\g\simeq V^{\ast}\to Sym_{\le 2}^{quant}(W)$ such that modulo $\hbar$ it coincides with  the monomorphism 
$\g\to Sym_{\le 2}(W)  $. 
We will call the former monomorphism quantum Airy structure.

\end{defn}

We have quantized Hamilitonians 
$$\widehat{H}_i=-\hbar\partial_i+\sum_{j,k}a_{ijk}x^jx^k+
2\hbar\sum_{j,k}b_{ij}^kx^j\partial_k+\hbar^2\sum_{j,k}c_i^{jk}\partial_j \partial_k+\hbar\epp_i\cdot 1,$$
where $\epp=(\epp_i\in \K)_{1\le i\le n}$ are some constants. 

In the coordinate-free notation  $\epp$ is an element of $V$. The constraint on tensor $\epp$ follows from the relation saying that
\begin{equation}\forall {i_1,i_2}\,\,\,\,\mbox{Coeff}_{\,\hbar\cdot 1}\left(\frac{1}{\hbar}\left(\widehat{H}_{i_1} \widehat{H}_{i_2}
-\widehat{H}_{i_2} \widehat{H}_{i_1}\right)-2\sum_k (b^k_{i_2 i_1}-b^k_{i_1 i_2}) \widehat{H}_k\right)=0. \label{additional quantum Airy constraint}\end{equation}
In the coordinate-free form, this constraint means vanishing of  certain purely quadratic expression in $A,B,C,\epp$ with values in $\wedge^2 V$.
Together with classical constraints (\ref{classical Airy relations}) we obtain  an explicit algebraic description of quantum Airy structures.

The gauge transformation formula (\ref{classical gauge change})
 extends to the quantum case as 
 \begin{equation}\epp_i\,\,\,\,\,\to \,\,\,\,\,\epp_i +\sum_{jk} a_{ijk} s^{jk}.  \label{quantum gauge change}\end{equation}
 (Apply transformation $\hbar\partial_i\to \hbar \partial_i,\,\,x^i\to x^i+\sum_j s^{ij} \hbar \partial_j$ where $s=(s^{ij})_{1\le i,j \le n}$ is a symmetric matrix).

 For any classical Airy structure on a {\it finite-dimensional} space $V$ there exists a {\it canonical} choice of its quantization\footnote{In the next section we will talk about Airy structures in the infinite-dimensional case where the analogous assertion is false in general.}. Namely, we set
 \[\widehat{H}_i^{can}:=-\hbar\partial_i+\sum_{j,k}a_{ijk}x^jx^k+
\hbar\sum_{j,k}b_{ij}^kx^j\partial_k+\hbar\sum_{j,k}b_{ij}^k \partial_k x^j+\hbar^2\sum_{j,k}c_i^{jk}\partial_j \partial_k.\]
In other words, we set $\epp_i^{can}:=\sum_j b^j_{ij}$. A general solution to (\ref{additional quantum Airy constraint}) is
\begin{equation}\epp=\epp^{can}+\chi,\mbox{ where }\chi\in V=\g^*\mbox{ is any character of }\g\Leftrightarrow \chi\in H^1(\g,\K).\label{finite quantum via characters}\end{equation}

 \subsection{Airy structures in the infinite-dimensional case}
 
 Four tensors $A,B,C, \epp$ satisfy four relations, see (\ref{classical Airy relations}) and (\ref{additional quantum Airy constraint}).

Notice that since all four relations include acyclic tensors, our definitions of Airy structures (both classical and quantum) make sense in the case of any $\K$-linear tensor category.
In particular, they make sense in the case when $dim\,V=\infty$. Notice that in the infinite-dimensional case Lie algebra $\mathfrak{g}$
is naturally endowed with linear topology being dual to an ordinary vector space $V$ (endowed with the discrete topology). In other words, $V$ is a Lie coalgebra.

Let us analyze constraint  (\ref{additional quantum Airy constraint}). Rewrite it as 
\begin{equation}
\forall {i_1,i_2}\,\,\,\,2\sum_{j,k}(c_{i_1}^{jk}a_{i_2jk}-c_{i_2}^{jk}a_{i_1jk})=\sum_k g_{i_1 i_2}^k \epp_k.
\label{rewritten quantum Airy constraint}\end{equation}
where $g_{i_1 i_2}^k:=2(b_{i_2 i_1}^k-b_{i_1 i_2}^k\!)$ are structure constants of $\mathfrak{g}$.
One sees that the  RHS   is the $(i_1,i_2)$-component  of the coboundary of the continuous cochain $(\epp_k)$ on the Lie algebra $\g=V^{\ast}$.

Let us look more closely on the expression in the LHS in (\ref{rewritten quantum Airy constraint}). As the indices $i_1,i_2$ are fixed,
we see that $a_{i_1jk}, a_{i_2jk}\in Sym^2(V)$, and $c_{i_1}^{jk}, c_{i_2}^{jk}\in \widehat{Sym^2}(V^{\ast})$, where 
the hat notation over a  vector space indicates the
natural topological completion.\footnote{We hope that the reader  will not mix it with the hat indicating
quantization, like in the previous section.}

Notice that the Lie algebra $gl(W)$ carries a natural $2$-cocycle $\psi$
(``Japanese cocycle''). 
 Namely, splitting $W=V\oplus V^{\ast}$ allows us to represent any linear  map $S$ as a block $2\times 2$ matrix $(s_{ij})_{1\le i,j\le 2}$.  Then the cocycle $\psi$ can be written as  $\psi(S,T)=Tr(t_{21}s_{12}-s_{21}t_{12})$.
Comparing it with the above expression for the LHS of (\ref{rewritten quantum Airy constraint}), we see that the latter is twice of the pull-back of the Japanese $2$-cocycle to the Lie algebra $\g$ under the composition of natural morphisms
\[\g\to \widehat{Sym}_{\le 2} W\epi \widehat{Sym^2}(W)\simeq sp(W)\mono gl(W).\]

Therefore  existence of a quantization of a classical Airy structure in these terms corresponds to triviality of the pullback of the cohomology class 
$\psi$.


Ansatz (\ref{finite quantum via characters}) shows  that in the finite-dimensional case one can recover the $1$-cochain $\epp$ from the other data by taking $\epp_i=\sum_{j}b_{ij}^j+\chi({H}_j)$, where $\chi: \g\to \K$ is any character.  In the infinite-dimensional case this formula does not make sense since it involves taking the trace of the tensor $B$. In this case the {\it continuous}  $1$-cochain $\epp$ is a new datum for the quantum Airy structure.

In general, if a quantization of a classical Airy structure given by a Lie algebra $\g$ exists, 
all such quantizations form a torsor over the space of characters $Hom(\g,\K)=H^1_{cont}(\mathfrak{g},\K)$.  In particular, if $H^1_{cont}(\mathfrak{g},\K)=H^2_{cont}(\mathfrak{g},\K)=0$ then there exists a unique  quantization of
the classical Airy structure.

Finally, we remark that
the tensors $A,B,C, \epp$ satisfying the above four quadratic relations make $V$ into an algebra
over a certain PROP.
 Obviously, (quantum or classical) Airy structures  form a category (there is a straightforward notion of a morphism). Also, for any family of Airy structures on spaces $(V_i)_{i\in I}$ the direct sum space $\oplus_{i\in I} V_i$ carries a  natural direct sum Airy structure.

\subsection{Cyclic $DQ$-modules}

We start with considerations which hold in a bigger generality than the framework of Airy structures. The corollary of these considerations will be the (non-trivial) fact that the polydifferential forms $\omega_{g,n}(z_1,...,z_n)$ defined via the topological recursion are
symmetric with respect to the variables $z_1,...,z_n$.

Let us first assume that $V$ is a finite-dimensional vector space over a field $\K$ of characteristic zero. For convenience, we introduce the coordinates $x^i, 1\le i\le \dim V$ on $V^{\ast}$ as we discussed previously\footnote{In this section we {\it do not} denote the dimension of $V$ by $n$, the letter $n$ is reserved for the  homogeneity degree of polynomials in $x^\bullet$.}. In order to keep track of the invariant (= coordinate-free) meaning, we use the natural  identification of the algebra of formal power series $\K[[x^\bullet]]$  with $\prod_{n\ge 0} Sym^n(V)$, by associating with every series $f$ the collection of tensors (Taylor coefficients)
$$Taylor_n(f)\in\left( V^{\otimes n}\right)^{Sym_n}\subset V^{\otimes n},\,\,\,Taylor_n(f)_{i_1,\dots,i_n}:=\partial_{i_1}\dots \partial_{i_n}(f)(0).$$

Let ${\bf D}$ be the completed algebra of differential operators on $V^{\ast}$. More precisely, ${\bf D}=\K[[\hbar]]\,[[x^\bullet]]\,[[\hbar\partial_\bullet]]$ endowed with the standard product. Let us introduce a grading on ${\bf D}$ such that $deg\,\hbar =2, deg\, x^i=deg\,\hbar\partial_i=1$. 
Denote by ${\bf I}\subset {\bf D}$ the closed $2$-sided ideal spanned by elements of degree $\ge 2$.

Let now $\widehat{H}_i, 1\le i\le \dim V$ be elements of ${\bf D}$ which satisfy the following properties:

1) $\widehat{H}_i=-\hbar\partial_i+P_i$, where $P_i\in {\bf I}$;

2) for any $i\ne j$ the commutator $[\widehat{H}_i,\widehat{H}_j]$ belongs to the left ideal generated by $\hbar \widehat{H}_\bullet$,
 i.e. there exists (non-unique) ``structure constants'' $g_{i,j}^k\in {\bf D}$ such that
 $$\frac{1}{\hbar} [\widehat{H}_i,\widehat{H}_j]=\sum_k g_{i,j}^k \widehat{H}_k.$$

 Notice that the case when $\widehat{H}_i, 1\le i\le \dim V  $ form a basis of a Lie algebra over $\K$, i.e. $g_{i,j}^k\in \K\subset {\bf D}$ (and, more generally, over $\K[[\hbar]]$: $g_{i,j}^k\in \K[[\hbar]]\subset {\bf D}$) is a particular case of the above condition.
 
Having a formal expression 
$$\psi=\exp(S)=\exp\left(S_0/\hbar+\sum_{g\ge 1}\hbar^{g-1}S_g\right),$$
with 
$$S_g=\sum_{n\ge 0} S_{g,n}\in \prod_{n\ge 0}Sym^n(V)=\K[[x^\bullet]],\forall\,\, g\ge 0\,\,\,$$
 we define a continuous automorphism $T_\psi: {\bf D}\to {\bf D}$ which corresponds informally to the conjugation $D\mapsto \psi^{-1} D \psi$.
Explicitly, the automorphism $T_\psi$ is defined on the generators by the
formulas $T_\psi(x^i)=x^i, T_\psi(\hbar\partial_i)=\hbar\partial_i+\partial_i(S_0+\sum_{g\ge 1}\hbar^{g}S_g)$. Of course this is motivated by the well-known formula $e^{-f} \partial_i e^f=\partial_i+\partial_i(f)$.
By definition we have $T_{\psi\phi}=T_\psi T_\phi$, where the product $\psi\phi$ is defined naturally, as the exponent of the sum of two expressions. The formula $\widehat{H}(e^Sf)=e^S(T_{e^S}(\widehat{H})(f))$ which is valid for the differential operators $\widehat{H}$ and functions $S,f$ motivates the following definition.

\begin{defn} We say that the above formal expression $\psi$ satisfies the system of equations $\widehat{H}_i(\psi)=0$ if $T_\psi(\widehat{H}_i)(1)=0,  1\le i\le \dim V $.

\end{defn}

In the language of quantum mechanics, $\psi$ is a   wave function in WKB form.

\begin{thm} Suppose that the Hamiltonians $\widehat{H}_i, 1\le i\le \dim V$ satisfy the above conditions 1) and  2).
Then there exists a unique  solution to the system of equations $\widehat{H}_i(\psi)=0, 1\le i\le \dim V$
of the form $\psi=exp(S_0/\hbar+\sum_{g\ge 1}\hbar^{g-1}S_g)$, where $S_0=\sum_{n\ge 3} S_{0,n}$ and for $g\ge 1$
$S_g=\sum_{n\ge 1}S_{g,n}$ with  $S_{g,n}\in  Sym^{ n}(V)$ for any $g,n\ge 0$ (in particular $S_{0,\le 2}=S_{\ge 0,0}=0$).

\end{thm}
 {\it Proof.}  We will proceed by induction in $g$, and for each $g$ we will proceed by induction in $n$. We start with $g=0$ case. 
 The equation for $g=0$ says that for any $i$ 
 $$H_i(x^\bullet,\partial_\bullet S_0)=0\in \K[[x^\bullet]].$$
 Here the classical limit $H_i\in \K[[x^\bullet,y_\bullet]]$ is defined as the image of $\widehat{H}_i$
  in ${\bf D}/\hbar{\bf D}\simeq \K[]x^\bullet]]\,[[\hbar \partial_\bullet]]$, when we identify variables $y_i$ with the images of $\hbar\partial_i,\,1\le i \le \dim V$ in the quotient algebra ${\bf D}/\hbar{\bf D}$.
   The condition 1) implies that $H_i=-y_i$ modulo quadratic and higher order terms in $x^\bullet,y_\bullet$.
  The classical limit of the commutation relation 2) implies that $H_\bullet$ generate a Poisson ideal. Therefore, the formal germ at 0 of a submanifold $L$ in $\K^{2\dim V}=V\oplus V^*$
  given by equation $H_i(x^\bullet, y_\bullet)=0$ is coisotropic. It is of dimension $\dim V$, hence it is Lagrangian. This implies that it is given by the graph of $d S_0$ for
  some formal power series $S_0\in \K[[x_\bullet]]$. We can assume that the constant term is zero. The tangent space to $L$ at $0$ is given by equations $y_i=0,\,\,1\le i \le\dim V$,
  hence $S_0$ starts with a cubic term.
  
  It is useful to formulate the proof purely in terms of linear algebra using the induction on $n\ge 3$ (the order of the  term in the Taylor expansion).
  In what follows we will use a shorthand notation where  the upper index describes the constraints on the  homogeneity degree of monomials in variables $x^\bullet$ or $y_\bullet$.
  First, let us prove that there exists a unique cubic polynomial $S_{0,3}=x^{=3}\in \K[x^\bullet]$ such that 
  $H_i(x^\bullet,\partial_\bullet S_{0,3})=x^{\ge 3}\,\,\forall i$.
  Indeed, we have 
  $$H_i=-y_i+H_{i,2}+x^{\ge 3}+x^{\ge 1} y^{=1}+x^{\ge 0} y^{\ge 2},$$
  where $H_{i,2}=x^{=2}$ is a homogeneous polynomial in $x$ of degree $2$. The equation on $S_{0,3}$ says that $-\partial_i S_{0,3}+H_{i,2}=0\,\,\forall i$.
  The Poisson bracket $$\{H_{i_1},H_{i_2}\}=\partial_{i_1}H_{i_2,2}-\partial_{i_2}H_{i_1,2}+x^{\ge 2}+x^{\ge 0} y^{\ge 1}$$
  belongs to the ideal generated by $H_\bullet$.  Hence $-\partial_{i_1}H_{i_2,2}+\partial_{i_2}H_{i_1,2}=0\,\,\, \forall i_1,i_2$, i.e. 1-form
  $\sum_i H_{i,2} dx^i$ is closed, therefore exact and given by $dS_{0,3}$ for a homogeneous cubic polynomial $S_{0,3}$.

  Now, let us assume that we have already  a solution of $H_i(x^\bullet,\partial_\bullet S_0)=0$ up to order $n\ge 3$, and want to extend it to a solution up to order $n+1$.
   Applying the conjugation by $\exp(S_0/\hbar)$ we may assume that the solution up to order $n$ vanishes, i.e.  
   $$S_0=x^{\ge (n+1)},\,\,\,H_i(x^\bullet,\partial_\bullet S_0)=x^{\ge (n+1)}\,\,\forall i\,\,.$$
   This implies that $H_i$ has  the form
   $$H_i=-y_i+H_{i,n}+x^{\ge (n+1)}+ x^{\ge 1} y^{=1}+ x^{\ge 0} y^{\ge 2}.$$
    Similarly to the previous  arguments, the Poisson bracket has the form
   $$\{H_{i_1},H_{i_2}\}=-\partial_{i_1}H_{i_2,n}+\partial_{i_2}H_{i_1,n}+x^{\ge n}+x^{\ge 0} y^{\ge 1}.$$
   Hence homogeneous 1-form $\sum_i  H_{i,n}dx^i$ is closed and can be written as $d S_{0,n+1}$ where $S_{0,n+1}=x^{=(n+1)}$.
   This proves the $n$-th step of induction for the case $g=0$.
   
   Finally, we proceed with the induction on $g$. Assume that for $g\ge 1$ and $n\ge 0$ we have a solution $S$ up to $\hbar^{g-1} x^{\ge (n+1)} + \hbar^{\ge g} x^{\ge 1}$. At beginning of our  induction, when $n=0$, we may assume that $S_{g,0}=0$. Indeed adding a multiple of $\hbar^{g-1}$ to $S$ does not change the constraint $\widehat{H}_i(\exp(S))=0$. Applying the conjugation, we may assume that
   $$S=\hbar^{g-1} x^{\ge (n+1)} + \hbar^{\ge g} x^{\ge 1},\,\,\, \widehat{H}_i (\exp(S))=0 \mbox{ modulo } \hbar^g x^{\ge n}+ \hbar^{\ge (g+1)} x^{\ge 0},$$
   where $\exp(S):=1+S+S^2/2!+\dots $ is a well-defined invertible element of $\K[[\hbar]][[x^\bullet]]\subset {\bf D}$
   (unlike the formal expression $\psi=\exp(S_0/\hbar+\dots)$ when $S_0\ne 0$).
 The latter condition on $\widehat{H}_i, S$  implies that
 $$ \widehat{H}_i=-\hbar\partial_i+\hbar^g x^{\ge n}+\hbar^{\ge (g+1)}x^{\ge 0}+\hbar^{\ge 0} x^{\ge 1}(\hbar\partial)^{\ge 1}
 + \hbar^{\ge 1} x^{\ge 0} (\hbar\partial)^{\ge 1}.$$
 Let us denote by $H_{g;i,n}\in \K[x^\bullet]$ the homogeneous term  of the type $\hbar^g x^{=n}$ in $ \widehat{H}_i$.
 The commutator has the form
 $$[\widehat{H}_{i_1} ,\widehat{H}_{i_2}]=-\hbar \partial_{i_1} H_{g;i_2,n}+\hbar \partial_{i_2} H_{g;i_1,n}+\mbox{ higher order terms},$$
 where the higher order terms are of the form
 $$\hbar^{g+1} x^{\ge n}
 +\hbar^{\ge (g+2)}x^{\ge 0}+\hbar^{\ge 0} x^{\ge 1}(\hbar\partial)^{\ge 1}
 + \hbar^{\ge 1} x^{\ge 0} (\hbar\partial)^{\ge 1}.$$

 The constraint 2)
 $$[\widehat{H}_{i_1} ,\widehat{H}_{i_2}]\in \hbar \sum_j {\bf D}\cdot \widehat{H}_j$$ 
 implies that $\hbar \partial_{i_1} H_{g;i_2,n}=\hbar \partial_{i_2} H_{g;i_1,n}$, hence there exists a unique homegenous
  $S_{g,n}\in \K[x^\bullet]$ of degree $n+1$ such that $H_{g;i,n}=\hbar^g \partial_i S_{g,n}\,\,\forall i$.
  Then $S=\hbar^{g-1} S_{g,n}+\dots$ is a solution up to higher order terms
  $\hbar^{g-1} x^{\ge (n+2)} + \hbar^{\ge g} x^{\ge 1}$. We proved the induction step $(g,n)\to (g,n+1)$, hence proved the theorem. $\blacksquare$

   Let us discuss algebraic meaning of the above theorem.
   First of all, the algebra $\bf D$ is a topologically free module over the central subalgebra $\K[[\hbar]]$. Furthermore $[{\bf D},{\bf D}]\subset \hbar{\bf D}$.
    We claim that the condition 2) on the collection of elements $\widehat{H}_i$ is equivalent to the constraint {\it only } on the left ideal $J:=\sum_i {\bf D}\cdot \widehat{H}_i\subset {\bf D}$:
    $$J={\bf D} J,\,\,\,[J,J]\subset \hbar J.$$
    The implication in one direction is obvious. Conversely, assuming condition 2) we have for  $\forall a_1,a_2\in {\bf D},\,\,\forall i_1,i_2$
    $$[a_1\widehat{H}_{i_1} ,a_2\widehat{H}_{i_2}]=[a_1\widehat{H}_{i_1},a_2]\,\widehat{H}_{i_2}+
   a_2a_1[\widehat{H}_{i_1},\widehat{H}_{i_2}]+ a_2\,[a_1,\widehat{H}_{i_2}]\,\widehat{H}_{i_1} \in \hbar  J.$$
   Theorem means that under conditions 1), 2) the cyclic module ${\cal M}:={\bf D}/J$ is isomorphic to the canonical holonomic ${\bf D}$-module
   associated with the germ of Lagrangian submanifold ${graph }(dS_0)$. In particular it is a topologically free  $\K[[\hbar]]$-module.
   
   In the case (relevant to Airy structures) when $\g:=\oplus_i \K\cdot  \widehat{H}_i$ is closed under Lie bracket,  the fact that cyclic module
   $
   \mathcal{E}:={\bf D}/\sum_i {\bf D}\widehat{H}_i$ 
   is topologically free over $\K[[\hbar]]$, has an easy alternative proof. Indeed, let us consider the $\hbar$-enveloping algebra $U_{\hbar}(\g)$ of $\g$ (i.e. the commutation relations
 are $ab-ba=\hbar[a,b], a,b, \in \g$). It is a deformation quantization of $Sym^{\ast}(\g)$ 
 endowed with the Lie-Kirillov-Kostant Poisson bracket. The natural completion  $\widehat{U}_{\hbar}(\g)$ 
 is a topologically
 free $\K[[\hbar]]$-algebra, which is naturally  embedded to the topologically free $\K[[\hbar]]$-algebra $\K[[x_i,\hbar \partial_i,\hbar]]$
 of $\hbar$-differential operators (it is an embedding modulo $\hbar$ which implies an embedding over $\K[[\hbar]]$). 
 Moreover, $\bf D$ is topologically free as right  $U_{\hbar}(\g)$-module. Finally, 
 $${\mathcal E}={\bf D}/\sum_i {\bf D}\widehat{H}_i={\bf D}\otimes_{U_{\hbar}(\g)} \K[[\hbar]],$$
where we treat  $\K[[\hbar]]$ as the trivial rank one
 representation of $\widehat{U}_{\hbar}(\g)$. Hence $\mathcal E$ is topologically free over $\K[[\hbar]]$.
  This argument generalizes immediately to the case of Lie algebra over $\K[[\hbar]]$.

 \subsection{Application to quantum Airy structures: abstract topological recursion}\label{section:abastract recursion}

 Assume that our Hamiltonians $\widehat{H}_i, 1\le i\le n$ come from a quantum Airy structure.
  The general result from the previous subsection implies that there exists a unique WKB solution $\exp(S_0/\hbar+\dots)$.
  We obtain a collection of symmetric tensors
  $$S_{g,n}\in  Sym^n(V),\,\,g=0,n\ge 3\mbox { or } g\ge 1,n\ge 1,$$
  $$ S_{g,n;i_1,\dots,i_n}:=\partial_{i_1}\dots \partial_{i_n} S_g (0).$$

The equation $\widehat{H}_i(\exp(S))=0$ for $S=\sum_g \hbar^{g-1} S_g,\,\,S_g=\sum_n S_{g,n}$ reads as
$$ \hbar \partial_i S=\sum_{jk} a_{ijk}\,x^j x^k+ 2\sum_{jk}b_{ij}^k\, x^j \,\hbar\partial_k S+$$
$$+ \sum_{jk}c_i^{jk}\,\hbar \partial_j(S) \hbar\partial_k(S)+  \sum_{jk}c_i^{jk}\,\hbar^2 \partial_j\partial_k(S)+\hbar \epp_i.$$

Expanding the RHS of the above equation into series in $\hbar,x^\bullet$ we obtain  recursive formulas for $\partial_i S_{g,n}$ in terms of 
 tensors $A,B,C,\epp$ and previously known $S_{g',n'}$ with either $g'<g$ or $g'=g,n'<n$.
   It is {\it not obvious} a priori that the RHS which should be equal by induction to the 1-form $d S_{g,n}=\sum_i \partial_i S_{g,n} dx^i$ is in fact a closed 1-form (which would imply  exactness).
   Theorem guarantees that this is the case.
   
   The recursion starts with 
  the cases $g=0,n=3$ and $g=1,n=1$. Then
  $$\mbox{at the level of tensors }S_{0,3}=2A\,\,\Longleftrightarrow \mbox{ as a polynomial } S_{0,3}=\frac{1}{3}\sum_{ijk} a_{ijk}x^i x^j x^k;$$
  $$\mbox{at the level of tensors }S_{1,1}=\epp\,\,\Longleftrightarrow \mbox{ as a polynomial } S_{1,1}=\sum_i \epp_i x^i.$$
  
  In general, for all other possible values of $g,n$ (i.e. for
  $g=0,n\ge 4$, or $g=1,n\ge 2$ or $g\ge 2, n\ge 1$) we have
  $$S_{g,n;i,i_1,\dots,i_{n-1}}=2\sum_{\alpha=1}^{n-1}   \sum_k b^k_{i i_\alpha} S_{g,n-1; k, i_{\{1,\dots,,n-1\}\setminus\{\alpha\}}}+$$
  \begin{equation}+\sum_{g_1+g_2=g; J_1\sqcup J_2=\{1,\dots,n-1\}}\sum_{jk} c_i^{jk} S_{g_1,\#J_1+1;j, i_{J_1}} 
  S_{g_2,\#J_2+1; k,i_{J_2}}+\label{Abstract recursion}\end{equation}
  $$+\sum_{jk} c_i^{jk} S_{g-1,n+2;j,k,i_1,\dots,i_{n-1}},$$
  where at the RHS we use only terms $S_{g',n'}$ with either $g'=0,n'\ge 3 $ or $g'\ge1, n'\ge 1$. Also, here for any subset
  $J\subset\{1,\dots,n-1\}$ such that $J=\{\alpha_1,\dots,\alpha_{\# J}\}$ we denote by $i_J$ the sequence $i_{\alpha_1},\dots,
  i_{\alpha_{\# J}}$. This sequence {\it depends} on the order in which we enumerate elements of $J$. Nevertheless this dependence is irrelevant for the  recursion formula because tensors $S_{g,n}$ are symmetric.
 
 We call the above formula {\bf abstract topological recursion}.
  
\begin{exa} Suppose we are given a homogeneous cubic polynomial $P_3$ 
as well as the non-degenerate quadratic form $P_2$, both in variables $y_\bullet$.

Consider the $1/\hbar$-Fourier transform
$$I(\hbar,x)=\int_{Z}e^{-{1\over{\hbar}}xy}e^{{1\over{\hbar}}(P_2(y)+P_3(y))}dy.$$
Here we understand $I(\hbar,x)$ in terms of formal Feynman diagram expansion as $\hbar\to 0$. The ``integration cycle'' $Z$ is the germ at the origin of a real closed $n$-dimensional ball in ${\bf C}^n$  such that $Re(P_2+P_3)_{|Z-\{0\}}<0$.

Then  one has
\[\log(I(\hbar,x))\sim const+ \sum_{g\ge 0,n\ge 0, 2-2g-n<0}\hbar^{g-1}x^nS_{g,n}\] with the coefficients computed by Feynman rules. This means that
 \[S_{g,n}=\sum_{G(g,n)}R_{G(g,n)}/|Aut(G(g,n)|\,,\] where  
the sum is taken over all connected $3$-valent graphs
of genus $g$ with $n$ numbered tails, each tail is assigned the corresponding variable $x^\bullet$. The degree $n$ homogeneous polynomial 
$R_{G(g,n)}(x)$ 
is obtained by  formal application of the Feynman rules (i.e. every edge is assigned the inverse to the quadratic form $P_2$, while each $3$-valent vertex is assigned the cubic form $P_3$). Then our theorem claims that $S_{g,n}$ satisfy the abstract topological recursion formulas.
The corresponding Lie algebra $\g$ is abelian, and the quantum Airy structure is obtained by the gauge transformation (see (\ref{classical gauge change}) and (\ref{quantum gauge change}) )
from the following Ansatz:
\[\widehat{H}_i:=-\hbar\partial_i+\sum_{j,k}a_{ijk}x^jx^k,\,\,\,\widehat{H}_i\left(\exp\left(\hbar^{-1} S_{0,3}(x)\right)\right)=0.\]
Here $S_{0,3}=\frac{1}{3}\sum_{ijk} a_{ijk} x^i x^j x^k $ is obtained by a change of variables in $P_3$ from $y_\bullet$ to $x^\bullet$ via the inverse of the quadratic form $P_2$.

\end{exa}

\subsection{Primitive Airy structures}

The terminology below will be used both in classical and quantum cases. For that reason sometimes we will skip the words ``classical'' or ``quantum'',
and will speak about Airy structures, etc.
 The mathematical content of this section is almost empty, we introduce it only to handle the fact that the topological recursion as originally defined by Eynard and Orantin, does {\it not} satisfy the axiomatics of Airy structures. Nevertheless, the data of the conventional TR contain an Airy substructure. Hence by our general theory of Airy structures the conventional TR produces a collection of symmetric tensors (usually denoted by $\omega_{g,n}$).

\begin{defn} A {\bf classical pre-Airy structure} on vector space $V$ given by  a tuple $(A,B,C)$ such that 
$$A\in V^{\otimes 3}, \,\,B(V)\subset V\otimes V, \,\,C(V\otimes V)\subset V.$$
A {\bf quantum pre-Airy structure} on $V$ is given by  a tuple $(A,B,C,\epp)$ such that $(A,B,C)$ is a classical pre-Airy structure and $\epp\in V$ is a vector.
A {\bf substructure} of a (quantum or classical) pre-Airy structure is given by a vector subspace $U\subset V$ such that 
$$A\in U^{\otimes 3}, \,\,B(U)\subset U\otimes U, \,\,C(U\otimes U)\subset U.$$
and $\epp\in U$ in the quantum case.
\end{defn}
Obviously, any quantum pre-Airy structure on $V$  induces a classical pre-Airy structure.

\begin{defn}  A (quantum or classical)
pre-Airy structure on space $V$ is called {\bf primitive} if it does not contain a nontrivial substructure (i.e. the one with $U\ne V$).
\end{defn}
Obviously, if $U$ is  an Airy substructure, then it carries the induced classical Airy structure. Equivalently, an Airy substructure is the same as a morphism
 of Airy structures which is an inclusion at the level of underlying vector spaces.
It is easy to see that   any (quantum or classical) pre-Airy structure contains a unique primitive  substructure.

In both quantum or classical cases, an Airy structure is a particular case of a pre-Airy structure.
 A substructure (in the sense of pre-structures) of an Airy structure is again an Airy structure.

If a classical Airy structure has an  Airy substructure  then the generating series $S_0$ from the Theorem 2.4.2 is  constant along fibers of the linear surjection $V^*\epi U^*$. Equivalently, for a basis in $V$ such that $U$ is spanned by a part of the basis, the series  $S_0$  depends only on the corresponding subset of the coordinates in $V^*$. Classical Hamiltonians corresponding to elements of 
$(V/U)^*\subset V^*$ are linear (i.e. they do not contain non-trivial quadratic forms).

Geometrically, the existence of a non-trivial classical Airy substructure implies that the germ $L$ of the corresponding Lagrangian submanifold in an affine symplectic space $M$ lies in  a non-trivial affine coisotropic submanifold $N\subset M$. Hence $L$ is the pullback
 of a Lagrangian germ $L'\subset M'$, where $M'$ is a  symplectic affine space which is the quotient of $N$ by the constant  foliation given by the skew-orthogonal $T_N^\perp \subset T_N$ (equivalently, by the kernel of the symplectic form restricted to $N$).

Similarly to the above, one can speak about substructures and primitivity of {\it quantum} Airy structures, by adding the condition $\epp \in U$. Quantum Hamiltionians corresponding to
$(V/U)^*\subset V^*$ are just derivations along corresponding directions, and the equation on wave function means that $\exp(S)$ does not depend on a part of coordinates.
 
It follows from the above considerations that if one uses a quantum Airy structure only as a technical tool for producing the series $S_g,\,g\ge 0$ then it suffices to work with any 
quantum Airy substructure, in particular with the primitive one. 

\begin{rmk} \label{modification pre Airy}The modification formulas \eqref{classical gauge change} and \eqref{quantum gauge change} make sense for pre-Airy structures as well, where the symmetric marrix $s=(s^{ij})$ controlling the gauge change is understood in coordinate-free manner as a symmetric bilinear form on $V$. Obviosuly, if $U\subset V$ is a substructure, then the restriction of bilinear form $s$ to $U$ gives a gauge change of the induced structure on $U$.
\end{rmk}
   
   \section{Comparison with topological recursion}
   
   \subsection{Reminder: topological recursion}\label{section: TR}
   
   Here we recall the classical Eynard-Orantin formalism for topological recursion.
   Let $\Sigma$ be a smooth complex curve (not necessarily connected, e.g. a disjoint union of discs).
   We assume that we are given a finite collection of distinct points $(r_{\alpha}\in \Sigma)_{\alpha\in Ram}$ where $Ram$ is a finite set, a  collection of germs of involutions $\sigma=\sigma_{\alpha}$ defined in the vicinity of $r_\alpha$ such that 
   $$\sigma(r_\alpha)=r_\alpha,\,  \sigma_{|T_{r_\alpha} \Sigma}= -id.$$
   In what follows we will consider ``tensor fields'' on the product of copies of $\Sigma$. They are  sections of tensor products of pull-backs the line bundles $K_\Sigma$ and $K_\Sigma^{ -1}$, where $K_\Sigma=T_\Sigma^{\ast}$ is the canonical bundle of $\Sigma$. For example, the expression  $dz_1\dots dz_n$ on $\C^n$ where $\Sigma=\C$ and  $(z_i)_{i=1,\dots,n}$ are coordinates on $\Sigma^n$  is not the top degree differential form, it is a polydifferential. Almost all our tensor fields will be polydifferentials, i.e. the pull-backs of $K_\Sigma^{-1}$ will be absent (exception is the recursion kernel below). We denote by $(\pi_i)_{i=1,\dots,n}$ the projections $\Sigma^n\to \Sigma$ to factors.
   
   The input data for the topological recursion consist of:
   \begin{enumerate}
   \item The {\it Bergman kernel} $\omega_{0,2}\in \Gamma(\Sigma\times \Sigma\setminus\{ diag\},\pi_1^*K_\Sigma\otimes \pi_2^*K_\Sigma)$
    such that $\omega_{0,2}=\omega_{0,2}(z_1,z_2)$ is symmetric under permutation $z_1\leftrightarrow z_2$ and has second order pole
    at the diagonal with the leading coefficient of the polar part equal to $1$: in local coordinates near every point
    $$\omega_{0,2}(z_1,z_2)=\frac{dz_1\,dz_2}{(z_1-z_2)^2}+\mbox{ regular terms for } z_1 \mbox { close to }z_2.$$
    \item A germ of $1$-form $\omega_{0,1}=\omega_{0,1}(z_1)$ defined in a vicinity of each point $r_\alpha$, satisfying the condition that $$ \omega_{0,1}(z_1)-\omega_{0,1}(\sigma(z_1))\mbox{ vanishes with the order $2$ at }r_\alpha.$$
In particular the form $\omega_{0,1}$ vanishes at $r_\alpha$.
   \end{enumerate}
   In a typical example, $\Sigma$ is an affine curve in the plane $\C\times \C$ endowed with coordinates $(x,y)$, such that the projection to $x$-coordinate $\C$ has only double ramification points $r_\alpha$ (hence the notation $Ram$ for the set of ramification points). The  local involution $\sigma_\alpha$ interchanges two branches of the projection $x:\Sigma\to \C$ near $r_\alpha$, and (the germ of) the 1-form $\omega_{0,1}$ is $y\,dx$.
   
   \begin{defn} The {\bf recursion kernel} $K=K(z_1,z_2=z)$ is a collection labeled by $\alpha \in Ram$ of sections of $\pi_1^*K_\Sigma\otimes \pi_2^*K_\Sigma^{-1}$ defined 
    for arbitrary
   $z_1\in \Sigma$ such that $z_1\ne  r_\alpha$ for \emph{all} $\alpha\in Ram$, and for arbitrary  $z$ lying in the vicinity of $r_\alpha$ for {\emph{some}} $\alpha$, $z\ne z_1,r_\alpha$, given by
   the formula 
   $$K(z_1,z):=\frac{-1}{2}\cdot \frac{\int_{z'=\sigma(z)}^{z'=z} \omega_{0,2}(z_1,z')}{\omega_{0,1}(z)-\omega_{0,1}(\sigma(z))},$$
   where in the numerator we use arbitrary integration path from $z'=\sigma(z)$ to $z'=z$ in the vicinity of  $r_\alpha$.
  \end{defn}
  
\begin{exa} The most important for us example is the {\it Airy curve} $\Sigma=\{(x,y)\in \C^2|\, x=y^2\}$ parametrized by the coordinate $z$ so that $x=z^2,y=z$.
In this case we have only one ramification point $z=0$. The involution is given by $\sigma(z)=-z$.
The most natural choice of the Bergman kernel is $\omega_{0,2}(z_1,z_2):=\frac {dz_1\, dz_2}{(z_2-z_1)^2}$.
Also, the 1-form $\omega_{0,1}$ is $\omega_{0,1}(z):=y(z) dx(z)=2z^2 dz$.

The integral  in the definition of the recursion kernel is
 $$\int_{z'=-z}^{z'=z} \frac{dz'\,dz_1}{(z'-z_1)^2}=dz_1\cdot \int_{z'=-z}^{z'=z} d\frac{1}{z_1-z'}=dz_1\cdot\left( \frac{1}{z_1-z}-\frac{1}{z_1+z}\right)=\frac{2z\, dz_1}{(z_1^2-z^2)}$$
$$\Longrightarrow K(z_1,z):=\frac{1}{4(z^2-z_1^2) z} \frac{dz_1}{dz}.$$  \label{Airy curve example}
\end{exa}

Finally, given the input data $\omega_{0,1},\omega_{0,2}$ one defines a collection of new tensor fields $\omega_{g,n}$
 for $g=0,n\ge 3$ or $g\ge 1,n\ge 1$ which are sections of $\pi_1^* K_\Sigma\otimes\dots \otimes \pi_n^* K_\Sigma$ on $(\Sigma-\{r_\alpha|\,\alpha\in Ram\})^n$ inductively by the formula
 $$\mbox{for }g=0,n\ge 3\mbox{ or }g\ge1, n\ge 1  \,\,\,\,\,\,\,\,\,\,\,\,\,\,\,\,\,  \omega_{g,n}(z_1,z_2,\dots,z_n)=$$
 \begin{equation}\sum_\alpha Res_{z=r_\alpha} K(z_1,z)\cdot
  \omega_{g-1,n+1} (z,\sigma(z),z_2,\dots,z_n)+\label{core topological reduction}\end{equation}
  $$+\sum_{\alpha}\sideset{}{^*}\sum_{\substack {g_1+g_2=g \\ I_1\sqcup I_2=\{2,\dots,n\} }} Res_{z=r_\alpha} K(z_1,z)\cdot \omega_{g_1,1+\#I_1}(z,z_{I_1})\cdot \omega_{g_2,1+\#I_2}(\sigma(z),z_{I_2}).$$
The first term in the RHS appears only if $g\ge 1$. The symbol $\sum^*$  in the second term indicates that we {\it do not} allow cases $g_1=0,\#I_1=0$ or $g_2=0,\#I_2=0$, i.e. only $\omega_{0,2}$, $\omega_{0,\ge 3}$ and $\omega_{\ge1, \ge 1}$ can appear in the RHS.
Notice that $\omega_{0,1}$ does not appear explicitly  in the RHS, it participates in the whole formalism only via the denominator $\omega_{0,1}(z)-\omega_{0,1}(\sigma(z))$ for the recursion kernel.

  An amazing property of the above formalism is that the tensor fields $\omega_{g,n}$ defined recursively are invariant under the action of the symmetric group $S_n$.
  The recursion formulas for cases $g=0,n=3$ and $g=1,n=1$ are somewhat different from the rest:
  \begin{equation}\label{A-formula general}\mbox{\bf A-type }:\omega_{0,3}(z_1,z_2,z_3)=2 \sum_{\alpha}Res _{z= r_\alpha} K(z_1,z)\, \omega_{0,2}(z,z_2)\,\omega_{0,2}(\sigma(z),z_3),\end{equation}
  \begin{equation}\label{epsilon-formula general}\mbox{\bf ${\epp}$-type }:\omega_{1,1}(z_1)=\sum_\alpha Res _{z= r_\alpha} K(z_1,z) \,\omega_{0,2}(z,\sigma(z)).\end{equation}
In further cases $g=0,n\ge 4$, or $g=1,n\ge 2$ or $g=2,n\ge 1$ in the RHS of the
 recursion formula for $\omega_{g,n}$ one has {\it two}
 types of summands: 
 \begin{enumerate}
 \item  {\bf B-type}: the terms including decomposition $g=g_1+g_2,\{2,\dots,n\}=I_1\sqcup I_2$ with {\it either} $g_1=0,\#I_1=1$ or $g_2=0,\#I_2=1$. These terms are given by the action of a linear operator  whose argument $f(z)\,dz$ is a collection of germs of meromorphic $1$-forms at points $r_\alpha$,
  \begin{equation}\label{B-formula general}f(z)dz\mapsto \sum_\alpha Res _{z= r_\alpha}K(z_1,z)\,\omega_{0,2} (\sigma(z),z_2)\, f(z) dz. 
  \end{equation}
  \item {\bf C-type}: all other terms, corresponding to a bilinear operator on pairs $(f(z_1)\,dz_1,g(z_2) \,dz_2)$ of collections of germs of meromorphic forms as above,
  \begin{equation}\label{C-formula general}f(z_1)dz_1\otimes g(z_2) dz_2\mapsto \sum_\alpha Res _{z= r_\alpha}K(z_1,z)\,f(z) dz \,
  g(\sigma(z))\,d(\sigma(z))\,.
  \end{equation}
  \end{enumerate}
  We can recast the topological recursion by saying  that $\omega_{0,3},\omega_{1,1}$ are the {\it initial data} (of correspondingly A-type and $\epp$-type), and the recursive steps use
   only operations of B-type or C-type. Comparing formula \eqref{core topological reduction} with the formula \eqref{Abstract recursion} from Section \ref{section:abastract recursion}, we see that the perfect match: 
$$\omega_{g,n}\longleftrightarrow S_{g,n},\,\,\,\,g=0,n\ge 3\mbox{ or }g\ge 1,n\ge 0\,.$$ Nevertheless, in order to identify the topological recursion with the formalism of Airy structures, one has to choose an appropriate space $V$
    of 1-forms on $\Sigma-
  \{r_\alpha|\,\alpha \in Ram\}$.

 \subsection{Topological recursion, the  case of Airy curve}\label{TR Airy case}
 
 Recall from Example~\ref{Airy curve example} Bergman and recursion kernels for the Airy curve:
 $$\omega_{0,2}(z_1,z_2):=\frac {dz_1\, dz_2}{(z_2-z_1)^2},\,\,\,\,K(z_1,z):=\frac{1}{4(z^2-z_1^2) z} \frac{dz_1}{dz}.$$
 
 Let $V:=z^{-1} \C[z^{-1}] \frac{dz}{z}$. We endow it with the basis
 $$e_n:= z^{-n} \frac{dz}{z},\,\,\,n=1,2,\dots$$
 Define 4 tensors of topological recursion 
$${A}_{TR}\in Sym^3(V),{B}_{TR}\in Hom(V,V\otimes V),{C}_{TR}\in Hom(V\otimes V,V)^{S_2},{\epp}_{TR}\in V$$
via formulas (see \eqref{A-formula general},\eqref{epsilon-formula general},\eqref{B-formula general}, \eqref{C-formula general}):
 \[\begin{array}{ll}
 {A}_{TR}:=   & Res_{z} K(z_1,z) \,\omega_{0,2}(z,z_2)\, \omega_{0,2}(-z,z_3),    \\
 {B}_{TR}: f(z) dz\mapsto   & Res_{z}  K(z_1,z)\, \omega_{0,2}( z,z_2)\, f(-z) d(-z) ,                        \\
 {C}_{TR}: f(z)dz\otimes g(z)dz\mapsto   &Res_{z}  K(z_1,z)\, f(z) dz\, g(-z)d(-z),              \\
{\epp}_{TR}:= & Res_{z}  K(z_1,z)\, \omega_{0,2}(z,-z),
\end{array}\]
where $Res_z:=Res_{z=0}$.
Let us calculate the structure constants of these tensors in the basis $(e_n)_{n\ge 1}$. We have:
$${A}_{TR}=-\frac{1}{4}Res_{z}  \frac {1}{(z_1^2-z^2) z}\frac {dz_1}{dz}\frac{dz\,dz_2}{(z_2-z)^2}\frac{-dz \,dz_3}{(z_3+z)^2}=$$
$$=\frac{dz_1\,dz_2\,dz_3}{4} \,Res_{z}  \frac{dz}{z}\left(\frac{1}{z_1^2}+\frac{z^2}{z_1^4}+\dots\right)\left(\frac{1}{z_2^2}+\frac{2z}{z_2^3}+\dots\right)\left(\frac{1}{z_3^2}-\frac{2z}{z_3^3}+\dots\right)=$$
$$=\frac{1}{4}\frac{dz_1}{z_1^2}\frac{dz_2}{z_2^2}
\frac{dz_3}{z_3^2}=\frac{1}{4}\,e_1\otimes e_1\otimes e_1,$$

$${B}_{TR}: e_n=z^{-n}\frac{dz}{z}\mapsto-\frac{1}{4} Res_{z}  \frac{1}{(z_1^2-z^2)z}\frac{dz_1}{dz}\frac{dz\,dz_2}{(z_2-z)^2}(-z)^{-n}\frac{dz}{z}=$$
$$=-\frac{1}{4}dz_1\,dz_2\,Res_{z}  \frac{dz}{z}(-1)^n z^{-n-1} \left(\frac{1}{z_1^2}+\frac{z^2}{z_1^4}+\dots\right)\left(\frac{1}{z_2^2}+\frac{2z}{z_2^3}+\frac{3z^2}{z_2^4}+\dots\right)=$$
$$= -\frac{1}{4} dz_1\,dz_2\,(-1)^n \sum_{n_1,n_2\ge 0; \,2n_1+n_2=n+1}  (n_2+1)\, z_1^{-2n_1-2}\, z_2^{-n_2-2}=$$
$$=-\frac{1}{4} \sum_{n_1,n_2\ge 0; \,2n_1+n_2=n+1} (-1)^{n_2+1} (n_2+1) \, e_{2n_1+1}\otimes e_{n_2+1},$$

$${C}_{TR}: e_{n_1}\otimes e_{n_2}=\frac{dz}{z^{n_1+1}}\otimes \frac{dz}{z^{n_2+1}}\mapsto -\frac{1}{4}Res_{z} 
\frac{1}{(z_1^2-z^2)z}\frac{dz_1}{dz}\frac{ dz\,d(-z)}{ z^{n_1+1}\,(-z)^{n_2+1}}=$$
$$=-\frac{1}{4}(-1)^{n_2} \,dz_1\,Res_{z}  \frac{dz}{z} z^{-(n_1+n_2+2)}\left(\frac{1}{z_1^2}+\frac{z^2}{z_1^4}+\dots\right)=$$
$$=\left\{ \begin{array} {ll} 0 & \mbox{ if }  n_1+n_2\in 1+2\Z\\
                      -\frac{1}{4}              (-1)^{n_2}\,e_{2m+1} & \mbox{ if }  n_1+n_2+2=2m\ge 2,
\end{array}\right.$$

$${\epp}_{TR}=  -\frac{1}{4}Res_{z}  \frac{1}{(z_1^2-z^2)z}\frac{dz_1}{dz} \frac{ dz\,d(-z)}{ 4z^2}= -\frac{1}{4} dz_1\, Res_{z}  \frac{dz}{z}\frac{-1}{4z^2}\left( \frac{1}{z_1^2}+\frac{z^2}{z_1^4}+\dots\right)=$$
$$= \frac{dz_1}{16\,z_1^4}=\frac{1}{16}e_3.$$

For the future reference we combine above formulas:
\begin{align}
     & A_{TR}=\frac{1}{4}e_1\otimes e_1\otimes e_1,\label{A Formulas for TR for Airy curve}\\
&B_{TR}:e_n\mapsto \frac{1}{4}\sum_{n_1,n_2\ge 1;\,n_1+n_2=n+3;\, n_1\,{odd}}\,(-1)^{n_2+1}n_2\,e_{n_1}\otimes e_{n_2},\label{B Formulas for TR for Airy curve}\\
&C_{TR}:e_{n_1}\otimes e_{n_2}\mapsto \begin{cases} 0 & \mbox{ if }  n_1+n_2\mbox{ odd}\\
                      \frac{1}{4}              (-1)^{n_2+1}\,e_{n_1+n_2+3} & \mbox{ if }  n_1+n_2\mbox{ even},\end{cases}  \label{C Formulas for TR for Airy curve} \\           
       &           \epp_{TR}=\frac{1}{16}e_3.\label{epsilon Formulas for TR for Airy curve}
        \end{align}          

We will see below that the data $(V,{A}_{TR},{B}_{TR},{C}_{TR},{\epp}_{TR})$ define a quantum pre-Airy structure on $V$, which is {\it not} a quantum Airy structure. In order to fix the problem let us consider
$V_{odd}:=z^{-1}\C[z^{-2}]\,\frac{dz}{z}\subset V$, which is the vector subspace spanned by $e_1,e_3,e_5,\dots$. We will see that the restrictions of the tensors ${A}_{TR},{B}_{TR},{C}_{TR},{\epp}_{TR}$ to $V_{odd}$ give rise to a quantum Airy structure.

\subsection{Airy structure from residue constraints}

\label{sec: canonical Airy structure}
Let us consider the infinite-dimensional symplectic vector space
$$W_{Airy}:=\{\eta \in \C(\!(z)\!)\,dz\,|\,Res_{z}\, (\eta)=0\}=
\left\{J\cdot \frac{dz}{z}\,\left |\,J= \sum_{n\in \Z;n\ne 0;n\ll+\infty
} J_n z^{-n} \right.\right\}$$
with the symplectic pairing given by
\begin{equation}\label{symp pairing}\langle df_1, df_2\rangle:=Res_{z}\, (f_1 \,df_2)\Longleftrightarrow\left\langle  z^{-n}\frac{dz}{z}, z^{n}\frac{dz}{z}\right\rangle=\frac{1}{n}\Longleftrightarrow
\{J_n,J_{-n}\}=n.\end{equation}
Let us introduce a family of residue  constraints on the elements $\eta$ which are  polynomials of  degree $\le 2$ in the coordinates\footnote{In formulas below  $x^n$ means the power of the variable $x$. The reader should not confuse it with the notation for the coordinates $(x^1,...,x^n)$ in Section 2.} $x,y$ defined below\footnote{The geometric meaning of the coordinates and the constraints will be explained below in Section \ref{meaning res}.
At the moment  the reader should consider all that as a ``black box'' producing an Airy structure.}:
$$\forall n\ge 1 \,\,Res_{z}\,(y \,x^n \,dx)=0,\,\,\,\,\,\,\,\,\,\,\,\forall m\ge 0\,\, Res_{z}\,(y^2\,x^m\,dx)=0,$$
where 
$$x:=z^2\in \C(\!(z)\!),\,\,\,y=y(\eta):=z-\frac{J}{2x}=z-\frac{J}{2z^2},\,\,\,J=\frac{\eta}{z^{-1}dz}.$$
Equivalently, $y(\eta) \,dx= z d(z^2)-\eta$.

The first constraint $Res_{z}\,(y \,x^n \,dx)=0$ gives
$$Res_{z} (\,\frac{J}{z^2}\, z^{2n}\,z\, dz)=Res_{z}(\, \frac{dz}{z}\, J \, z^{2n})=0
\Longleftrightarrow \,J_{2n}=0.$$
The second constraint $Res_{z}\,(y^2 \,x^m \,dx)=0$  for $m\ge 0$ gives
$$Res_{z}\, \left( z-\frac{J}{2z^2}\right)^2 z^{2m}\,z\,dz= Res_{z} \,\frac{dz}{z}\left(- J\,z^{2m+1}+\frac{1}{4} J^2  z^{2m-2}\right)=0\Longleftrightarrow$$
$$-J_{2m+1}+\frac{1}{4}\sum_{k+l=2m-2;k,l\ne 0} J_k J_l=0.$$
Therefore, the first and the second residue constraints can be written together  as a system of equations 
$(H_n=0)_{n\ge 1}$ where:
 $$\begin{array} {l} H_1:= -J_1+\frac{1}{4}J_{-1}^2+\frac{1}{2}\left(J_{-3} J_1+J_{-4} J_2+\dots\right),
 \\H_2:=-J_2,\\
 H_3:=-J_3+\frac{1}{2}\left(J_{-1} J_1+J_{-2} J_2+\dots\right),\\
 H_4:=-J_4,\\
 H_5:=-J_5+\frac{1}{4}J_1^2+\frac{1}{2}\left(J_{-1}J_3 +J_{-2}J_4+\dots\right),\\
 H_6:=-J_6,\\
 H_7:=-J_7+\frac{1}{4}J_2^2+\frac{1}{2}J_1 J_3+\frac{1}{2}\left(J_{-1} J_5 +J_{-2} J_6+\dots\right),\\
 \dots
 \end{array}$$
 
\begin{rmk}\label{Tate}
 The infinite-dimensional topological vector or affine spaces considered in this paper are Tate spaces, i. e. they are direct sum of vector or affine spaces endowed with discrete topology and dual to such spaces. One can also say that Tate spaces are vector spaces endowed with the topology of the vector space of Laurent series. Formal germ at $0\in W$ of the Tate space $W$ (e.g. $W=W_{Airy}$) is understood as a functor on the category of finite-dimensional commutative local Artin algebras:
$$R\mapsto {m}_R^{\ast}\otimes W,$$
where ${m}_R\subset R$ is the maximal ideal.
In this way one can also define the notion of subvariety, intersection of subvarieties, etc.
\end{rmk}

 \begin{lmm} The complete topological vector subspace of the  space $Sym_{\le 2}(W)$ spanned by $(H_n)_{n\ge 1}$ is closed under the Poisson bracket given in coordinates $(J_k)_{k\ne 0}$ by $\{J_k,J_{-k}\}=k$. Moreover, there exists a unique lift of $\prod_{n\ge 1} \C\cdot H_n$ to a Lie algebra of $\hbar$-differential operators.
\end{lmm}
{\it Proof.} It is convenient to add a central variable
$J_0$ and 
 denote for $\forall n\in \Z$
$$L_n:=\frac{1}{4}\sum_{k+l=2n}J_k J_l.$$
In the end we will specialize results to the hypersurface $J_0=0$.
Functions $L_\bullet$ form a Lie algebra with respect to the Poisson bracket:
$$\{L_{n_1},L_{n_2}\}=\frac{1}{16}\sum_{k_1+l_1=2n_1,k_2+l_2=2n_2}
\{J_{k_1}J_{l_1},J_{k_2} J_{l_2}\}=$$
$$=\frac{1}{4}\sum_{k_1+l_1=2n_1,-k_1+l_2=2n_2} k_1 J_{l_1}J_{l_2}=$$
$$=\frac{1}{4}\sum_k k J_{2n_1-k} J_{2n_2+k}=\frac{1}{8}\sum_{k_1+k_2=2(n_1+n_2)}
(( 2n_1-k_1)-(2n_2-k_2))J_{k_1}J_{k_2}=$$
$$=\frac{1}{8}\sum_{k_1+k_2=2(n_1+n_2)} 2(n_1-n_2) J_{k_1}J_{k_2}=
(n_1-n_2) L_{n_1+n_2}.$$
The reader can recognize here the structure constants of the Lie algebra $\C[[x]]\partial_x$ of formal vector fields on the line.
Also, let us calculate the bracket 
$$\{L_n,J_m\}=\frac{1}{4} \sum_{k+l=2n}\{J_k J_l,J_m\}=\frac{1}{2}(-m)J_{2n+m}.$$
We have 
$$H_{2n+3}=-J_{2n+3}+L_n\,\forall n\ge -1,\,\,\,\,\,\,H_{2m}=-J_{2m}\,\forall m\ge 1.$$
This formulas implies
\begin{equation}\begin{array}{l}\left\{H_{2n_1+3},H_{2n_2+3}\right\}=(n_1-n_2)H_{2n_1+2n_2+3}, \\
\left\{H_{2n+3},H_{2m}\right\}=(-m)H_{2n+2m} ,\\
\{H_{2m_1},H_{2m_2}\}=0.
\end{array}
\label{Lie structure constants}
\end{equation}
The above calculation shows that the Lie algebra $\mathfrak g=\prod_{n\ge 1}\C\cdot H_n$ is isomorphic to the quotient of the algebra of differential operators of order $\le 1$ on $\C[[x]]$ by constants $\C\cdot 1$:
$$H_{2n+3}\leftrightarrow  -x^{n+1}\partial_x\,\,\forall n\ge -1,\,\,\,\,\,\,\,\,\,\,H_{2m}\leftrightarrow x^m\,\,\,\,\forall m\ge 1.$$
Let us introduce coordinates $x_1,x_2,\dots$ on an infinite-dimensional space, and realize quantum operators
$(\widehat{J}_n)_{n\in \Z,n\ne 0}$ as
$$\widehat{J}_n=\hbar \frac{\partial}{\partial x_n},\,\,\,\widehat{J}_{-n}:=n x_n\,\,\,\forall n\ge 1.$$
The direct calculation shows that the same commutation relations as in Equation (\ref{Lie structure constants}) hold for
$$\begin{array}{ll}\widehat{H}_{2n+3} & :=-\widehat{J}_{2n+3}+\frac{1}{4}\sum_{k+l=2n;k,l\ne 0} \widehat{J}_k\widehat{J}_l,\,\,\,n=-1,1,2,3,\dots, \\
\widehat{H}_3 &:=-\widehat{J_3}+\frac{1}{2}(\widehat{J}_{-1}\widehat{J}_1+\widehat{J}_{-2}\widehat{J}_2+\dots)+\frac{\hbar}{16},\\
\widehat{H}_{2m}&:=-\widehat{J}_{2m},\,\,\,m\ge 1.
\end{array}$$
More precisely, instead of the Poisson bracket we use the Lie bracket $(ab-ba)/\hbar$.

The uniqueness of the quantum lift follows from the vanishing $H^1_{cont}(\mathfrak{g},\C)$. The existence was proven above via explicit formulas.\footnote{Notice that  $H^2_{cont}(\mathfrak{g},\C)\ne 0$
as $\mathfrak{g}$ has a non-trivial central extension which is the Atiyah Lie algebra $\C[[x]]\partial_x\oplus  \C[[ x]]$, hence the existence of a lift is not automatic.}
 $\blacksquare$

 Let us split $W=W_{Airy}$ into the direct sum of two complementary Lagrangian subspaces:
 $$W=V^*\oplus V,\,\,\,V^*:=\C[[z]] dz, \,\,V:=z^{-2}\C[z^{-1}] dz.$$
 This splitting gives us a quantum Airy structure on $V$.
More precisely, the germ of the quadratic Lagrangian submanifold $L=L_{Airy}$ defined by the system of equations $H_n=0, n\ge 1$ has the tangent space at $0$ given by the system of linear equations $J_1=J_2=....=0$, which is isomorphic to $V^{\ast}$.
 
From the above splitting we derive the tensors $A,B,C,\epp$. E.g. the tensor $B$ is encoded in the bracketed sum in the formulas for $(H_n)_{n\ge 1}$.

In the basis $(e_n)_{n\ge 1}$ we have
\begin{equation}\begin{array}{l}
A=\frac{1}{4}e_1\otimes e_1\otimes e_1,\\
\\
B:e_n\mapsto \frac{1}{4}\sum_{n_1,n_2\ge 1;\,n_1+n_2=n+3;\, n_1\,{odd}}\,n_2\,e_{n_1}\otimes e_{n_2},\\
\\
C:e_{n_1}\otimes e_{n_2}\mapsto \left\{ \begin{array} {ll} 0 & \mbox{ if }  n_1+n_2\mbox{ odd}\\
                      \frac{1}{4}              \,e_{n_1+n_2+3} & \mbox{ if }  n_1+n_2\mbox{ even},\end{array}\right.\\
                      \\
                      \epp=\frac{1}{16}e_3.
\end{array}\label{Formulas for structure for Airy curve}\end{equation}
One can write formulas for tensors $B,C$ using residues:
\begin{equation}
\begin{array}{ll}
 
 {B}: f(z) dz\mapsto   & -Res_{z} \, K(z_1,z)\, \omega_{0,2}( z,z_2)\, f(z) dz,                         \\
\\
 {C}: f(z)dz\otimes g(z)dz\mapsto   &-Res_{z}  \,K(z_1,z)\, f(z) dz\, g(z)dz.           

\end{array}\label{modified TR}\end{equation}

\subsection{Common Airy substructure and modified topological recursion}\label{sec: modified}

 We see (formulas \eqref{A Formulas for TR for Airy curve}, \eqref{B Formulas for TR for Airy curve}, \eqref{C Formulas for TR for Airy curve}, \eqref{epsilon Formulas for TR for Airy curve} versus \eqref{Formulas for structure for Airy curve}) that in the example of Airy curve we have {\it two} collections of tensors $(A_{TR},B_{TR},C_{TR},\epp_{TR})$ and $(A,B,C,\epp)$
 which can be used for the recursive definition of $\omega_{g,n}$ with $g=0,n\ge 3$ or $g\ge1,n\ge 1$.
 
 It turns out  that they produce the {\it same} tensors $\omega_{g,n}$.
It will follow from the fact that both collections of tensors (i.e. the corresponding pre-Airy structures) define  quantum Airy structures on the vector subspace
 $$V_{odd}:=\oplus_{m\ge 0} \,\C\cdot e_{2m+1}\subset V,$$
 and moreover these quantum structures coincide, since the restrictions of the above tensors on $V_{odd}$ coincide.
 
Furthermore $V_{odd}$ endowed with this quantum Airy structure is a  primitive Airy structure (and the corresponding classical Airy structure is also primitive).

 For this primitive {\it classical} structure the existence and uniqueness of the quantum lift is obvious, as $H^1_{cont}(\mathfrak{g}_0,\C)=H^2_{cont}(\mathfrak{g}_0,\C)=0$ for the Lie algebra $\mathfrak{g}_0\simeq \C[[x]]\partial_x$ of formal vector fields on the line (sometimes it is called Witt Lie algebra).

 \begin{rmk}
 One can characterize $V_{odd}$ as a subspace of $V$ intrinsically. Namely, it is the space of such meromorphic 1-forms from $V$ that their direct image
 to the $x$-line under the projection $z\mapsto x=z^2$ has no poles. 
\end{rmk}
  
  Subspaces $V,V_{odd}\subset \C(\!(z)\!)\, dz$ are naturally associated with kernels $\omega_{0,2}(z_1,z)$ and $K(z_1,z)$ in the following way.
   Let us expand these kernels in the iterated Laurent series in the domain $0\ll |z| \ll|z_1|\ll 1$:
 \[ \omega_{0,2}(z_1,z)_{|z|\ll |z_1\!|}  = \left(\frac{1}{z_1^2}+\frac{2z}{z_1^3}+\frac{3z^2}{z_1^4}+\dots\right)dz_1 dz \in \C(\!(z_1)\!)(\!(z)\!)\,{dz_1}{dz},\]
 \[K(z_1,z)_{|z|\ll |z_1\!|} =  \frac{-1}{4}\cdot\left(\frac{z^{-1}}{z_1^2}+\frac{z}{z_1^4}+\frac{z^3}{z_1^6}+\dots\right)\frac{dz_1}{dz}\in 
  \C(\!(z_1)\!)(\!(z)\!)\,\frac{dz_1}{dz}.\]
 We see that $V^{z\to z_1},V^{z\to z_1}_{odd}$ are the {\it minimal} subspaces of $\C(\!(z_1)\!)\,dz_1$ such that
 $$\omega_{0,2}(z_1,z)_{|z|\ll |z_1\!|}\in V(\!(z)\!)\,dz,\,\,\,\,\, K(z_1,z)_{|z|\ll |z_1\!|}\in V_{odd}(\!(z)\!)(dz)^{-1}.$$
 Here the notation $z\to z_1$ means that we replace variable $z$ by $z_1$ and realize $V,V_{odd}$ as subspaces of  $\C(\!(z_1)\!)\,dz_1$ 
 instead of  $\C(\!(z)\!)\,dz$.

Let us define  the recursion based on tensors $A,B,C,\epp$ and using the residue formulas from equation (\ref{modified TR}). We call it the {\bf modified
topological recursion}. Explicitly, it is given by:
$$\text{initial data: } \omega_{0,3}=2A=\frac{1}{2}\frac{dz_1\,dz_2\,dz_3}{z_1^2z_2^2z_3^2},\,\omega_{1,1}=\epsilon=\frac{1}{16}\frac{dz_1}{z_1^4},$$
$$\mbox{for }g=0,n\ge 4\mbox{ or }g=1, n\ge 2\mbox{ or }g\ge 2, n\ge 1 \,\,\, \,\,\,\,\,\,\,\,\,\,\,\ $$

$$\omega_{g,n}(z_1,z_2,\dots,z_n)=$$
 $$-\sum_\alpha Res_{z= r_\alpha} K(z_1,z)\cdot
  \omega_{g-1,n+1} (z,z,z_2,\dots,z_n)+$$ 
  $$-\sum_{\alpha}\sideset{}{^*}\sum_{\substack {g_1+g_2=g \\ I_1\sqcup I_2=\{2,\dots,n\} }}
  Res_{z=r_\alpha} K(z_1,z)\cdot \omega_{g_1,1+\#I_1}(z,z_{I_1})\cdot \omega_{g_2,1+\#I_2}(z,z_{I_2}).$$
The first sum in the RHS appears only if $g\ge 1$. Also, similarly to the Eynard-Orantin topological recursion, the notation $\sum^*$ in  the second sum means that we  do not allow cases $g_1=0,\#I_1=0$ or $g_2=0,\#I_2=0$. In fact, the above formula holds also in the a priori forbidden case $g=0,n=3$ (but not for $g=1,n=1$ as $\omega_{0,2}(z,z)$ is ill-defined).

Let us summarize  differences between the modified topological recursion and the usual one:

 a)  the overall minus sign;
 
 b) the absence of the local involution in the RHS
 of the formula for the modified topological recursion;
 
 c) lesser range of validity of the modified topological recursion (i.e. $\omega_{0,3}$ and $\omega_{1,1}$ are treated in the modified recursion as  initial data). 

\begin{rmk}
An interesting (and a priori unexpected) corollary of the above discussion is the fact that the tensors $\omega_{g,n}$ derived from the usual (i.e. Eynard-Orantin) topological recursion satisfy also the modified topological recursion. This observation seems to be new.
\end{rmk}
 
\subsection{Comparison: TR versus Airy structures in general}

In general,   in the formalism of topological recursion, one can replace the  curve $\Sigma$ by a disjoint union of small discs containing ramification points $\{r_\alpha\}$.
Moreover, one can pick a local coordinate $z$ on each disc (uniquely up to  $\Z/3\Z$=action $z\mapsto \xi z,\,\,\xi^3=1$) in such a way that 
$$\sigma(z)=-z,\,\,\,\,\,\,\,\,\,\omega_{1,0}(z)-\omega_{1,0}(-z)=4z^2\,dz.$$
The remaining input is the Bergman kernel $\omega_{0,2}$ which can be written as the sum of canonical kernels $dz_1 dz_2/(z_1-z_2)^2$ on each disc,
plus a regular (= without poles) $S_2$-symmetric section of $\pi_1^* K_\Sigma \otimes \pi_2^* K_\Sigma$.
  As a matrix-valued $2$-differential
 \[ \omega_{0,2}(z_1,z)=\left(\frac{\delta_{\alpha_1,\alpha_2}\,dz_1dz}{(z_1-z)^2}+\sum_{n_1,n_2\ge 1} P^{\alpha_1,\alpha_2}_{n_1,n_2}\, z_1^{n_1} z^{n_2} \frac{dz_1 dz}{ z_1 z}\right)_{\alpha_1,\alpha_2\in Ram},\]
 where $(P^{\alpha_1,\alpha_2}_{n_1,n_2}\in \C)_{\alpha_1,\alpha_2\in Ram;n_1,n_2\ge 1}$ obey the symmetry 
$P^{\alpha_1,\alpha_2}_{n_1,n_2}=P^{\alpha_2,\alpha_1}_{n_2,n_1}$.

Define the Tate topological vector space 
\[W:=W_{Airy}^{Ram}=\C^{Ram}\otimes\{\gamma \in \C(\!(z)\!)\,dz\,|\,Res_{z}( \gamma)=0\}\subset\C^{Ram}(\!(z)\!)dz.\]
It has topological basis $([\alpha]\otimes z^n dz/z)_{\alpha \in Ram, n\ne 0}$ where $([\alpha])_{\alpha \in Ram}$ is the coordinate basis of $\C^{Ram}$. We endow $W$ with the symplectic pairing
\begin{eqnarray*}\langle [\alpha_1]\otimes df_1,[\alpha_2]\otimes df_2\rangle&:=&\delta_{\alpha_1,\alpha_2} Res_{z}( f_1 df_2)\,\,\Longleftrightarrow \\
\langle [\alpha_1]\otimes z^{n_1}\frac {dz}{z},  [\alpha_2]\otimes z^{n_2}\frac {dz}{z}\rangle&:=&\frac{1}{n_2}\,\delta_{\alpha_1,\alpha_2}\delta_{n_1,-n_2}. \end{eqnarray*}
Define a Lagrangian subspace $W_+:=\C^{Ram}[[z]] dz\subset W$, which is  topologically  spanned by the base vectors $[\alpha]\otimes z^n\frac{dz}{z}$ with ${n\ge 1}$. It is complementary to the Lagrangian subspace $W_-\simeq W_+^*$ spanned by $[\alpha]\otimes z^{-n}\frac{dz}{z}$ with ${n\ge 1}$.
\begin{lmm} Under the above assumptions  formal Bergman kernels (formal in the sense that we do not put convergence condition on the regular part of $\omega_{0,2}$) are in one-to-one correspondence with Lagrangian subspaces $V\subset W$ complementary to  $W_+$.
The correspondence is defined as follows: expand $\omega_{0,2}(z_1,z)$ into the iterated Laurent series:
$$\omega_{0,2}(z_1,z)_{|z|\ll |z_1\!|}\in \C^{Ram}\widehat{\otimes} \C^{Ram} (\!(z_1)\!)(\!(z)\!)\,dz_1 dz.$$
Then  one has 
\begin{equation} \omega_{0,2}(z_1,z)_{|z|\ll |z_1\!|}\in V^{z\to z_1}\widehat{\otimes} \C^{Ram} [[z]] \,dz\label{Bergman and V}.\end{equation}
\end{lmm}
{ \it{Proof.}} In the above notation one has
\[ \omega_{0,2}(z_1,z)=\sum_\alpha \frac{[\alpha]\otimes[\alpha] \,dz_1 dz}{(z_1-z)^2}+\omega_{0,2}^{reg},\]
where $\omega_{0,2}^{reg}\in \widehat{Sym^2}W_+$ (recall that the hat means topological completion). The polar part is equal to 
\[ \sum_{n\ge 1,\alpha} \left(  [\alpha]\otimes  n  z^{-n}\frac{dz}{z}\right) \otimes  \left(  [\alpha]\otimes  z^{n}\frac{dz}{z}\right)=(\mbox{canonical pairing})^{-1}\in W_-\widehat{\otimes} W_+ .\]
The regular part is given by 
\[\omega_{0,2}^{reg}=\sum_{\alpha_1,\alpha_2}\,\sum_{n_1,n_2\ge 1} P^{\alpha_1 \alpha_2}_{n_1 n_2} \left( [\alpha_1] \otimes z_1^{n_1} \frac{dz_1}{z_1}\right)\otimes \left(  [\alpha_2] \otimes z_2^{n_2} \frac{dz_2}{z_2}\right)\,.\]
 Any Lagrangian subspace $V\subset W$ complementary to $W_+$ is the graph in $W\simeq T^* W_-$  of the differential of a quadratic form on $W_-$. Such a form corresponds to  an element of
 $\widehat{Sym^2}W_+$. $\blacksquare$
 
 From now on we denote by $V$ the Lagrangian subspace of $W$ associated to the given  Bergman kernel $\omega_{0,2}$ whose regular part could be even formal (non-convergent) matrix-valued power series in two variables.\footnote{ In fact the whole formalism of the topological recursion can be generalized straightforwardly to the  formal case. We are going to use this observation without further comments.} Space $V$ has canonical basis $(e^{\alpha}_n)_{\alpha\in Ram;n\ge 1}$:
\begin{equation}          \label{basis in V}
e^{\alpha}_n:=[\alpha]\otimes\,z^{-n}\frac{dz}{z}+\sum_{\alpha';n'\ge 1}\frac{1}{n} P^{\alpha', \alpha}_{n' ,n}\, [\alpha'] \otimes z^{n'}\frac{dz}{z}\,.\end{equation}
Vector $e^\alpha_n$ is the unique vector in $V$ such that
\[e^{\alpha}_n =[\alpha]\otimes z^{-n}\frac{dz}{z}+\left(\text{regular part }\in \C^{Ram}[[z]] \,dz\right)\,.\]
It will be convenient to write Bergman  kernel as
\begin{equation}\label{Bergmann kernel convenient}\omega_{0,2}(z_1,z)=\sum_{\alpha;n>0} (e_n^{\alpha})^{z\to z_1}\otimes \left([\alpha]\otimes nz^{n-1}{dz}\right)\end{equation}
where notation $(e_n^{\alpha})^{z\to z_1}$ means the element of $\C^{Ram}(\!(z_1)\!)\,dz_1$ obtained by the replacement of variables $z\to z_1$ from  $e_n^{\alpha}\in V \subset \C^{Ram}(\!(z)\!)\,dz$. 

The 
 recursion kernel $K(z_1,z)$ exists in general only as a $(Ram\times Ram$)-matrix with coefficients in 
 $\C(\!(z_1)\!)(\!(z)\!)\,\frac{dz_1}{dz}$. 
  Let us expand  the  recursion kernel in double series:
 \[ K(z_1,z)_{|z|\ll |z_1\!|}= \sum_\alpha \,[\alpha]\otimes[\alpha]\cdot\frac{-1}{4}\cdot\left(\frac{z^{-1}}{z_1^2}+\frac{z}{z_1^4}+\dots\right)\frac{dz_1}{dz}+\]
 \[+\sum_{\alpha_1,\alpha_2}[\alpha_1]\otimes [\alpha_2]\sum_{n_1,n_2\ge 1; n_2\,\,odd}\frac{-1}{4n_2}P^{\alpha_1,\alpha_2}_{n_1,n_2}z_1^{n_1} z^{n_2-2}\frac{dz_1}{dz}. \]
 This formula implies
\begin{equation}\label{recursion kernel convenient}
K(z_1,z)=-\frac{1}{4}\sum_{\alpha;n\ge 1\,\,odd}(e_n^{\alpha})^{z\to z_1}\,\otimes\left([\alpha]\otimes \frac{z^{n-2}}{dz}\right)
\end{equation}
Denote by $W_{odd}\subset W$  the coisotropic  subspace of meromorphic 1-forms on $\Sigma$ whose pushfoward to the quotient $\Sigma/\sigma$
 has no poles at the images of ramification points $(r_\alpha)_{\alpha\in Ram}$. Subspace $W_{odd}$ is the  topological span of base vectors
  $[\alpha]\otimes z^n\frac{dz}{z}$ with $n\ge 1$ arbitrary and $n\le -1$ odd. In particular, $W_+$ is a subspace of $W_{odd}$.
Let 
$$V_{odd}:=V\cap W_{odd}$$
be the subspace of $V$ spanned by basis vectors $(e^{\alpha}_n)_{\alpha\in Ram, n\ge 1\,\,odd}$.

  Combining together the above formulas we conclude that the residue formulas from the conventional topological recursion  (see  formulas
  \eqref{A-formula general},\eqref{epsilon-formula general},\eqref{B-formula general}, \eqref{C-formula general} from Section~\ref{section: TR}) give tensors (no completions for the tensor products):
   \[\begin{array}{ll}
   A_{TR}^{univ}\in   &V_{odd}\otimes V\otimes V,\\
   B_{TR}^{univ}:&W\to V_{odd}\otimes V,\\
   C_{TR}^{univ}:&W\otimes W\to V_{odd},\\
   \epp^{univ}_{TR}\in& V_{odd}.
   \end{array}\]
   
Restricting $B^{univ}_{TR}$ and $C_{TR}^{univ}$ to $V$ we see that $(A_{TR}^{univ},B_{TR}^{univ},C_{TR}^{univ},\epp^{univ}_{TR})$ induce a pre-Airy structure on $V$. We denote the corresponding tensors by
$(A_{TR}^{(1)},B_{TR}^{(1)},C_{TR}^{(1)},\epp^{(1)}_{TR})$.  Below we will calculate these tensors explicitly.
\begin{lmm} Tensor $A^{(1)}_{TR}$ is equal to $\frac{1}{4}\sum_{\alpha} e_1^{\alpha}\otimes e_1^{\alpha}\otimes e_1^{\alpha}$.
\label{lmm A}
\end{lmm}
{\it Proof.}
The definition of tensor $A^{univ}_{TR}$ by the formula \eqref{A-formula general} gives
\begin{align*}  A^{univ}_{TR}=&-\frac{1}{4}\sum_{\alpha}Res_z\left(e^{\alpha}_1\frac{1}{z \,dz}+\dots\right)\otimes\left(e^{\alpha}_1 dz+\dots\right)\otimes\left(-e^{\alpha}_1 dz+\dots\right)\\
=& \,\frac{1}{4}\sum_{\alpha} e_1^{\alpha}\otimes e_1^{\alpha}\otimes e_1^{\alpha}\,.\end{align*}

In particular, we see that $A^{univ}_{TR}$ is a {\it symmetric} tensor.$\blacksquare$

\begin{lmm} Vector $\epp^{(1)}_{TR}$  is equal to the sum of two  terms
\[\sum_\alpha\left(\frac{e^{\alpha}_3}{16}+\frac {P^{\alpha,\alpha}_{1,1} e^{\alpha}_1}{4}\right).\]
\label{lmm epsilon}
\end{lmm}
{\it Proof.} By definition {epsilon-formula general} we have
\[\epp^{(1)}_{TR}=-\frac{1}{4}\sum_{\alpha;n\ge 1\,\,odd}e^\alpha_n\times Res_z \left(
\frac{z^{n-2}}{dz}\,Coeff_{[\alpha]\otimes[\alpha]}\big(\omega_{0,2}(z,-z)\big)\right).\]
For each value of $\alpha\in Ram$ the canonical pole part of $\omega_{0,2}$ gives $e^\alpha_3/16$, and the regular part $\omega_{0,2}^{reg}$  gives $P_{1,1}^{\alpha,\alpha}e_1^\alpha/4$. $\blacksquare$

\begin{lmm} Operator $B^{(1)}_{TR}$ maps $e^{\alpha}_n$ to the sum of two terms
\[   \frac{1}{4}\sum_{n_1,n_2\ge 1;\,n_1+n_2=n+3;\, n_1\,{odd}}\,(-1)^{n_2+1}n_2\,e_{n_1}^{\alpha}\otimes e_{n_2}^{\alpha}+\frac{1}{4}\sum_{\alpha'}  P^{\alpha',\alpha}_{1,n}\,e^{\alpha'}_1\otimes e^{\alpha'}_1\,.\]
\label{lmm B}
\end{lmm}
{\it Proof.} For any $\alpha'\in Ram,n'\in \Z_{\ne 0}$ tensor $B^{univ}_{TR}([\alpha']\otimes z^{n'}\frac{dz}{z})\in V\otimes V$ is given (using \eqref{B-formula general}) by formula
\begin{equation}\label{354}
-\frac{1}{4}\sum_{n_1\ge 1\,\,odd}\,\,
\sum_{n_2\ge 1}\,e_{n_1}^{\alpha'}\otimes
e_{n_2}^{\alpha'}\cdot Res
\left(\frac{z^{n_1-2}}{dz}\frac{(-1)^{n_2}}{n_2}z^{n_2-1} dz\cdot z^{n'}\frac{dz}{z}\right).
\end{equation}
In particular, we see that this expression vanishes if $n'\ge 2$. Hence, in the calculation of $B^{(1)}_{TR}(e^\alpha_n)$ we should use only the polar term $(n'=-n$)
\[e^\alpha_n=[\alpha]\otimes z^{-n}\frac{dz}{z}+\dots\] (giving the main part of the formula), and constant terms ($n'=1$)
\[e^\alpha_n=\dots+\sum_{\alpha'} P^{\alpha',\alpha}_{1,n}[\alpha']\otimes z^1\frac{dz}{z}+\dots\]
contributing in \eqref{354} only for $n_1=n_2=1$.
$\blacksquare$
\begin{lmm} Operator $C^{(1)}_{TR}$ maps $e^{\alpha_1}_{n_1}\otimes e^{\alpha_2}_{n_2}$ to the sum of four terms:
\begin{enumerate}
\item $\frac{1}{4}(-1)^{n_2+1}e_{n_1+n_2+3}^{\alpha_1}$ if $\alpha_1=\alpha_2$ and $n_1+n_2$ even,
\item $\frac{1}{4}\sum_{n_2'\ge1; 2+n_2'-n_1\in 2\Z_{\ge 0}} \frac{(-1)^{n_2'+1}}{n_2'}\,P^{\alpha_1,\alpha_2}_{n_2',n_2}\, e^{\alpha_1}_{3+n_1-n_2'}$,
\item $\frac{1}{4}\sum_{n_2'\ge1; 2+n_2'-n_1\in 2\Z_{\ge 0}} \frac{(-1)^{n_2'+1}}{n_2'}\,P^{\alpha_1,\alpha_2}_{n_2',n_2} \,e^{\alpha_1}_{3+n_1-n_2'}$,
\item $\frac{1}{4}\sum_{\alpha} P^{\alpha,\alpha_1}_{1,n_1} \,P^{\alpha,\alpha_2}_{1,n_2} \,e^{\alpha}_1$.
\end{enumerate}
\label{lmm C}
\end{lmm}
{\it Proof.} For any $\alpha_1'\ne \alpha_2'$ and any $n_1',n_2'\in\Z_{\ne 0}$ we have
\[C_{TR}^{univ}\big(([\alpha_1']\otimes z^{n_1'}\frac{dz}{z})\otimes
([\alpha_2']\otimes z^{n_2'}\frac{dz}{z})\big)=0\,,\]
while for coinciding ramification points $r_{\alpha_1'}=r_{\alpha_2'},\,\,\,\,\alpha_1'=\alpha_2'=:\alpha'$ we have (using \eqref{C-formula general}):
\[C_{TR}^{univ}\big(([\alpha']\otimes z^{n_1'}\frac{dz}{z})\otimes
([\alpha']\otimes z^{n_1'}\frac{dz}{z})\big)=\]
\[=-\frac{1}{4}\sum_{n'\ge 1\,\,odd} e^{\alpha'}_{n'}\cdot Res\left(
\frac{z^{n'-2}}{dz} \,z^{n_1'-1} dz \,(-1)^{n_2'+1} z^{n_2'-1}dz\right)\,.\]
Therefore, we get a general formula
\begin{equation*}
C_{TR}^{(1)}\big(([\alpha_1']\otimes z^{n_1'}\frac{dz}{z})\otimes
([\alpha_2']\otimes z^{n_2'}\frac{dz}{z})\big)=\end{equation*}
\begin{equation}=\begin{cases}

\frac{(-1)^{n_2'}}{4} e^{\alpha_1'}_{3-n_1'-n_2'} &\text{if }\alpha_1'=\alpha_2',\,\,2-n_1'-n_2'\in 2\Z_{\ge 0}\\
0&\text{otherwise}\end{cases}\label{formula for general C}
\end{equation}
Using \eqref{basis in V}  write $e^{\alpha_1}_{n_1}\otimes e^{\alpha_2}_{n_2},\,\,\,n_1,n_2\ge 1$ as
\begin{align}\label{line1}& +([\alpha_1]\otimes z^{-n_1}\frac{dz}{z})\otimes ([\alpha_2]\otimes z^{-n_2}\frac{dz}{z})\\
\label{line2}&+ \sum_{\alpha_2';n_2'\ge1}\frac{1}{n_2'} \,P^{\alpha_2',\alpha_2}_{n_2',n_2}\,([\alpha_1]\otimes z^{-n_1}\frac{dz}{z})\otimes(
[\alpha_2']\otimes z^{n_2'}\frac{dz}{z})\\
\label{line3}&+\sum_{\alpha_1';n_1'\ge1}\frac{1}{n_1'}\, P^{\alpha_1',\alpha_1}_{n_1',n_1}\,([\alpha_1']\otimes z^{n_1'}\frac{dz}{z})\otimes(
[\alpha_2]\otimes z^{-n_2}\frac{dz}{z})\\
\label{line4}&+\sum_{\alpha_1'\alpha_2';n_1',n_2'\ge1}\frac{1}{n_1'\,n_2'} \,P^{\alpha_1',\alpha_1}_{n_1',n_1}\,P^{\alpha_2',\alpha_2}_{n_2',n_2}
\,([\alpha_1']\otimes z^{n_1'}\frac{dz}{z})\otimes(
[\alpha_2']\otimes z^{n_2'}\frac{dz}{z}).
\end{align}
Applying \eqref{formula for general C} to four lines \eqref{line1},\eqref{line2},\eqref{line3},\eqref{line4} we obtain the statement of the lemma (notice that for term \eqref{line4} the only relevant part will be $\alpha_1'=\alpha_2'=:\alpha$ and $n_1'=n_2'=1$).
$\blacksquare$

Let us consider the quantum pre-Airy structure  $(A^{(0)}_{TR},B^{(0)}_{TR},C^{(0)}_{TR},\epp^{(0)}_{TR})$ on the vector space $W_-$ which is equal to the direct sum of copies labeled by set $Ram$ of the canonical quantum pre-Airy structure from Section 
\ref{TR Airy case}. The vector space $W_+$ is the tangent space at $0$ to the  germ of the associated Lagrangian subvariety.
 We have shown  in Lemma 3.5.1 that a choice of the formal Bergman kernel corresponds to the choice of a complementary Lagrangian subspace $V$ to $W_+$.
 Hence it gives a way to modify the pre-structure  $(A^{(0)}_{TR},B^{(0)}_{TR},C^{(0)}_{TR},\epp^{(0)}_{TR})$  by a gauge equivalence (see formulas \eqref{classical gauge change} and \eqref{quantum gauge change} from Section 2). The role of symmetric matrix $s=(s^{ij})$ in these formulas is played  by $\omega_{0,2}^{reg}\in \widehat{Sym^2}W_+$
  understood as an infinite symmetric matrix:
 $$P:=\left(P^{\alpha_1,\alpha_2}_{n_1,n_2}\right)_{\alpha_1,\alpha_2\in Ram; n_1,n_2>0},\,\,\,\mbox{ index set} =Ram\times \Z_{>0}.$$
 We identify vector spaces $W_-$ and $V$  via
 \[W_-\stackrel{\sim}{\longleftarrow} W/W_+\stackrel{\sim}{\longrightarrow} V\,.\]
 
 \begin{thm} The gauge transformation of  $(A^{(0)}_{TR},B^{(0)}_{TR},C^{(0)}_{TR},\epp^{(0)}_{TR})$ as described above, coincides
 with  $(A^{(1)}_{TR},B^{(1)}_{TR},C^{(1)}_{TR},\epp^{(1)}_{TR})$.\label{TR general}
  \end{thm}
 
 {\it Proof.} Just read formulas  from Lemmas \ref{lmm A}, \ref{lmm epsilon}, \ref{lmm B}, \ref{lmm C}, as well as formulas \eqref{A Formulas for TR for Airy curve}, \eqref{B Formulas for TR for Airy curve}, \eqref{C Formulas for TR for Airy curve}, \eqref{epsilon Formulas for TR for Airy curve} for the standard pre-Airy structure, and compare with corresponding expressions \eqref{classical gauge change} and \eqref{quantum gauge change}. $\blacksquare$
 
 \begin{cor} Pre-Airy structure $(A^{(1)}_{TR},B^{(1)}_{TR},C^{(1)}_{TR},\epp^{(1)}_{TR})$ contains a primitive quantum Airy substructure on the subspace $V_{odd}$. 
 \end{cor}
 
 {\it Proof.} Here it is useful to recall  (essentially trivial) Remark \ref{modification pre Airy}.  The quantum pre-Airy structure $(A^{(0)}_{TR},B^{(0)}_{TR},C^{(0)}_{TR},\epp^{(0)}_{TR})$ on $V$ contains the quantum Airy substructure given by the restriction of all tensors to  $V_{odd}$. The gauge transformation from the above Theorem transforms this Airy substructure to the Airy substructure of the pre-Airy structure $(A^{(1)}_{TR},B^{(1)}_{TR},C^{(1)}_{TR},\epp^{(1)}_{TR})$. This proves the result. $\blacksquare$

There is a variation of the above theorem. Namely, consider the quantum Airy structure $(A^{(0)},B^{(0)},C^{(0)},\epp^{(0)})$ on $W_-$  which is the direct sum of copies labeled by set $Ram$ of the canonical quantum Airy structure from Section \ref{sec: canonical Airy structure}. Also, denote by $(A^{(1)},B^{(1)},C^{(1)},\epp^{(1)})$ the quantum pre-Airy structure on $V$ given by the modified recursion as in  Section \ref{sec: modified}.  
\begin{thm} Pre-Airy structure $(A^{(1)},B^{(1)},C^{(1)},\epp^{(1)})$ is in fact an Airy structure, and it is obtained from $(A^{(0)},B^{(0)},C^{(0)},\epp^{(0)})$ by the gauge transformation associated with $\omega_{0,2}^{reg}\in \widehat{Sym^2}W_+$. Moreover, this structure has the primitive substructure on subspace $V_{odd}$, and this substructure coincides with the one from Theorem \ref{TR general}.
\end{thm}
 {\it Proof.} The fact that as a pre-structure $(A^{(0)},B^{(0)},C^{(0)},\epp^{(0)})$ transforms to $(A^{(1)},B^{(1)},C^{(1)},\epp^{(1)})$ follows from the calculations almost identical to those in the proof of  Theorem \ref{TR general}, with the only correction that factors $(-1)^{n_2+1},(-1)^{n'_2+1}$  are removed. The rest of statements of the Theorem follows immediately. $\blacksquare$

\section{Spectral curves}

\subsection{Poisson surfaces with foliation}

Any complex analytic Poisson surface which does not consists entirely of $0$-dimensional symplectic leaves, contains an open dense symplectic leaf. The complement to the symplectic leaf is the set of zeros of the Poisson structure. Hence it is an analytic divisor, which consists of $0$-dimensional symplectic leaves.

\begin{defn} \label{def:good}
A {\bf good foliated Poisson surface} is  a complex  Poisson surface $P$ endowed with an analytic 1-dimensional foliation 
$\mathcal F$, such that the Poisson structure is non-degenerate on an open dense symplectic leaf $S\subset P$, the effective divisor $D$ of zeros of the Poisson structure  is  the disjoint union of smooth components with positive multiplicities, and $\mathcal F$ is tangent to $D$.
\end{defn}

The inverse to the non-degenerate Poisson structure on $S$ is the symplectic form.
Hence we can (and will) denote the data from the above Definition by $(P,\omega,\mathcal{F})$, where $\omega$ is the symplectic form on $S=P-D$ understood as a  meromorphic $2$-form on $P$ with poles at $D$.

\begin{defn} Let $(P,\omega,\mathcal{F})$ be a good foliated Poisson surface.
A smooth complex curve $\Sigma\subset P$ is called a $D$-{\bf transversal} spectral curve if
it is a smooth connected compact complex curve which intersects $D$ at finitely many points $\{p_\beta\}_{\beta\in Pol}$, and the intersection 
is transversal. 

Such curve is called $(\FF,D)$-{\bf transversal}\footnote{Probably a better name would be $\FF$-Morse.} if 
$\Sigma$ is tangent to leaves of $\mathcal{F}$ at finitely many points $\{r_\alpha\}_{\alpha\in Ram}$,
with the tangency order $1$ (i.e. locally near the tangency
point it looks like the curve $x=y^2$ tangent to the line $x=0$). Here $Pol$ (resp. $Ram$) is the notation for the above finite set  of ``poles'' (resp. ``ramificiation points''). 
\footnote{There are further  possibilities for generalizations here. In particular we can consider the case of immersed rather than embedded curves $\Sigma$.}
\end{defn}
We denote by $\BB_0$ (resp. $\BB_1$) the moduli space of $D$-transversal (resp.  $(\FF,D)$-transversal) curves.
The moduli spaces are defined in terms of complex analytic families over complex analytic bases in the usual way.
Notice that the families of sets $Pol=Pol(\Sigma)$ and $Ram=Ram(\Sigma)$ are locally constant over $\BB_0$ or $\BB_1$ (here we speak about $D$-transversal or resp. $(\FF,D)$-transversal spectral curves depending on the class of curves $\Sigma$). 

For any $D$-transversal spectral curve $\Sigma$ we denote the corresponding point by $[\Sigma]\in \BB_0$ (and, sometimes, simply by $\Sigma$ in order to alleviate the notation).

Obviously, $\BB_1$ is an open subset of $\BB_0$ (dense in typical situation). We will see below (see Lemma \ref{unobstructedness}) that  $\BB_0$
(and hence $\BB_1$) is smooth.

Notice that for a  $(\FF,D)$-transversal curve,  we can choose near each tangent point $r_\alpha$ local analytic coordinates $(x,y)$ on $S$ such that
the point $r_\alpha$ has coordinates $x=0,y=0$, the symplectic form on $S$ is equal to $\omega=dx\wedge dy$, the leaves of  foliation $\mathcal{F}$ are given by
the equations $x=const$ and the spectral curve $\Sigma$ is given by the equation $x=y^2$. Such coordinates are unique modulo analytic automorphism $x\mapsto \xi, y\mapsto \xi^{-1}y$, where $\xi=e^{2\pi i/3}$ is the primitive cubic root of $1$. We will call such coordinates {\it canonical}.

\begin{exa}
Let $C$ be a smooth projective curve of genus $g\ge 2$.
Let $\pi:T^{\ast}C\to C$ be the standard projection, and $\Sigma\subset T^{\ast}C$  be a smooth compact  connected curve. Take $P=S:=T^{\ast}C$ and $\mathcal{F}$ to be the foliation by cotangent fibers. Then $Pol=\emptyset$ and $\Sigma$ is an $D$-transversal spectral curve. For any integer $d\ge 1$, the connected component of $\BB_0$ corresponding to smooth spectral curves $\Sigma$ such that $\pi_{|\Sigma}:\Sigma\to C$ has degree $d$, is a Zariski open part of the base  of Hitchin fibration for the group $GL(d)$. Recall that this base is 
\[\BB=\oplus_{i=1}^d \Gamma\left(C,(\Omega^1_c)^{\otimes i}\right)\simeq \C^{\,d^2(g -1)+1}\,\]
and it can be interpreted as the linear system of compact effective divisors in $T^{\ast}C$ of degree $d$ over $C$.
The complement $\BB-\BB_0$ is the discriminant hypersurface parametrizing singular spectral curves. Zariski open subvariety $\BB_1\subset \BB_0$ parametrizes smooth $\Sigma$ such
 that $\pi_{|\Sigma}$ has only double ramification points.
\end{exa}
This example suggests that in general the moduli space $\BB_0$ should be considered as an open subspace of a larger space $\BB$ parametrizing compact effective divisors in $P$.  It is not clear whether the unobstructedness result (Lemma \ref{unobstructedness} from the next section) extends to the larger context
of effective divisors.
\begin{exa}
More generally, let $C$ be a smooth projective curve endowed with irregular data in the sense of [KoSo].  The irregular data consist of a collection of marked points, and Puiseux polynomials (called singular terms in the loc. cit.) with multiplicities.
They encode the irregular singular behavior of Higgs fields or connections at the marked points. In the case of Higgs bundles the irregular data modify the definition of the spectral curve from the previous example.

By [KoSo] one can construct canonically a Poisson surface $P$, which  is obtained by a sequence of blow-ups of the fiberwise
compactification of $T^{\ast}C$. 

Then there are two possibilities for the choice of $P$ and the symplectic
leaf $S\subset P$.  Both possibilities are described in the loc. cit. :

1) If one does not fix eigenvalues of the Higgs field (by definition, the residues of the restriction of the Liouville 1-form to $\Sigma$) then there is a choice of $P$ which has  an open dense symplectic leaf $S$  isomorphic to $T^{\ast}C$. Then the spectral curve $\Sigma$ is defined as the  compact closure in $P$ of  a curve in $S$. Eigenvalues of the Higgs field correspond to intersections of $\Sigma$ with $P-S$.

2) If one fixes the eigenvalues of  the Higgs field, then one makes additional blow-ups in 1) at the intersection points of $\Sigma$ with the divisor $P-S$. In this case $P$ contains an open symplectic leaf $S$ which is strictly bigger than $T^{\ast}C$.
The spectral curve $\Sigma\subset S$ is compact.

The Poisson surface $P$ endowed with the mermorphic foliation  $\mathcal{F}$ induced from the foliation along fibers $T^*C\to C$ is not in general a good foliated Poisson surface in the sense of definition \ref{def:good}. One has to remove various types of ``bad'' points: double points of divisor $D=P-S$, undeterminacy points of $\mathcal{F}$, components of $D$ which are generically transversal to $\mathcal{F}$. If we are fortunate, the closure $\Sigma$ of a smooth spectral curve in $T^*C$ will  not meet bad points and will be a smooth compact curve transversal to $D$. In this case we will obtain a point in  $\BB_0$.

\end{exa}

Notice that in [KoSo] we considered compactified spectral curves which intersect divisor $D$ only along the components with multiplicity 1, while here we allow transversal intersections of $\Sigma$  with components $D$ of arbitrary multiplicity. One can make a sequence of blow-ups at the transversal intersection point with a component of multiplicity $k>1$ and arrive to the situation when the multiplicity of the intersection is equal to one. Nevertheless, this procedure will produce {\it different} moduli space of spectral curves. The conclusion is that the framework of the present paper is indeed more general than the one in [KoSo], even for the case of cotangent bundles.

\begin{exa} One can consider a ``multiplicative'' version of the previous example. Namely, let us start with  $S=\C^\ast\times \C^*$ endowed with the symplectic form ${dx\over{x}}\wedge {dy\over{y}}$ and  endow with it with the foliation $\FF$ with leaves $x=const$. Let us compactify our symplectic surface to a toric Poisson surface $P$ and make 
a chain of Poisson blow-ups similar to those mentioned in the previous example. After removing bad points as in the previous example, we obtain a good foliated Poisson surface with all components of $D$ having multiplicity one. The corresponding spectral curves are naturally related to $q$-difference equations. Similarly to the previous example they could happen to be smooth and $D$-transversal, hence giving points in $\BB_0$. 
 
\end{exa}


\subsection{Deformation theory of spectral curves}\label{defspectral}

We are interested in the deformation theory of $\Sigma$ in the class of $D$-transversal spectral curves. The deformation complex is $\R\Gamma(\Sigma, N_\Sigma)$, where $N_\Sigma=(T_P)_{|\Sigma}/T_\Sigma$ is the normal bundle to $\Sigma$ in $P$. Using the symplectic structure on $S\subset P$ we identify 
\[N_\Sigma\simeq{{\Omega}}^1_\Sigma(D_\Sigma),\,\,\,D_\Sigma:=D\cap\Sigma=\sum_{\beta\in Pol}k_\beta p_\beta,\]
where  ${{\Omega}}^1_\Sigma$ is the sheaf of holomorphic $1$-forms on $\Sigma$, and $k_\beta\ge 1$ is the local multiplicity of $D$ at $p_\beta$, i.e. the order of pole of $\omega$
 at divisor $D$ near the point $p_\beta$.

\begin{lmm}\label{unobstructedness}
 The deformation theory of $\Sigma$ is unobstructed. 
 \end{lmm}
 
{\it Proof.} In the case when $Pol\ne \emptyset$ we have  $H^1(\Sigma, {{\Omega}}^1_\Sigma(D_\Sigma))=0$ by Serre duality. Hence the deformation theory of $\Sigma$ is  unobstructed.

In the more classical  case when $Pol=\emptyset$ the unobstructedness is presumably known. Alternatively, one can adapt the proof of  Proposition 8.2.1 in [KoSo], where a similar statement (unobstructedness) was proven  
 for compact complex Lagrangian submanifolds in quasi-projective varieties. In the present framework we do not require $S$ to be quasi-projective, it is just an analytic surface. Nevertheless, the proof from the loc. cit. still works in our case, because the degeneration of the Hodge-to-de Rham spectral sequence (see Step 2 of the proof in the loc.cit.) holds for a family
of compact complex analytic curves over the spectrum of a local Artin algebra.
The Lemma is proved. $\blacksquare$

 By Lemma \ref{unobstructedness}  the complex analytic   moduli space space $\BB_0$  of  $D$-transversal curves is a disjoint union of complex analytic manifolds. The dimension of the the germ $\BB_{[\Sigma]}$ of  $\BB_0$ at point $[\Sigma]$ is 
 \[\dim_\C\BB_{[\Sigma]}=\text{rk}\, H^0(\Sigma,\Omega^1_\Sigma(D_\Sigma))=\begin{cases} g &\text{ if }  Pol=\emptyset\\
  g+\sum_{\beta\in Pol}k_\beta-1&\text{ if }Pol\ne \emptyset\end{cases}\,.\]

\section{Affine symplectic connection and local embedding of the moduli space of spectral curves}

\subsection{Generalities on local Lagrangian embeddings}\label{embedding general}

Let us recall the notion of an affine flat connection on a vector bundle.

\begin{defn} For a complex analytic manifold $X$ and an analytic vector bundle $\H$ on $X$, an {\bf affine flat connection} on $\H$ is  a pair $(\nabla,\phi)$ where $\nabla$ is an analytic flat linear connection, and $\phi\in \Gamma(X,\Omega^1_X\otimes \H)$ is such that $\nabla(\phi)=0$. 
\end{defn}
An affine flat connection is the same as a flat connection $\nabla+\phi$ with values in the Lie algebra of affine transformations. Hence it gives an identification of fibers $\H_x$ as {\it affine spaces} for all points $x\in X$ close to any reference point $x_0$. Under this identification the zero point $0_x\in \H_x$ is mapped to a certain point in $\H_{x_0}$. Therefore, we obtain a germ of an analytic map
\[\text{germ of }X \text{ at  }x_0\mapsto \text{germ of }\H_{x_0} \text{ at  }0.\]
This map is an immersion if and only if $\phi$ is a monomorphism (if it is understood as a map of vector bundles $T_X\to \H$). Finally, let us assume that the
 bundle $\H$ carries a non-degenerate skew-symmetric pairing $\langle\bullet,\bullet\rangle$ which is covariantly constant, and map $\phi:T_X\to \H$ is a linear Lagrangian embedding. Then $(\nabla,\phi)$ induces a  Lagrangian embedding from the germ of $X$ at $x_0$ to the symplectic affine space $\H_{x_0}$. Moreover, these embeddings form a formally flat family, as we vary the reference point $x_0$. 
 
 \begin{rmk} \label{criterium closed form} Suppose that 
  with   any point $x_0\in X$ and any  pair of points $x_1,x_2$  sufficiently close to  $x_0$, we associate a vector $v(x_1,x_2)_{x_0}\in \H_{x_0}$, which is covariantly constant if we vary $x_0$ (and keep points $x_1,x_2$ fixed), and satisfying the cocycle condition
 \[v(x_1,x_2)_{x_0}+v(x_2,x_3)_{x_0}=v(x_1,x_3)_{x_0}\,.\]
 It is easy to see that such  data are equivalent to a \emph{closed} 1-form $\phi$ with values in a flat vector bundle $\H$ on $X$. More precisely, suppose we have a linear flat connection $\nabla$ and $\nabla$-covariant 1-form $\phi$. Given a triple of points $(x_0,x_1,x_2)$ which are sufficiently close to each other we can find a contractible open subset $U\subset X$ containing all of them. Then on $U^3\subset X^3$ we have a trivial bundle $pr_{(x_0,x_1,x_2)\mapsto x_0}^{\ast}(\H)$, endowed with the pull-backs of $\nabla$ and $\phi$. It can be trivialized, such that  the fiber of the bundle is identified with $\H_{x_0}$. Then for a short path $\gamma_{x_1,x_2}$ joining $x_1$ and $x_2$ we define $v(x_1,x_2)_{x_0}=\int_{\gamma_{x_1,x_2}}\phi\in \H_{x_0}$.
Conversely, one can reconstruct $\phi$ from $v$ by taking the first derivative of $v(x_1,x_2)_{x_0}$ w.r.t. $x_2$
at $x_2=x_1$.
 \end{rmk}
 
 \subsection{Local Lagrangian embedding of $\BB_0$}\label{HSigma section}
 
 Our goal is to apply the general construction from the previous subsection to the case $X=\BB_0$ and obtain for every smooth spectral curve $\Sigma$ the corresponding germ of the Lagrangian embedding. The construction will depend on the foliation $\FF$, but only through the collection of its $(2k-1)$-germs near components of $D$ at which $\omega$ has pole of order $k\ge 1$.
 
First, we define a symplectic vector bundle such that its fiber at $[\Sigma]\in \BB_0$ is the hypercohomology
\begin{equation}\label{formula H Sigma}
\H_\Sigma={\mathbb H}^0(\Sigma, Cone(d_{dR}:{{\OO}}_\Sigma(-D_\Sigma)\to {{\Omega}}^1_\Sigma(D_\Sigma))),
\end{equation}
where $D_\Sigma:=\sum_{\beta\in Pol} k_\beta p_\beta$ is an effective divisor on $\Sigma$ as in Section \ref{defspectral}. Symplectic pairing
 $\langle\bullet,\bullet\rangle$ on  $\H_\Sigma$ is given by the Serre duality.
 It is easy to see that 
 \[ \text{dim}\, \H_\Sigma=\begin{cases} 2g &\text{ if }  Pol=\emptyset\\
  2g+2\sum_{\beta\in Pol}k_\beta-2&\text{ if }Pol\ne \emptyset\end{cases}\,.\] 
  The 1-form $\phi$ interpreted as a morphism $T_{\BB_0}\to \H$ is given by the linear map
  \begin{equation}\label{morphism phi}T_{[\Sigma]}\BB_0=H^0(\Sigma, {\Omega}^1_\Sigma(D_\Sigma))\to \H_\Sigma\end{equation}
  which is  the boundary map in the standard long exact sequence. Equivalently, it comes  from the exact sequence of morphisms of cones
  \[Cone(0\to {{\Omega}}_\Sigma^1(D_\Sigma))\to Cone({{\OO}}_\Sigma(-D_\Sigma)\to {{\Omega}}_\Sigma^1(D_\Sigma))\,.\]
It is easy to see that $\phi$  always gives a Lagrangian embedding.
 
It is left to construct a linear flat symplectic connection on the bundle  $\H$ and prove that $\phi$ is covariant. In the case $Pol=\emptyset$ there is an obvious candidate, namely, the Gauss-Manin connection. For $Pol\ne \emptyset$
   the story  is more complicated.
     
     Let us use analytic topology, and consider a disjoint union of open discs 
     \[\sqcup_{\beta\in Pol}\D_{p_\beta}\subset \Sigma,\,\,\,\,p_\beta \in \D_{p_\beta}\,\,\forall \beta\]
     containing points $p_\beta$.  Let us choose a coordinate $x=x_{(\beta)}$ on each disc $\D_{p_\beta}$ for each $\beta$, such that  $x_{(\beta)}$ vanishes at the point $p_\beta$.
      Using these local coordinates we will construct an isomorphism 
      \begin{equation}\label{iso psi}\psi:\H_\Sigma\stackrel{\sim}{\longrightarrow}\Big(H^1(\Sigma -\{p_\beta\}_{\beta \in Pol},i_\epsilon(Pol);\Z)\otimes \C\Big) \oplus \bigoplus_{\beta\in Pol} \left(\oplus_{i=-k_\beta}^{k_\beta-2} \C\cdot x^i dx/\C x^{-1}d x \right)
    , \end{equation}
    where $i_\epsilon:Pol\hookrightarrow {\Sigma}-\{p_\beta\}_{\beta \in Pol}$ is the embedding which associates with $\beta$ the point in $\D_{p_\beta}$ with $x_{(\beta)}$-coordinate equal to $\epsilon, \,0<\epsilon\ll 1$. The definition makes sense as for small $\epsilon$ the relative cohomology groups are canonically identified.

      In order to describe the morphism $\psi$ it will be convenient to use the following explicit presentation of $ \H_\Sigma$ in terms of global differential forms with arbitrary poles at a finite {\it non-empty} collection of distinct points $(r_i\in \Sigma)_{i\in I}$, which are also distinct from points $(p_\beta)_{\beta \in Pol}$\footnote{The points $r_i$ at this time have nothing to do with ramification points $r_\alpha$. But our notation is suggestive: we will later choose the set $(r_i\in \Sigma)_{i\in I}$ in such a way that it will coincide with the set $(r_\alpha)_{\alpha\in Ram}$.}:
      \begin{equation}\label{represent H}\H_\Sigma\simeq \{\eta\in \Gamma(\Sigma,\Omega^1_\Sigma(D_\Sigma+\sum_i \infty r_i)|\,\forall i\,\,Res_{r_i} \eta=0\}/d \Gamma(\Sigma, \OO_\Sigma(-D_\Sigma+\sum_i \infty r_i)).\end{equation}
        For a class $[\eta]\in \H_\Sigma$ represented by meromorphic 1-form $\eta$, we define $\psi([\eta])$ as the sum of two terms $$\psi([\eta])=\psi_1(\eta)+\psi_2(\eta)$$ 
 described such as follows.

The term $\psi_2(\eta)$ is defined as the  sum over $\beta\in Pol$ of the image of the formal power series decomposition of $\eta$  in the coordinate $x$ at the point $p_\beta\in \D_{p_\beta}$, under the natural map
         \[  x^{-k_\beta}\C[[x]]dx \to \oplus_{i=-k_\beta}^{k_\beta-2} \C\cdot x^i dx/\C x^{-1}d x                                                   ,\,\,\,x=x_{(\beta)}\,\,.\]
       The term $\psi_1(\eta)\in H^1(\Sigma -\{p_\beta\}_{\beta \in Pol},i_\epsilon(Pol);\Z)\otimes \C$ is defined via the pairing with relative homology classes such as follows.
        Let $\gamma=\gamma(\epsilon)$ be a $1$-chain on $\Sigma-\{p_\beta\}_{\beta \in Pol}$ with the boundary lying in the finite set $i_{\epsilon}(Pol)$ and depending continuously on the parameter $\epsilon$ (and hence representing the same relative homology class).
        Let us assume that this chain   avoids the set  
        $\{r_i\}_{i\in I}$. Then the integral $\int_{\gamma(\epsilon)}\eta$ considered as a function of $\epsilon$ has an  expansion as
         $\epsilon \to 0$ of the form:
\begin{equation}\label{expansion}
         \int_{\gamma(\epsilon)}\eta=\sum_{i=-N}^\infty c_i \epsilon^i + c_{log}\cdot \log(\epsilon),\,\,\,N=\max_\beta(k_\beta)-1\,.
\end{equation}
         Then we define the cohomology class $\psi_1(\eta)$ by the formula
\begin{equation}\label{constant}
         \psi_1(\eta)([\gamma]):=\text{ coefficient }c_0\text{ in the above expansion}\,.
\end{equation}
         \begin{prp} Cohomology class $\psi_1(\eta)$ does not depend on the choice of representatives $\gamma(\epsilon)$. Furthermore it depends 
         only on the class $[\eta]$ in $\H_\Sigma$. Also, $\psi_1(\eta)$ depends only on the  class $[\eta]$ in $\H_\Sigma$.
         \end{prp}
        The proof is a straightforward check,  and is left as an exercise to the reader. The above Proposition implies that the map $\psi$ is well-defined.
        \begin{thm}\label{psi iso} Map $\psi$ defined in  \eqref{iso psi} is an isomorphism.
        \end{thm}
        {\it Proof.}  The dimensions of the RHS and LHS in \eqref{iso psi} are equal to each other, hence it is sufficient to check that $Ker\, \psi=0$.
         Let us assume that $\psi([\eta])=0$ for some meromorphic 1-form $\eta$ representing the cohomology class $[\eta]$ via  \eqref{represent H}.
           Then $\eta$ should have zero residue at each point $p_\beta$ (otherwise it will produce a non-zero class in $H^1(\Sigma -\{p_\beta\}_{\beta \in Pol},i_\epsilon(Pol);\C)$). Vanishing of other Taylor coefficients (guaranteed by the condition $\psi_2(\eta)=0$) implies that
           \[\eta\in \Gamma(\Sigma,\Omega^1_\Sigma(-\sum_{\beta \in Pol} (k_\beta-1)p_\beta+\sum_i \infty r_i)),\,\,\,Res_{r_i} (\eta)=0\,\,\,\forall i\,.\]
We see that $\eta$ is non-singular at all points $p_\beta, \beta\in Pol$. Therefore in the formula (\ref{expansion}) there are no terms $c_i, i<0$ and no term $c_{log}$. Hence we can define
$c_0$ as the limit as $\varepsilon\to 0$ of $\int_{\gamma(\varepsilon)}\eta$.
           The condition $\psi_1(\eta)=0$ implies that $\eta$ represents zero cohomology class in
           $H^1(\Sigma-\{r_i\}_{i\in I},\{p_\beta\}_{\beta\in Pol};\C),$
           hence it is equal to $df$ for some function $f$ with poles at $\{r_i\}_{i\in I}$ and vanishing at $\{p_\beta\}_{\beta\in Pol}$.
             Moreover, vanishing of $\psi_2(\eta)$ implies that for every $\beta\in Pol$ the function $f$ vanishes at $p_\beta$ with multiplicity
             $\ge k_\beta$. Therefore $[\eta]=0\in \H_\Sigma$. $\blacksquare$
             
            Now we are ready to define the desired flat connection $\nabla$ on the vector bundle $\H$ on the moduli space $\BB_0$ of smooth $D$-transversal spectral curves. For a given curve $[\Sigma] \in \BB_0$  and every $\beta\in Pol:=Pol(\Sigma)$ let us choose a local coordinate $x_{(\beta)}$ on a small disc 
            parametrizing leaves of the foliation $\FF$ which are close to $p_\beta$. We also assume that $x_{(\beta)}$ vanishes  on the leaf passing through $p_\beta$. This local coordinate defines 
            a local coordinate on the disc $\D_{p_\beta}'\subset\Sigma'$ associated to $\beta\in Pol(\Sigma)=Pol(\Sigma')$ for any spectral curve $\Sigma'$ which is  sufficiently close to $\Sigma$. We identify $\psi(\H_\Sigma)$ and $\psi(\H_{\Sigma'})$ using homotopy equivalence or surfaces with marked points (Gauss-Manin connection)
            \[ H^1(\Sigma -\{p_\beta\}_{\beta \in Pol},i_\epsilon(Pol);\Z)\simeq H^1(\Sigma' -\{p_\beta'\}_{\beta \in Pol},i_\epsilon'(Pol);\Z),\]
            and using the {\it identity} isomorphism on spaces $ \oplus_{i=-k_\beta}^{k_\beta-2} \C\cdot x^i dx/\C x^{-1}d x,  \, x=x_{(\beta)}$ for each $\beta$. 
            \begin{thm} \label{conn} The above construction gives a flat connection $\nabla$ on the bundle $\H$ which does not depend on the choices of local coordinates $x_{(\beta)}$, and also on the choice of auxiliary punctures $(r_i)_{i\in I}$. Moreover, the connection depends only  on the collection of $(2k-1)$-jets of foliation $\FF$ near components of divisor $D=P-S$ where symplectic form $\omega$ has pole of order $k$. Furthermore, the connection 
            $\nabla$ preserves the natural symplectic pairing on $\H$, and the $\H$-valued 1-form $\phi$  defined in \eqref{morphism phi} is closed.
            \end{thm}
	      {\it Proof}. The independence of the choice of the auxiliary finite nonempty set $\{r_i\}_{i\in I}\subset \Sigma -\{p_\beta\}_{\beta \in Pol}$ is obvious, as one can pass from this subset to arbitrary larger finite subset. The independence  of the choice of local coordinates $\{x_{(\beta)}\}_{\beta \in Pol}$ follows from the fact that a change of coordinates $x_{(\beta)}\to \widetilde{x}_{(\beta)}$ modifies 
            the isomorphism $\psi=\psi_1+\psi_2$ in certain universal way which is independent of a spectral curve $\Sigma'$ close to $\Sigma$.

The covariance of the symplectic pairing follows from the description of connection in terms of Hamiltonian reduction in Section 7. Notice that considerations of Section 7 are given for $\Sigma\in \BB_1$. On the other hand, as one can show that $\BB_1$ is dense in $\BB_0$, it follows that the symplectic pairing is covariant by the argument of continuity. Alternatively, one can describe the symplectic structure on $\H_\Sigma$ in terms of the RHS of the formula \eqref{iso psi} and deduce the covariance from that description.
            

           The fact that the 1-form $\phi$ is closed, follows from Remark \ref{criterium closed form} as well as  the alternative description of the local embedding of the germ of $\BB_0$ at $[\Sigma]$ to the affine space $\H_\Sigma$ given in the next subsection.   The symplectic pairing is explained in Section 7, see Remark \ref{sympl pairing}. Technically considerations there deal with the open dense subspace $\BB_1\subset \BB_0$. Continuity arguments show that the result holds for $\BB_0$ as well.
              $\blacksquare$
        
   \begin{rmk} In Section 8 we will use a slightly different approach which gives an alternative (and maybe more natural) proof of Theorem \ref{conn}.
   \end{rmk}

   \subsection{Connection and $1$-form via regularized integrals}\label{regularized integral}

In this subsection we will use  Remark \ref{criterium closed form} in order to prove that $\nabla(\phi)=0$.

For that let us take three points $x_0=\Sigma, x_1=\Sigma_0, x_2=\Sigma_1\in \BB_0$ which are sufficiently close to each other
and consider an analytic  path $\Sigma(\lambda), \lambda\in [0,1]$ such that $\Sigma(0)=\Sigma_0,\Sigma(1)=\Sigma_1$.

Following Remark \ref{criterium closed form} we would like to define the corresponding vector $v(\Sigma_0,\Sigma_1)_{\Sigma}\in \H_{\Sigma}$ which satisfies the cocycle condition. As we see the curve $\Sigma$ plays an auxiliary role in the considerations, so we can rename $\Sigma_0$ to $\Sigma$ and consider just a pair of sufficiently close spectral curves. 
The path $\Sigma(\lambda)$ can be thought of as an analytic map $f: T\to S$ of the total space $T$ of the analytic family of spectral curves $\Sigma(\lambda)$ with the above boundary conditions as $\lambda=0,1$.
Then $\eta(\lambda):={{f^{\ast}(\omega)}\over {d\lambda}}$ is a family of meromorphic $1$-forms on $\Sigma(\lambda)$.

Recall that the RHS of the formula \eqref{iso psi} consists of two summands. Let us consider the first one, which denote by $H_{top}$. Relative cohomology classes in that formula are dual to the corresponding relative homology classes.  

Let $T^{\prime}:=T-\cup_{\beta\in Pol}\{(\lambda,p_\beta(\lambda))\}, \lambda\in [0,1]$. We denote by $Q_\varepsilon$ the subset it $T^{\prime}$ given by $(0,\Sigma(0))\cup (1,\Sigma(1))\cup (\lambda, i_\varepsilon(\beta)(\lambda))$, where $\lambda\in [0,1]$. Here we use the self-explained notation like $(\lambda, p), \lambda\in [0,1], p\in \Sigma(\lambda)$ for points of $T$,
There is a natural isomorphism $g: H_1(\Sigma -\{p_\beta\}_{\beta \in Pol},i_\epsilon(Pol);\Z)\simeq H_2(T^{\prime},Q_\varepsilon;\Z)$ which follows from the observation that families of points $p_\beta(\lambda), i_\varepsilon(\beta)(\lambda), \lambda\in [0,1]$ is a union of ${|Pol|}$ closed intervals,
and for the interval $[0,1]$ we have the natural isomorphism $H_1([0,1],\{0\}\cup\{1\};\Z)\simeq \Z$. 

Let $[\delta(\varepsilon)]=g([\gamma(\varepsilon)])$ be the class of a relative $2$-chain corresponding to the class of the $1$-chain $\gamma(\varepsilon)$. Then $\int_{f^{\ast}([\delta(\varepsilon)])}f^{\ast}(\omega)$ has  the asymptotic expansion (as $\varepsilon\to 0$) given by the RHS of the formula (\ref{expansion}).
Take the constant terms of this expansion (i.e. the coefficient $c_0$). Then we obtain the functional
$c_0=c_0(\Sigma_0,\Sigma_1)$ on $H_{top}$.

In order to take care about the second summand in the RHS of the isomorphism (\ref{iso psi}) let us use as before the  coordinate $x(\lambda)=x_{(\beta)}(\lambda)$ on the disc $\D_{p_\beta}(\lambda)$ containing $p_\beta(\lambda)$, and identify all coordinates $x(\lambda)$ for $\lambda\in [0,1]$. 
Let us denote this coordinate by $x$. Then we take the image of the Laurent expansion in $x$  of $\int_{[0,1]}\eta(\lambda)d\lambda $ in $\oplus_{i=-k_\beta}^{k_\beta-2} \C\cdot x^i dx/\C x^{-1}d x$ and after that take the direct sum over all $\beta\in Pol$.

Let us denote the sum of the above two summands by $v(\Sigma_0,\Sigma_1)_\Sigma$. It is clear from the construction that for a triple of spectral curves $\Sigma_0, \Sigma_1, \Sigma_2$ we have
$v(\Sigma_0,\Sigma_1)_\Sigma+v(\Sigma_1,\Sigma_2)_\Sigma=v(\Sigma_0,\Sigma_2)_\Sigma$. Indeed having $2$-chains adjacent one to another one we naturally compose them into a single $2$-chain, which gives the additivity of the corresponding integrals over the chains. Similarly the Laurent expansions are added because the corresponding integrals of $1$-forms $\eta(\lambda)$ are added. 

\begin{prp} The map which assigns $\eta(0)$ to the tangent vector at $\lambda=0$  to the path $\Sigma(\lambda)$
coincides with $\psi$ (after the above identification with the RHS of the formula \eqref{iso psi}).

\end{prp}

{\it Proof.} The map $\psi$ assigns to a tangent vector from $T_{[\Sigma]}\BB_0$ the corresponding $1$-form on $\Sigma$. By construction $\eta(0)$ can be identified with the same. $\blacksquare$

This concludes the alternative description of the connection and the $1$-form as well as the proof of the Theorem \ref{conn}.

\begin{rmk}
In the simplest case $Pol=\emptyset$ considerations of this subsection have simple (and well-known) meaning.
Indeed,  the family of the first homology groups $H_1(\Sigma^{\prime},\Z)$ form a local system over a neighborhood of $[\Sigma]\in \BB_0$. The map $\kappa: l\mapsto \int_{l}\omega, l\in H_1(\Sigma^{\prime},\Z)\simeq H_1(\Sigma,\Z) $
defines a local embedding of $\BB_0$ into the symplectic vector space $H^1(\Sigma, \C)$. The integer homology groups of the spectral curve  carry a linear flat connection along $\BB_0$, since they do not change when we deform $\Sigma$. The closed $1$-form $\kappa$ defines an embedding of the germ $\BB_\Sigma$ into ${\mathcal H}_\Sigma\simeq H^1(\Sigma,\C)$.
\end{rmk}

\section{Formal discs and universal Airy structure}

\subsection{Space of formal discs}\label{sec: space of discs}

Let $(S,\omega)$ be a complex  analytic symplectic surface endowed with an analytic 1-dimensional foliation $\FF$ (i.e. we consider only the sympletic leaf in a good foliated Poisson surface, and endow the leaf with the induced foliation).

For every $n\ge 2$ denote by  $Discs^{(n)}$  the space of $n$-jets of embedded formal discs in $S$ which are tangent to a leaf of foliation $\FF$ with the tangency  order equal to one. Technically speaking, a point of  $Discs^{(n)}$ is given by a point $s\in S$, an epimorphism from the local ring  $\OO_s$ (germs of analytic functions at $s$) to a  local Artin algebra $R$ isomorphic to the algebra of truncated polynomials $\C[z]/(z^{n+1})$, such that
  the image of the local coordinate on the space of leaves belongs to $\mathfrak{m}_R^2-\mathfrak{m}_R^3$, where $\mathfrak{m}_R\subset R$ is the maximal ideal.
  
  The space $Discs^{(n)}$ is a smooth variety fibered over  $S$ with fibers isomorphic to $\C^\ast\times \C^{n-2}$. Also, we have an obvious projective system of fibrations
  \[Discs^{(2)}\leftarrow Discs^{(3)}\leftarrow\dots \]
  We define $Discs$ to be the projective limit, understood e.g. as a functor on finite-dimensional analytic spaces.
In the algebro-geometric framework we assume that $S$ is an algebraic surface and line subbundle $T_\FF\subset T_S$ is algebraic as well. Then $Discs$ is a scheme of infinite type.

Our goal in this section will be to construct germs of Lagrangian embeddings of $Discs$ into infinite-dimensional affine symplectic spaces, in terms of the  general formalism from Section \ref{embedding general}.

For a point $t\in Discs$ we denote by $\widehat{\D}_t$ the corresponding formal disc embedded in $S$. Algebra of function $\OO(\widehat{\D}_t)$ is isomorphic {\it non-canonically} to $\C[[z]]$. Denote  by $UnivDiscs$ the total space of the fibration over $Discs$ whose fiber at $t$ is fomal scheme $\D_t$. We have following universal maps:
\begin{itemize}
\item [(i)] fibration $UnivDiscs\to Discs$;
\item [(ii)] section $Discs\to UnivDiscs$ of the above fibration, associating to the disc its unique tangency point to $\FF$;
\item [(iii)] $UnivDiscs\to S$, the universal embedding $\widehat{\D}_t\hookrightarrow S$ for $t\in Discs$.
\end{itemize}

As in the case of $D$-transversal spectral curves, we have a natural isomorphism
\begin{equation}\label{iso tangent  for discs}T_t Discs=\Gamma\big(\widehat{\D}_t,(T_S)_{|\widehat{\D}_t}/T_{\widehat{\D}_t}\big)\simeq \Omega^1(\widehat{\D}_t).\end{equation}
Denote by $UnivDiscs^*$ the fibration over $Discs$ with fiber $\widehat{\D}_t^*$ at $t\in Discs$ given by the {\it punctured} formal disc, i.e. $\widehat{\D}_t$ with the removed tangency point. Algebra of functions on $\widehat{\D}_t^*$  is isomorphic to $\C(\!(z)\!)$.

We will define a canonical non-linear flat connection on the fibration $UnivDiscs^*\to Discs$, using the foliation $\FF$  (and {\it not} using the symplectic form $\omega$). 
In general, in the case of {\it  non-punctured} discs, if one wants to make a non-linear connection on $UnivDiscs\to Discs$ along a tangent vector $v\in T_t Discs$, one has to make a lift
\[v\in \Gamma\big(\widehat{\D}_t,(T_S)_{|\widehat{\D}_t}/T_{\widehat{\D}_t}\big)\rightsquigarrow \widetilde{v}\in  \Gamma\big(\widehat{\D}_t,(T_S)_{|\widehat{\D}_t}\big)\,.\]
In the case of punctured discs, one should find a lift in $\Gamma\big(\widehat{\D}_t^*,(T_S)_{|\widehat{\D}_t^*}\big)$.
 We define the {\bf canonical} lift to be a unique lift which is the section of the subbundle $T_{\FF}\subset T_S$.
  Such a lift exists and unique because $\FF$ is transversal to $\widehat{\D}_t^*$. As a result we obtain a universal non-linear connection $\nabla_{univ}$
   on the fibration $UnivDiscs^*\to Discs$. 
   \begin{exa} Let us consider in local coordinates $(x,y)$ on $S$ such that $T_\FF=\C\dot \partial_y$ (leaves of $\FF$ are given by $x=const$).
   Consider 1-parameter family of embeddings $(f_\lambda)_{\lambda\in \C}$ of the standard coordinate formal disc into $S$ with the first order tangency to $\FF$:
   \[\text{parameter } \lambda\in \C \mapsto \Big(f_\lambda:Spf (\C[[z]])\to Spec (\C[x,y]),\,\,\,z\mapsto (x=z^2+\lambda, y=z)\Big)\,.\]
   In this way we obtain an analytic  path $\C\to Discs, \,\lambda\mapsto t(\lambda)$ is in $Discs$. The tangent vector $\partial_\lambda$ to the path can be lifted naively using identification of standard coordinates $z$ on discs $\widehat{\D}_{t(\lambda)}$. The corresponding section of
    $(T_S)_{\widehat{\D}_{t(\lambda)}}$ is given by $\partial_x$. The {\it canonical} lift of $\partial_\lambda$ is different, it is equal to $-\frac{1}{2z}\partial_y$.
    The difference between the canonical and the naive lifts is a meromorphic section of $(T_S)_{\widehat{\D}_{t(\lambda)}}$ which is tangent to $\widehat{\D}_{t(\lambda)}$,
     and corresponds to the vector field $-\frac{1}{2z}\partial_z$ on $\widehat{\D}_{t(\lambda)}^*$:
     \[T_{f_\lambda}:-\frac{1}{2z}\partial_z\mapsto -\frac{1}{2z}\partial_y-\partial_x\,.\]
   \end{exa}

\begin{prp} The non-linear connection $\nabla_{univ}$ defined above is flat.
\end{prp}
This statement can be proven by direct calculation. If one replaces formal punctured discs by actual embedded annuli transversal to the foliation, 
 an analogous statement will be almost a tautology.
 
 Let us return to the formal case. Lie algebra of vector fields on a punctured formal disc $\widehat{\D}_t$ naturally acts on the symplectic vector space
  \[\H_t:=\{\gamma\in \Omega^1(\widehat{\D}_t^*)\,|\,Res(\gamma)=0\},\]
  where the symplectic form is given by \eqref{symp pairing} (and is coordinate-independent). Therefore, the non-linear connection $\nabla_{univ}$ gives rise to a flat linear symplectic connection on the infinite-dimensional\footnote{Technically, one should consider $\H$ as an inductive limit of projective limits of finite-dimensional vector bundles.}  bundle $\H$ with fibers defined above.
  
  Define a $\H$-valued 1-form $\phi$ on $Discs$ as the natural embedding (here we use only the symplectic form $\omega$, and not the foliation):
  \[T_t Discs=\Omega^1(\widehat{\D}_t)\hookrightarrow Ker\left(Res: \Omega^1(\widehat{\D}_t^*)\to \C\right)=\H_t\,.\]
  \begin{thm} Form $\phi$ defined above is closed under the covariant derivative, and defines a Lagrangian embedding of the tangent bundle
   $T_{Discs}$ into the symplectic bundle $\H$.
  \end{thm}
{\it Proof.} The fact that $\phi$ gives a Lagrangian embedding is obvious, as in local coordinates $\C[[z]] dz$ is a Lagrangian subspace in
 $Ker\left(Res:\C(\!(z)\!)dz\to \C\right)$. The fact that $\phi$ is closed can be checked by direct calculation. 
 
 Alternatively, it can be explained by direct construction of a germ of Lagrangian embedding in the non-formal case and the Remark \ref{criterium closed form}. Namely, consider a germ of a finite-dimensional family of embedded analytic discs $\D_{t(\lambda)}$ where $\lambda\in \Lambda$ (parameter space) is close to $\lambda_0$.
  then one can pick an open domain $\mathcal{U}\subset S$ close to the tangency point of  $\D_{t(\lambda_0)}$ such that 
  $ \mathcal{U}\cap \D_{t(\lambda_0)}$ is an annulus $A\simeq \{z\in \C\,|\,\epsilon_1<|z|<|\epsilon_2\}$ for $0<\epsilon_1<\epsilon_2\ll 1$,
  and $\mathcal{U}$ is {\it fibered} over $A$ with fibers being leaves of $\FF$.
    Symplectic structure $\omega$ gives an {\it affine} structure on fibers $\mathcal{U}\to A$, and moreover makes fibers to be parallel to the cotangent bundle $T^* A$\footnote{This is a general feature of  Lagrangian foliations in arbitrary dimensions.}. Then for any $\lambda,\lambda'$ close
     to $\lambda_0$ we have two sections $ \mathcal{U}\cap \D_{t(\lambda)}$ and $ \mathcal{U}\cap \D_{t(\lambda')}$ of the fibration
     $\mathcal{U}\to A$. Their difference is a 1-form $\eta_{\lambda_1,\lambda_2}$ on $A$ (here we use affine structure along leaves of $\FF$).
     
\begin{lmm}      
The contour integral of this 1-form over a non-trivial loop in $A$ is equal to zero.
\end{lmm}

{\it Proof.}  Indeed, one can choose locally near the tangency point of $\D_{t(\lambda_0)}$ a local system of coordinates $(x,y)$ (canonical coordinates) such that $\omega=dx\wedge dy$, and the foliation is given by $x=const$. The contour integral of  $ydx$ over a loop in  $\D_{t(\lambda)}$ vanishes for any $\lambda$ close to $\lambda_0$, because this loop bounds a disc. Therefore, the countour integral over a loop in $A$
 of the pullback $\xi_\lambda$ of $ydx$ via the map $A\to \mathcal{U}\cap \D_{t(\lambda)}$ (inverse to the projection) also vanishes.
         Finally, the 1-form $\eta_{\lambda_1,\lambda_2}$ is equal to the difference $\xi_{\lambda_1}-\xi_{\lambda_2}$. Therefore its contour integral vanishes as well. 
The Lemma is proved. $\blacksquare$
        
We see that in the analytic situation we obtain a flat familiy of local embeddings of the space of discs into infinite-dimensional
  affine space. The cocycle condition \ref{criterium closed form} for the closedness of our vector-valued 1-form follows from the trivial  identity
         \[\eta_{\lambda_1,\lambda_2}+\eta_{\lambda_2,\lambda_3}=\eta_{\lambda_1,\lambda_3}\,.\]
     This explains that $\nabla \phi=0$. $\blacksquare$
     
\subsection{Relation to the universal classical Airy structure}

 For a given symplectic surface $(S,\omega)$ endowed with (automatically Lagrangain) one-dimensional foliation $\FF$, define an infinite-dimensional space $Coord$ to be the space of tuples $(s,x,y)$ where $s\in S$ is a point, $(x,y)$ is a formal coordinate system at $s$ such that
 \begin{itemize}
 \item[(i)] symplectic form $\omega$ is equal to $dx\wedge dy$;
 \item[(ii)] tangent bundle $T_\FF$ to the foliation $\FF$ is spanned by $\partial_y$.
 \end{itemize}
 This space is projective limit of a tower of fibrations of finite-dimensional smooth varieties (like the space $Discs$ in the previous subsection).

  The tangent space to $Coord$ at any point is canonically isomorphic to the Lie algebra $\mathfrak{g}$ where:
  \[\mathfrak{g}:=\left\{\xi=a(x,y)\,\partial_x+b(x,y)\,\partial_y\,|\,a,b\in\C[[x,y]],\,\,\xi \text{ preserves }\omega\text{ and }T_\FF\right\}\,.\]
   As in the formal differential geometry of Gelfand and Kazhdan (see [GelKaz]), we have  a {\it  free transitive action} of $\mathfrak{g}$ on $Coord$.
    An easy calculation shows that $\mathfrak{g}$ is canonically isomorphic to the quotient of the Lie algebra of differential operators of order $\le 1$ on $\C[[x]]$ by the one-dimensional center consisting of constants, i.e. $\g\simeq \C[[x]]\partial_x\ltimes (\C[[x]]/\C)$, where
    \footnote{ Here $g(x)/const$ denotes the image of the element $g=g(x)$ of $\C[[x]]$ in the quotient space $\C[[x]]/\C$.  }
    \[f(x)\,\partial_x + g(x)/const\mapsto  f(x)\,\partial_x-y f'(x)\,\partial_y-g'(x)\,\partial_y=\{f(x)y+g(x),\cdot\}\,.\]
     Consider the map $\pi: Coord\to Discs$ given by 
\[(s,x,y)\mapsto \text{ formal disc with coordinate }z\text{ embedded via } \big(z\mapsto(x=z^2,y=z)\big)\,.\]
 Map $\pi$ is a $3:1$ cover. Indeed $(s,x'=\xi^2 x, y'=\xi )$ gives the same embedded formal disc as $(s,x,y)$ for any $\xi\in \C, \xi^3=1$.
 
 For any point $(s,x,y)\in Coord$ and the corresponding point $t=\pi(s,x,y)\in Discs$, the local transitive free action of Lie algebra $\mathfrak{g}$ on $Coord$ gives an isomorphism
 \[\mathfrak{g}\simeq T_{(s,x,y)}Coord\simeq T_{t} Discs\simeq \Omega^1(\widehat{\D}_t)\simeq \C[[z]]dz\,.\]
 \begin{prp} \label{prop:iso tangent forms}The above isomorphism does not depend on $(s,x,y)\in Coord$, and is given by the formula
\[ f(x)\,\partial_x + g(x)/const\mapsto d(f(z^2)z +g(z^2))\in \C[[z]] dz\,.\]
 \end{prp}
  The proof is given by direct calculation.
 
 The pullback $\pi^*\H$ is a canonically trivialized symplectic vector bundle, each fiber is identified with the space $W:={W}_{Airy}:=Ker(Res:\C(\!(z)\!)dz\to \C)$ from Section \ref{sec: canonical Airy structure}.
 
 \begin{thm} Consider the affine flat connection on the trivial bundle with fiber ${W}$ over  $Coord$ given by the homomorphism
  \[\mathfrak{g}\to spAff(W)= sp(W)\ltimes W\]
 which is associated with the standard classical Airy structure from Section  \ref{sec: canonical Airy structure}. Then, after the above identifications, it coincides with the affine flat connection constructed in the  Section \ref{sec: space of discs}.
 \end{thm}
{\it Proof.} We have to check first that the linear connections match. Let us calculate the linear part of the flat affine connection from Section \ref{sec: space of discs} following the definitions.  We obtain it as the composition of two homomorphisms
\[\mathfrak{g}\to \C(\!(z)\!)\,\partial_z\to sp(W)\]
where the second maps is the natural action by coordinate changes.
 The first map comes from the geometric description of the non-linear flat connection $\nabla_{univ}$ via foliation.
  Namely, with any element  $f(x)\partial_x + g(x)/const$ we associate the corresponding Hamiltonian vector field $v$ in $Spf (\C[[x,y]])$.
  Let us consider restriction of $v$ to the parabola $x=y^2$, and find another meromorphic vector field $v'$ on parabola which is vertical
  and such that $u=v-v'$ is tangent to the parabola. Then $u$ expanded in coordinate $z$ gives an element in $\C(\!(z)\!)\,\partial_z$ corresponding to $v\in \mathfrak{g}$. The direct calculation gives
  \[v_{\widehat{D}_t}=f(z^2)\,\partial_x-zf'(z^2)\,\partial_y-g'(z^2)\,\partial_y\,.\]
  The map from vector fields tangent to  parabola in coordinate $z$, to section of the pullback $T_S$ to the parabola is given by
  \[h(z)\,\partial_z\mapsto 2z h(z)\,\partial_x+h(z)\,\partial_y\,.\]
  Therefore, the vector field tangent to parabola which differ from $v$ by a purely vertical field is 
  \[u=\frac{f(z^2)}{2z}\,\partial_z\,.\]
  This gives a homomorphism $\mathfrak{g}\to sp(W)$. Direct calculation shows that it coincides with the linear part of the homomorphism from Section \ref{sec: canonical Airy structure}. Similarly, the shift part (i.e. 1-form $\phi$ with values in the trivial bundle with the fiber $W$) matches, as we identify the tangent space to $Coord$ with the space of 1-forms $\C[[z]]\,dz$ as in Proposition \ref{prop:iso tangent forms}.
$\blacksquare$

\subsection{Geometric meaning of residue constraints}\label{meaning res}

Let $t\in Discs$ and $\widehat{\D}_t$ the corresponding formal disc. We endow it with canonical coordinates
$(x,y)$, so that the foliation $\FF$ is given by the equation $x=const$ and the symplectic form $\omega=dx\wedge dy$. Disc $\widehat{\D}_t$
 is endowed with coordinate $z$ such that the embedding into $S$ is $z\mapsto (x=z^2,y=z)$.
 
Tate vector space $W=W_{Airy}=Ker (Res_z:\C(\!(z)\!)dz \to \C)$ is  a hyperplane in  $W':=\C(\!(z)\!)dz$. The formal germ $\widehat{W}'$ of $W^{\prime}$ at $0$ controls deformations of the punctured formal disc $\widehat{\D}_t^{\ast}$. Its formal subgerm $\widehat{W}$  is determined by the condition that the residue of the pull-back of the $1$-form $ydx$ to the deformed punctured formal disc is trivial. The  formal germ $\widehat{W}$ contains a formal germ of $t\in Discs$ corresponding to deformations of the non-punctured formal disc $\D_t$.\footnote{We understand all submanifolds as functors on commutative local Artin algebras, see Remark \ref{Tate}.}
The latter germ is \emph{contained} in the formal germ of the Lagrangian submanifold $L$ which by definition  is singled out in $\widehat{W}$ by the system of equations

$$Res_z(x^ky\,dx)=0, k\ge 1,
Res_z(x^ky^2dx)=0, k\ge 0\,.$$
These are exactly the residue constraints from Section \ref{sec: canonical Airy structure}.
They hold on the formal completion of $Discs$ since the residue of a holomorphic form is trivial.

Notice that the tangent space at zero $T_0L\simeq \C[[z]]dz$ and hence it coincides with the tangent space to $Discs$. Hence  $L$ coincides the formal neighborhood of $t\in Discs$ embedded in $\widehat{W}$. This explains
the residue constraints from Section \ref{sec: canonical Airy structure}.

We will also explain below the invariant meaning the residue constraints.

Let $U_1$ be the vector space of $1$-forms on the formal completion of
of $S$ at $(x,y)=(0,0)\in \D_t$ which are equal to zero on the leaves of the foliation $\FF$ with the order $\le 2$ (the latter notion makes sense even globally because of the intrinsically defined affine structure on the leaves). 

In the canonical coordinates $(x,y)$ an element of $U_1$ can be written as $F(x,y)dx$, where $F(x,y)$ is a polynomial of degree $\le 2$ in $y$ for a fixed $x$. Let $U_2\subset U_1$ be a vector subspace consisting of such $1$-forms $\mu$ that $d\mu=c\cdot dx\wedge dy, c\in \C$.

Spaces $U_1, U_2$ have topological bases 
\[\{x^{\ge 0} dx, x^{\ge 0} y dx, x^{\ge 0} y^2 dx\},\quad
\{x^{\ge 0} dx, ydx\}\]
respectively. 
Therefore , the quotient $U_1/U_2$ has the topological basis 
\[\{x^{\ge 1}  y dx,x^{\ge 0} y^2 dx\}\,.\]
 The pull-back of $\mu\in U_2$ to a deformed $\D_t^{\ast}$ has trivial residue (we deform the disc inside of $\widehat{W}\subset \widehat{W}^{\prime}$. Hence $U_1/U_2$ is naturally mapped to the space of formal functions on the formal neighborhood of $0\in \widehat{W}$. Vanishing of these functions is the same as the residue constraints from Section \ref{sec: canonical Airy structure}.

\section{Hamiltonian reduction and Holomorphic Anomaly Equation}

\subsection{Global Hamiltonian reduction}\label{global reduction}

In this subsection we explain how to treat the germ  $\BB_\Sigma$ as a Hamiltonian reduction of a germ of an  infinite-dimensional Lagrangian submanifold of a  Tate  affine symplectic space (see Remark \ref{Tate}) that we understand Tate spaces as functors on Artin algebras.

We start with a finite-dimensional model. Let $M$ be a finite-dimensional symplectic manifold, $G\subset M$ a coisotropic submanifold and $L\subset M$ a Lagrangian submanifold. Let $x\in L\cap G$ be a point satisfying the following equivalent properties:

\begin{itemize}
 \item [1)] $T_xL+T_xG=T_xM,$
\item[2)] $T_xL\cap (T_xG)^\perp=0,$

\end{itemize}
where $\perp$ means the skew-orthogonal with respect to the symplectic form.

Notice that the family of vector spaces $(T_xG)^\perp$ defines an integrable distribution of planes. Hence we have the corresponding foliation of $G$ with (locally) symplectic quotient.
The following proposition is well-known and is easy to prove.

\begin{prp}
The image ${\mathcal L}_x$ of the germ of ${ L}_x\cap G_x$ at $x$ is a germ of Lagrangian submanifold in the (germ of the) above symplectic quotient.
\end{prp}

This Proposition can be generalized to the case of Tate spaces and formal germs.

Let now $\Sigma$ be a $(\FF,D)$-transversal curve $\Sigma$ in  a good foliated Poisson surface. Let us choose local coordinates $(x,y)$ on $S$ and $z$ on $\Sigma$ near each ramification point so we get  locally a standard parabola, as in the previous section. The germ $\BB_\Sigma$  of the moduli space $\BB_1$ at $[\Sigma]$\footnote{We will use the notation $\BB_\Sigma$ and $\BB_{[\Sigma]}$ for this germ on the equal footing.}
 maps to the product  $Discs^{Ram}$ of copies of the space $Discs$ labeled by ramification points. We denote by $\widehat{\BB}_\Sigma$ the formal germ of $\BB_\Sigma$ at the point $[\Sigma]\in \BB_1$.
 Thus, we obtain a local immersion of the finite-dimensional formal germ $\widehat{\BB}_\Sigma$ into an infinite-dimensional affine vector space $W:=W^{Ram}_{Airy}$. On the other side, we have a natural germ $L$ of infinite-dimensional Lagrangian subvariety in $W^{Ram}_{Airy}$, given by the formal neighborhood of a point in $Discs^{Ram}$. Clearly $L=L_{Airy}^{Ram}$ is the  product of copies  of formal germs of infinite-dimensional Lagrangian submanifolds $L_{Airy}$ from Section 3.3.

In order to simplify the notation we will denote by $\Delta=\sum_{\alpha\in Ram}r_\alpha$  the ramification divisor, and as before we will denote by $D_\Sigma=\sum_{\beta\in Pol}k_\beta p_\beta$ the pole divisor.

  \begin{prp} \label{prop with discs}
  The finite-dimensional germ  $\widehat{\BB}_\Sigma$ is equal to the intersection of the infinite-dimensional germ of $Discs^{Ram}$
  with the coisotropic affine subspace $G_\Sigma\subset W^{Ram}_{Airy}$ defined as the space of $1$-forms
\begin{equation}\label{coisotropic space}
  G_\Sigma:=\{\eta\in \Gamma(\Sigma,\Omega^1_\Sigma(D_\Sigma+\sum_{\alpha\in Ram} \infty \,r_\alpha))|\,\,Res_{r_\alpha} (\eta)=0\}
\end{equation}
  expanded in canonical local coordinates $z$ near each ramification poijnt $r_\alpha$
 (compare with \eqref{represent H}). Moreover, the subbundle with the fiber $G_\Sigma$ of the trivial bundle with the fiber $W^{Ram}_{Airy}$, is covariantly constant as an affine subbundle. 
  \end{prp}
  {\it Proof.} We will explain below both statements in the complex analytic setting. There is a similar but less transparent proof in the framework of Tate spaces.

Notice that by our assumptions, the foliation $\FF$ is transversal to $\Sigma$ everywhere except of ramification points. Therefore, there exists an open domain ${\mathcal U}\subset P$ such that 
  \[{\mathcal U}\cap \Sigma=\Sigma^0:=\Sigma-\sqcup_{\alpha\in Ram} \D_{\alpha,\epsilon_1}\]
  where $\D_{\alpha,\epsilon_1}$ is a closed disc of radius $\epsilon_1$ in the standard coordinate $z=z_{(\alpha)}$ near $r_\alpha,\,\,\forall\alpha\in Ram$, and $\epsilon_1,\,0<\epsilon_1\ll1 $ is sufficiently small.
  Moreover, we may assume that the projection to the space of  leaves of $\FF$ gives a fibration of ${\mathcal U}$ over $\Sigma^0$.
    Any spectral curve $\Sigma'$ which is sufficiently close to $\Sigma$ gives a section $\Sigma'\cap {\mathcal U}$ of the fibration ${\mathcal U}\to \Sigma^0$,
     which differs form the section given by $\Sigma$ by a 1-form on $\Sigma^0$ with poles of order $\le k_\beta$ at $p_\beta,\,\beta\in Pol$.
      This follows from the fact that fibers of the fibration  ${\mathcal U}\to \Sigma^0$ are naturally affine lines parallel to the cotangent bundle to $\Sigma^0$, twisted by divisor $D_\Sigma$ (here we use the symplectic structure to ensure that fibers of a Lagrangian foliation carry affine structure).  This 1-form has contour integrals over all annuli $\epsilon_1<|z^{(\alpha)}|<\epsilon_2\ll 1$
      equal to zero (because these annuli are contained in discs). This gives  an element in the (analytic version of) $G_\Sigma$ .
      
       We see that spectral curves $\Sigma'$ can be described as deformations of annuli $\epsilon_1<|z^{(\alpha)}|<\epsilon_2$ along the foliation 
       $\FF$ satisfying the vanishing of residue condition as well as two other constraints:
       \begin{itemize}
       \item each annulus is contained in a holomorphic  disc (geometrically it means that we land into the product of infinite-dimensional Lagrangian subvarieties given by quadratic equations of the standard classical Airy structure)
       \item the collection of 1-forms on annuli extends to a global 1-form from the analytic version of $G_\Sigma$, equivalently the deformed annuli lie in a shifted section of the fibration ${\mathcal U}\to \Sigma^0$.
       \end{itemize} 
      This description identifies the formal germ  $\widehat{\BB}_\Sigma$ with the intersection of a Lagrangian subvariety coming from the classical Airy structure and an affine coisotrpic subvariety.
  $\blacksquare$

Next, observe that the $2$-term complex of sheaves
${\OO}_\Sigma(-D_\Sigma)\to {\Omega}^1_\Sigma(D_\Sigma)$ endowed with de Rham differential is quasi-isomorphic to the following complex of sheaves:
$${\OO}_\Sigma(-D_\Sigma+\ast \Delta)\to {\Omega}^1_\Sigma(D_\Sigma+\ast \Delta)\to \oplus_{r_\alpha\in Ram}\C\to 0,$$
where $\Delta=\sum_{\alpha\in Ram}r_\alpha$, first two terms are sheaves of functions/ forms having poles/zeros at the points $p_\beta$ bounded by $k_\beta$ as well as and poles of arbitrary order at the points $r_\alpha$ (the map is the de Rham differential), and second (surjective) arrow is the map which assigns to a differential form the sum of residues at all points $r_\alpha$.
Each sheaf in the above complex has trivial cohomology, except of (possibly) in degree zero. It follows that ${\mathcal H}_\Sigma\simeq G_\Sigma/d({\OO}_\Sigma(-D_\Sigma+\ast \Delta))$. Thus we have the natural
map $G_\Sigma\to \H_\Sigma$.

\begin{prp}\label{reduction} We have the following 
isomorphism of Tate spaces: 
$$d(H^0({\OO}_\Sigma(-D_\Sigma+\ast \Delta)))= G_\Sigma^\perp,$$
and hence $G_\Sigma/G_\Sigma^\perp\simeq \H_\Sigma$.

\end{prp}

{\it Proof.} 
 It suffices to check that
$d(H^0({\OO}_\Sigma(-D_\Sigma+\ast \Delta)))^\perp= G_\Sigma$. 

By definition the LHS consists of collections of forms $(\mu_\alpha)_{\alpha\in Ram}$ such that $\mu_\alpha\in \C(\!(z_{(\alpha)})\!)dz_{(\alpha)}$ (here $z_{(\alpha)}$ is the canonical local parameter at $r_\alpha$), satisfying the conditions
\begin{itemize}
\item $Res(\mu_\alpha)=0\,\,\forall \alpha$,
\item
 $\sum_\alpha Res(f\mu_\alpha)=0$ for any $f\in H^0({\OO}_\Sigma(-D_\Sigma+\ast \Delta))$.
 \end{itemize} Let us introduce the Tate vector space 
$G_\Sigma^{\prime}$ as the space of collections of $1$-forms  $\mu_\alpha, \alpha\in Ram$ as above, but without the condition $ Res(\mu_\alpha)=0\,\,\forall \alpha$.

Before proceeding further, let us recall the following well-known result about Tate spaces associated with vector bundles on curves (Serre duality for curves formulated in terms of local fields, see [T], Corollary to Theorem 5).

Let $\Sigma$ be a smooth projective curve over a field ${\bf k}$, and ${\mathcal E}$ be a vector bundle over $\Sigma$. Let $\{r_i\}_{i\in I}$ be a finite subset of $\Sigma({\bf k})$. Let us denote by ${\mathcal E}^{\prime}$ the vector bundle ${\mathcal E}^{\vee}\otimes \Omega^1_\Sigma$. For each $r_i$ denote by $K_{r_i}$ the completed local field of $r_i$. Each $K_{r_i}$ is isomorphic to the field $\K(\!(z_{(i)})\!)$, where $z_{(i)}$ is a local coordinate at $r_i$. For each $i\in I$ we have the finite-dimensional $K_{r_i}$-vector spaces $\widehat{\mathcal E}_{r_i}$ and $\widehat{\mathcal E}_{r_i^{\prime}}$ of 
sections over the formal punctured discs at $r_i$. They are Tate $\K$-vector spaces endowed with the topological perfect duality $(v,w)\mapsto Res_{r_i}(v,w)$.
 Direct sum $\oplus_{i\in I}\widehat{\mathcal E}_{r_i}$ contains the subspace $\Gamma(\Sigma-\{r_i\}, {\mathcal E})$, and similarly $\oplus_{i\in I}\widehat{\mathcal E}_{r_i}$ contains the subspace $\Gamma(\Sigma-\{r_i\}, {\mathcal E}^{\prime})$. The direct sums themselves are naturally dual to each other.

\begin{lmm} 

The duality between above Tate spaces makes the subspaces orthogonal to each other.
\end{lmm}

We will apply the Lemma to ${\mathcal E}={\OO}_\Sigma(-D_\Sigma), {\mathcal E}^{\prime}={\Omega}_\Sigma^1(D_\Sigma), \{r_i\}_{i\in I}=\{r_\alpha\}_{\alpha\in Ram}$.

Then the lemma implies that

$$G_\Sigma^{\prime}\simeq \{\eta\in \Gamma(\Sigma,\Omega^1_\Sigma(D_\Sigma+\sum_{\alpha\in Ram} \infty \,r_\alpha))\}.$$
Imposing the condition $Res_{r_\alpha}(\eta)=0$ on  elements of $G_\Sigma^{\prime}$ we obtain $G_\Sigma$.
This proves the Proposition.
$\blacksquare$

\begin{rmk}\label{sympl pairing}
The Proposition gives also an alternative explanation to the fact that the hypercohomology is a symplectic vector space. 
\end{rmk}

The Proposition \ref{reduction} gives rise to a formal symplectic vector space and formal connection only, since in the framework of Tate spaces we understand varieties and bundles only as functors on Artin algebras. 

Recall that we have an embedding of the germ $\BB_1$ (in fact even a bigger space $\BB_0$) into $\H_\Sigma$ derived from the Theorem \ref{psi iso}. Let us compare these two embeddings.

\begin{prp}\label{comparison}
The formal Lagrangian submanifold of $\H_\Sigma$ defined in this subsection is a formal completion at zero 
of the germ of the Lagrangian submanifold defined in Section 5.2.

\end{prp}

{\it Proof.} Recall the RHS of the formula \eqref{iso psi} which consists of two summands. Recall also the global description of $\H_\Sigma$ given in the formula (\ref{represent H}). Consider a complex manifold with boundary obtained from $\Sigma$ by removing small open discs $\D_{r_\alpha}, \alpha\in Ram$ which contain ramification points $r_\alpha$. We have the corresponding vector space of holomorphic $1$-forms, which are analytic  on the boundary, and which have trivial integrals over each boundary component. We can use this space of $1$-forms in order to obtain an presentation for $\H_\Sigma$ analogous to (\ref{represent H}):
\begin{equation}\label{new quotient}
\H_{\Sigma}=\{\eta\in \Gamma(\Sigma-\cup_\alpha \D_{r_\alpha}, \Omega^1_\Sigma(D_\Sigma))|\int_{\partial\overline{\D}_{r_\alpha}}\eta=0\}/d\Gamma(\Sigma-\cup_\alpha \D_{r_\alpha},\OO_\Sigma(-D_\Sigma)).
\end{equation}

Using the foliation $\FF$ we can identify as abstract complex manifolds with boundary 
$$\Sigma-\cup_\alpha \D_{r_\alpha}\simeq \Sigma^{\prime}-\cup_\alpha \D_{r_\alpha}^{\prime},$$
for any $\Sigma^{\prime}$ which is sufficiently close to $\Sigma$ and the corresponding $\D_{r_\alpha}^{\prime}$. Hence the RHS of the formula (\ref{new quotient}) is 
invariant with respect to the deformations of $\Sigma$. As a result  we obtain a linear flat symplectic connection described in terms of spectral curves with small discs removed.

Notice the RHS of the formulas \eqref{new quotient} and \eqref{represent H} are the only cohomologies in degree $1$ of a $2$-term complexes. There is a natural quasi-isomorphism of these complexes induced by the embedding $\Sigma-\cup_\alpha \D_{r_\alpha}\to \Sigma-\cup_\alpha \{r_\alpha\}$.
Moreover the isomorphsim $\psi$ from formula \eqref{iso psi} can be generalized to the description of $\H_\Sigma$ given in the formula \eqref{new quotient}.

This ensures equivalence of both descriptions of $\H_\Sigma$, which in turn ensures equivalence of two descriptions of the affine symplectic connection on $\H_\Sigma$.

Notice that the Hamiltonian reduction in the framework of Tate spaces should be understood
in terms of functors on Artin algebras (see Remark \ref{Tate}). Then we see that the embedding of the moduli space $\BB_1$ obtained from the  Proposition \ref{reduction} coincides (as a functor on Artin algebras)
with the description given in Section 5. This proves the statement. $\blacksquare$

\subsection{Quantum Airy structures and quantum Hamiltonian reduction of the wave function}

In the case $Pol=\emptyset$
let us choose a Lagrangian complement 
${\mathcal V}_\Sigma\subset H^1(\Sigma,\C)$ to the Lagrangian subspace given by the first term of the Hodge filtration
$F^1H^1_{dR}(H^1(\Sigma,\C)):=H^{1,0}(\Sigma)=T_\Sigma \BB\subset H^1(\Sigma,\C)={\mathcal H}_\Sigma$ (a splitting of the Hodge filtration is equivalent to a choice of 
Bergman kernel in the traditional language of topological recursion). 

Let us  consider the infinite-dimensional vector space 
${V}_\Sigma \subset G_\Sigma$ consisting of holomorphic $1$-forms $\gamma\in \Omega^1(\Sigma-\cup_{\alpha\in Ram}\{r_\alpha\})$, 
which are meromorphic on $\Sigma$, have trivial residue at all $r_\alpha$, and such that $[\gamma]\in {\mathcal V}_\Sigma$. 
It is easy to see that ${V}_\Sigma$ is a Lagrangian complement to the tangent space to the Lagrangian submanifold ${L}=L_{Airy}^{Ram}\subset W=W^{Ram}_{Airy}$.

The story is similar in the case $Pol\ne \emptyset$. In this case we choose a transversal Lagrangian complement ${\mathcal V}_\Sigma$ in $\H_\Sigma$ 
to the Lagrangian subspace  $H^0(\Sigma,\Omega^1_\Sigma(D_\Sigma))$. Its preimage $V_\Sigma:=(G_\Sigma\epi \H_\Sigma)^{-1}{\mathcal V}_\Sigma$ is a Lagrangian complement to $T_0 L$. 
Thus, we have in general the following result.

\begin{cor} The Lagrangian complement $V_\Sigma$ carries both classical and quantum Airy structures.

\end{cor}

It follows that we have the wave function $\psi=\exp(\sum_{g\ge 0}S_g\hbar^{g-1})$ associated with the cyclic $DQ$-module corresponding to ${L}$ and the Lagrangian complement $V_\Sigma$. If we don't choose the Lagrangian complement then we still have a cyclic $DQ$-module with a cyclic vector, i.e. a left ideal in the completed Weyl algebra. The choice of Lagrangian complement allows us to encode the left ideal in terms of the wave function.

\begin{rmk} Notice that in the global version of the story the space $V_{\Sigma,odd}$ is the subspace of the previously defined ``local'' space $V_{odd}$ which consists of such meromorphic forms $\gamma$ on $\Sigma-\cup_\alpha\{r_\alpha\}$ that
$\gamma+\sigma^{\ast}(\gamma)$ has no poles at each $r_\alpha$. Here $\sigma$ is the local involution from Section 3.1. The subspace $V_{\Sigma,odd}\subset V_\Sigma$ carries a quantum Airy substructure.

\end{rmk}

Let us explain how the Hamiltonian reduction gives us the cyclic $DQ$-module ${\mathcal E}_\Sigma$ with a cyclic vector quantizing the Lagrangian germ $\BB_\Sigma\subset {\mathcal H}_\Sigma$.  
This is a part of the following more general story.

Suppose we are given an affine symplectic vector space $M$ an affine coisotropic subspace $G\subset M$,
a point $x\in G$ and a formal germ of Lagrangian submanifold $L_x$ at $x$ satisfying the conditions 1) or 2) from Section 7.1. Let us work for simplicity over the field of complex numbers $\C$, although it is not necessary.

Notice that with the point $x$ we can associate the completed Weyl algebra $Weyl_x$ of the tangent space $T_xM$. It is a topologically free algebra over $\C[[\hbar]]$. Suppose we are given a cyclic module $\mathcal E$ over $Weyl_x$ with the cyclic vector $e\in {\mathcal E}$. Assume that modulo $\hbar$ we have ${\mathcal E}\simeq \OO(L_x)$ and $e=1\in \OO(L_x)$. By Proposition 7.1.1 we have a symplectic affine space $M^{\prime}=G/G^\perp$ and the Hamiltonian reduction $L_x^{\prime}\subset M^{\prime}$ of the germ $L_x$. We also have the completed Weyl algebra
$Weyl_x^{\prime}$ associated with $M^{\prime}$. Then $G$ is naturally embedded an affine Lagrangian subspace in the symplectic affine space $M^{\prime}\times \overline{M}$, where $\overline{M}$ denote $M$ with the opposite symplectic form. By general theory we have a {\bf canonical} cyclic DQ modules ${\mathcal E}_G^{can}$ with the canonical cyclic vector $e_G^{can}$ quantizing $G$.

We treat ${\mathcal E}_G^{can}$ as a $Weyl_x^{\prime}-Weyl_x$-bimodule.
\begin{defn} Quantum Hamiltonian reduction of $({\mathcal E},e)$ is the $Weyl_x^{\prime}$-module
${\mathcal E}^{\prime}= {\mathcal E}_G^{can}\otimes {\mathcal E}$ endowed with the cyclic vector $e^{\prime}=e_G^{can}\otimes e$.

\end{defn}
Clearly modulo $\hbar$ we obtain the Hamiltonian reduction $L_x^{\prime}$ of $L_x$.

The above construction works for our infinite-dimensional affine spaces and germs of formal Lagrangian manifolds. As a result we obtain a cyclic $DQ$-module ${\mathcal E}_\Sigma$ associated with the spectral curve $\Sigma$. In general a choice of Lagrangian complement to $T_x(L_x)$ allows us to translate the story from the language o cyclic $DQ$-modules to the language of wave functions. 
In the case of spectral curves this translation looks as follows.

Recall that we have a ``local'' symplectic vector space $W^{Ram}$ together with a germ of formal Lagrangian submanifold $L\subset W^{Ram}$, both depending only\footnote{Here we choose a canonical coordinate near each ramification point, determined a priori only up to multiplication by a third root of 1.} on the set $Ram$ of ramification points. We have also defined a coisotropic vector space
$G_\Sigma\subset W^{Ram}$ depending on $\Sigma$.

Then the quantization procedure give us  a cyclic $DQ$-module ${\mathcal E}_\Sigma$ with the cyclic vector $\psi_L$.
Let us split coordinates in each of the symplectic spaces, so that the dual coordinates in $W^{Ram}$ will be denoted by $(q,p)$ and in ${\mathcal H}_\Sigma$ by $(q^\prime,p^\prime)$ in such a way that $\psi_L=\psi_L(q)$ is the wave function corresponding to the quantized Airy structure.

We have  the wave function $\psi_{G_\Sigma}=\psi_{G_\Sigma}(q,q^\prime)=e^{Q_2(q,q^\prime)/\hbar}$ corresponding
to the coistoropic subspace $G_\Sigma$, 
where $Q_2$ is the quadratic form such that $graph(dQ_2)=G_\Sigma$.
Then the Hamiltonian reduction of $L\subset W^{Ram}$ to $\widehat{\BB}_\Sigma\subset {\mathcal H}_\Sigma$ at the level of wave functions becomes the integral operator 
$$\psi_L\mapsto \psi_{\widehat{\BB}}:=\int e^{Q_2(q,q^\prime)/\hbar}\psi_L(q)dq$$
understood as a formal integration with respect to the Gaussian measure.

\subsection{Holomorphic Anomaly Equation}\label{HAE section}

Let us recall the framework and notation from Section 5.1. Thus we have a tuple $(X,\H,\langle\bullet,\bullet\rangle,\nabla,\phi)$ consisting of an analytic manifold $X$, analytic vector bundle $\H$ over $X$, linear flat connection $\nabla$ on $\H$, $\H$-valued $\nabla$-flat $1$-form $\phi$, and a $\nabla$-covariant symplectic structure $\langle\bullet,\bullet\rangle$ on $\H$. 

Assume that for any point $x\in X$ and assume that the embedding of the germ of $X_x$ of $X$ at $x$ into $\H_{x}$ defined in Section 5.1 is Lagrangian (i.e. the image is a germ $L_x$ of Lagrangian submanifold in $\H_{x}$). Notice that the affine flat connection $\nabla+\phi$ induces a flat connection $\widehat{\nabla}$ on the vector bundle over $X$ with the fiber $Weyl_x$ (completed at zero  Weyl algebra of $\H_x$).  This comes from the fact that the group of affine symplectic transformations acts on the non-completed Weyl algebra.

\begin{defn} A {\bf quantization} of the tuple $(X,\H,\langle\bullet,\bullet\rangle,\nabla,\phi)$ is given by the a $\widehat{\nabla}$-covariantly constant family of closed left ideals $J_x\subset Weyl_x$ such that the cyclic module $Weyl_x/J_x$ is a quantization the Lagrangian germ $L_x:=\phi(X_x)$ (i.e. modulo $\hbar$ it gives the algebra of function $\OO(L_x)$, and the cyclic vector $1\in Weyl_x/J_x$ corresponds to $1\in \OO(L_x)$).
\label{def:quantization}

\end{defn}

Let us  choose a $\nabla$-covariantly constant subspace $V_x$ which is transversal to $L_x$ for all $x\in X$. Then the condition that $J_x$ form a covariantly constant family can be translated into a system of differential equations on the components $(S_{g,n})_x\in Sym^n(V_x)$ of the cyclic vector $\psi_x=\exp(\sum_{g,n} \hbar^{g-1}(S_{g,n})_x)$. The summation is over pairs $(g,n)$ with either $g=0,n\ge 3$ or $g\ge 1,n\ge 1$. This system of equations is called (quantum) {\it Holomorphic anomaly equation} (HAE for short). 

\begin{rmk} Originally HAE was formulated in the context of \emph{compact} Calabi-Yau 3-folds.
 In this case  $X$ (probably should be better denoted  by $\mathcal L$, as it is locally a Lagrangian cone)   is a $\C^*$-bundle over  a connected component of the moduli space $\mathcal M$ of projective algebraic 3-dimensional varieties $Y$ with $H^1(Y,\C)=0$ and trivial canonical class. Space $X={\mathcal L}$ parametrizes pairs $(Y,\Omega^{3,0}_Y)$ where   $\Omega^{3,0}_Y\in\Gamma(Y,\Omega^3_Y)-\{0\}$ a non-zero holomorphic volume form, and it is embedded locally as  holomorphic Lagrangian cone into sympletic vector space $H^3(Y,\C)$ via
  \[(Y,\Omega^{3,0}_Y)\mapsto [\Omega^{3,0}_Y]\in H^3(Y,\C)\,.\]
 It was predicted  in [BCOV]\footnote{In [BCOV] at each point $x=[(Y,\Omega^{3,0}_Y)]\in {\mathcal L}$  the Lagrangian complement to the tangent space $T_x {\mathcal L}=H^{3,0}\oplus H^{2,1}\subset H^3(Y,\C)$ was secretly chosen as the complex conjugate $H^{1,2}\oplus H^{0,3}$. E.~Witten in [Wi] explained the invariant meaning of   HAE as a quantization.} that this Lagrangian cone has a {\bf canonical} quantization, equivariant with respect to $
 \C^*$-action acting by rescaling
 \[ [\Omega^{3,0}_Y]\mapsto \lambda [\Omega^{3,0}_Y],\,\,\,\hbar\mapsto \lambda^2\hbar,\,\,\,\lambda\in \C^*\,.\]
 The rigorous construction of this quantization was proposed by K.~Costello in [Cos]. 
 Our framework here is more general, and is applicable e.g. to the moduli spaces of non-compact Calabi-Yau 3-folds, as in [KoSo].
\end{rmk}

Let us now set $X=\BB_1$. We have constructed the symplectic vector bundle $\H$ with the fiber $\H_\Sigma$ previously, same for the $1$-form $\phi$. It follows from considerations of the  Sections 7.1, 7.2 that
we have a tuple  $(X,\H,\langle\bullet,\bullet\rangle,\nabla,\phi)$ and its quantization.

Notice that the larger space $\BB_0\supset \BB_1$ does not depend on the choice of foliation $\FF$. Its local Lagrangian embeddings into 
 symplectic affine spaces depend only on finite jets of  $\FF$ at $D_\Sigma$.

\begin{que} Can this quantization be extended to $\BB_0$?

\end{que}

We expect that the answer is positive. Notice that since $\BB_0$ does not depend on the foliation then the corollary would be symplectic invariance in topological recursion which does not have so far a conceptual explanation (see [EynOr2]).

\begin{que} If the answer to the previous question is positive, does the quantization for $\BB_0$ depend
only on the collection  of $(2\cdot k_\beta-1)$-jets of $\FF$ at $D$ near points $p_\beta,  {\beta\in Pol}$ (see Section 5.2)?

\end{que}
In the next section we will ask a stronger question \ref{Question:  symplectic invariance}, positive answer to  which implies independence of the quantization of the choice of finite jets  of $\FF$.

Notice that although consideration of the Sections 7.1, 7.2 ensures that we have formal germs and quantizations only, the condition to be covariantly constant guarantees that everything is analytic (e.g. $S_{g,n}$ are analytic germs at each point $x\in X$).

\begin{rmk} Our definition \ref{def:quantization} of quantization can be  improved further. Namely, for any covariantly constant family of left ideals $J_x\subset Weyl_x$
 one can look for cyclic vectors (after choosing a Lagrangian complement to $T_x X \subset \H_x$) of the form $\exp(\sum_{g,n} \hbar^{g-1}(S_{g,n})_x)$ where the constraint $n>0$ is omitted. The ambiguity for the lift is a torsor over the abelian group
  $\hbar^{-1}\C[[\hbar]]$ (or, equivalently, the multiplicative group of formal expressions $\exp(\sum_{g\ge 0} (S_{g,0})_x \hbar^{g-1})$ where $(S_{0,0})_x,(S_{1,0})_x,\dots\in \C$ is a sequence of constants).  
  If we ignore terms $(S_{0,0})_x$ and $(S_{1,1})_x$, the rest does not depend on the choice of Lagrangian splitting, and one can speak about covariantly constant  formal wave functions modulo the factor $\exp(c_0\hbar^{-1}+c_1)$. Therefore, for each $g\ge 2$ we obtain a {\emph canonical} $\C$-torsor with flat connection over $X$. The monodromy gives for every $g\ge 2$ a homomorphism $\pi_1(X)\to \C$.  Combining them together we obtain an ``obstruction'' class in
   $H^1(X,\hbar\,\C[[\hbar]])$. \end{rmk}
\begin{defn}
We call an {\bf enhanced} quantization a choice of flat sections of the above torsors for all $g\ge 2$
\end{defn}
The enhanced quantization  exists if and only if the obstruction class 
 vanishes. In the  case of a conical quantization, like for compact Calabi-Yau 3-folds, there is a canonical conical flat trivialization of all above  torsors for $g\ge 2$, hence a canonical enhancement.

\begin{que} Does our quantization of the moduli space of spectral curves $\BB_1$ extend (in some canonical way) to an enhanced quantization?
\end{que}

\begin{rmk} The datum of a (non-enhanced) quantization corresponding to $g=0$ is canonically determined by the tuple $(X,\H,\langle\bullet,\bullet\rangle,\nabla,\phi)$, hence carry no more additional information. It is well-known that $g=1$ part of the quantization is equivalent to a choice of flat connection on virtual line bundle $(\Omega^{\dim X}_X)^{\otimes (1/2)}$ (well-defined locally up to involution $\times(-1)$). Equivalently, by taking the square, one gets a flat connection on the honest line bundle $\Omega^{\dim X}_X$, but it is still preferable to speak about its square root as the wave function is invariantly a formal section of  $(\Omega^{\dim X}_X)^{\otimes (1/2)}$. In the example of compact Calabi-Yau 3-folds the monodromy of the induced flat connection on  $\Omega^{\dim X}_X$ is non-trivial and is equal to multiplication by $ \exp(2\pi i\chi(Y)/24)$ along the loop in $\C^*$ in the conical direction. One may ask the question whether there exists a canonical flat section well-defined up to $\sqrt[24]{1}$ (as an analog of enhanced quantization for $g=1$).
\end{rmk}

\section{Semi-affine Lagrangian embeddings}

The goal of this section  is to define a generalization of local Lagrangian embedding 
\[\text{germ at }[\Sigma]\text{ of } \BB_0\hookrightarrow  \text{ germ at }0 \text{ of } \H_\Sigma\]
which  {\it does not} use  foliation $\FF$, and formulate  questions about the  irrelevance of the choice of $\FF$ for  the canonical quantization of locally embedded $\BB_1$ discussed in the previous section.

\subsection{Semi-affine symplectic structures}
\begin{defn} A {\bf semi-affine structure} on a complex  analytic symplectic  manifold $(M,\omega=\omega_M)$ 
  is given by a coisotropic foliation $\FF$ and an affine structure on leaves of $\FF$ 
 such that near each point there exists a system of local coordinates 
 \begin{equation}\label{coordinates}
 x_1,\dots,x_n,y^1,\dots,y^n;z_1,\dots,z_{2m},\,\,\omega=\sum_{i=1}^n dx_i\wedge dy^i+\sum_{j=1}^m dz_j\wedge dz_{j+m}\,,\end{equation}
 the leaves of $\FF$ are given by $x_\bullet=const$ and the affine structure on leaves is given by coordinates $y^\bullet,z_\bullet$.
 \end{defn}
 
 There are two particular cases: 
 \begin{itemize}
 \item[(1)] $n=0$, then $M$ is the only leaf of the foliation and carries an affine symplectic structure, 
 \item[(2)] $m=0$, then $M$  is endowed with a Lagrangian foliation, and the affine structure on leaves is automatic (e.g. $M=T^*X$ is a cotangent bundle endowed with the cotangent foliation).
 \end{itemize}
In general, local coordinates in the above definition are determined up to the action of the infinitesimal group of transformations whose Lie algebra $\g_{m,n}$ is given by Poisson brackets with analytic
functions of the form
\[F^{(1)}(x_\bullet)+\sum_{i=1}^n F^{(2)}_i(x_\bullet) y^i+\sum_{j=1}^m F^{(3)}_j(x_\bullet) z_j+\sum_{j_1=1}^m\sum_{j_2=1}^m  F^{(4)}_{j_1 j_2}(x_\bullet) z_{j_1}z_{j_2}\,.\]
In what follows we will use an alternative description of germs of semi-affine symplectic structures:
\begin{lmm} \label{lemma: Poisson fibrations}
A semi-affine structure on a {\emph{germ}} at $m\in M$ of symplectic manifold $M$ is equivalent to a Poisson fibration 
\[\pi:\text{germ of }M \text{ at } m\to \text{ germ of Possion manifold }P\text{ at }p=\pi(m)\]
such that $P$ is {\emph {regular}} (i.e. the rank of Poisson structure is constant),
 and symplectic leaves of $P$ near $p$ are endowed with affine structures compatible with the natural symplectic structures (i.e. the symplectic structure is constant in affine coordinates).
\end{lmm}
{\it Proof.} Our assumptions ({\it without} the use of affine structures) imply that one can choose local coordinates on $M$ at $m$ as in \eqref{coordinates} such that $\pi:M\to P$ is the coordinate projection 
\[\pi: (x_\bullet,y^\bullet,z_\bullet)\mapsto (x_\bullet,z_\bullet)\,.\]
The natural map $\pi_1$ from the germ  of $P$ at $p$ to the ``center'' of $P$ (the germ of the space of symplectic  leaves in $P$) is given by the coordinate
projection
\[\pi_1:(x_\bullet,z_\bullet)\mapsto (z_\bullet)\,.\]
The Lie algebra of germs of vector field at $m$ preserving $\omega_M$ and acting on $\pi:M\to P$ is given by Poisson brackets with analytic functions of the form
\begin{equation}\label{symmetries of Poisson map}
\sum_i f_i(x_\bullet) y^i+G(x_\bullet,z_\bullet)\,.\end{equation}
 In particular, it maps {\it epimorphically} to the Lie algebra of germs of Poisson vector fields on $P$ at $p$. The latter contains as a Lie subalgebra germs of symplectic vector fields in coordinates $z_\bullet$ depending arbitrarily on center coordinates $x_\bullet$. Therefore, we can assume that the affine symplectic structures on symplectic leaves of $P$ near $p$ (see the Lemma) are the {\it standard} ones in coordinates $z_\bullet$.  It is straightforward to see that Poisson brackets with functions from \eqref{symmetries of Poisson map} act as affine transformations in $z_\bullet$ if and only if $G=G(x_\bullet,z_\bullet)$ is at most quadratic in $z_\bullet$. Hence we obtain exactly the Lie algebra $\g_{m,n}$ of local symmetries of a semi-affine symplectic structure, where $n$ is number of $x_\bullet$-coordinates, and $2m$ is number of $z_\bullet$-coordinates.
$\blacksquare$

  \begin{defn}\label{def:flat semiaffine} 
  A {\bf flat family of semi-affine  symplectic germs} on a complex manifold $X$ is a holomorphically depending on $x\in X$  family of germs of semi-affine symplectic manifolds $(M_x,m_x)$ (with symplectic forms $\omega_x$, and germs of Poisson quotients $(P_x,p_x)$ as in 
  Lemma \ref{lemma: Poisson fibrations}), endowed with a non-linear flat connection preserving semi-affine symplectic structures, but not necessary based points $(m_x)_{x\in X}$. 
  \end{defn}
Such a family gives  for any $x_0\in X$ a map from the germ of $X$ at $x_0$ to the germ $M_{x_0}$ at $m_{x_0}$ given by $x\mapsto m_x$  for $x$ sufficiently close to $x_0$, where we idenitify locally fibers with $M_{x_0}$ via the flat connection.
\begin{defn} A flat family of  semi-affine  symplectic germs is said to give a {\bf flat family of local Lagrangian embeddings} if the above-defined maps of germs gives a Lagrangian embedding 
\[\text{germ of }X\text{ at } x_0\hookrightarrow  M_{x_0}\]
for any $x_0\in X$.
\end{defn}
This definition generalizes tuples $(X,\H,\langle\bullet,\bullet\rangle,\nabla,\phi)$ from Sections 5.1, \ref{HAE section} which give local Lagrangian embeddings into affine symplectic spaces. One can also generalize  straightforwardly   the previous definitions to the case when germs of semi-affine manifolds are replaced by {\bf formal} germs. 

\begin{rmk} Omitting semi-affine symplectic structures in the definition \ref{def:flat semiaffine}
we obtain a more basic notion of a {\bf flat family of germs} (or formal germs)  on a manifold $X$.  For example, if $X$ is endowed with a foliation $\FF$, then on $X$ we have a flat family of germs of the ``space of leaves'' of $\FF$.
\end{rmk}

There is an embedding of the Lie algebra $\g_{m,n}$ into the Lie algebra of $\hbar$-differential operators in $(x_i)_{1\le i \le n},(z_j)_{1\le j\le m}$ under which $y^j\mapsto \hbar\partial/\partial x_j, z_{l+m}\mapsto \hbar  \partial/\partial z_l, 1\le j\le n, 1\le l\le m$,
and the rest of variables are mapped to themselves. This give rise to a {\it canonical} completed Weyl algebra associated with a formal germ of a semi-affine symplectic manifold.
  
  For any flat family of  formal semi-affine  symplectic germs we obtain a bundle of completed Weyl algebras endowed with a flat connection.
  Similarly to Section \ref{HAE section} we define {\bf quantization} and {\bf enhanced quantization} of  a flat family of  formal local  Lagrangian embeddings.

    A particular example is the one of a global Lagrangian submanifold $X=L$ in a cotangent bundle $M=T^*Y$ endowed with the natural semi-affine structure given by the foliation by cotangent fibers. Taking formal germs of $M$ at points of $L$ we obtain a flat family of local Lagrangian embeddings. The quantization of this structure in our sense is the same as a cyclic $DQ$-module on $T^*Y$ with the classical limit $\mathcal{O}_L$ endowed with cyclic vector $1$, i.e. a WKB-module.
    
    \begin{defn} For a given symplectic manifold $(M,\omega)$ and two semi-affine structures given by coisotropic foliations $\FF,\FF'$ with 
     affine structures on leaves, we say that $\FF$ is {\bf weaker} than $\FF'$ if leaves of $\FF$ are contained in leaves of $\FF'$ and the inclusion of any leaf is an affine map.
    \end{defn}
    It is easy to see that in such a case the flat bundles of completed Weyl algebras associated with two  semi-affine structures are  \emph{canonically} isomorphic. For example, the Weyl algebra associated with the cotangent space $T^*V$ to an affine space  $V$ with the natural semi-affine structure of the cotangent bundle is canonically isomorphic to the Weyl algebra associated with $T^*V$ considered as  an affine symplectic space. Same is true for completed Weyl algebras.

\subsection{Canonical semi-affine structure on $\BB_0$}
The setup in this section is the following: 
$P$ is a complex analytic Poisson surface, with open dense symplectic leaf $(S,\omega=\omega_S)$ and such that the effective divisor $D$ of zeroes of Poisson structure is smooth. We denote as usual by $\BB_0$ the moduli space of $D$-transversal smooth connected compact curves $\Sigma\subset P$. 
Our goal here is to construct a flat family  of Lagrangian embeddings into  formal semi-affine symplectic
germs 
\[\text{germ }\BB_{[\Sigma]}\text{ of }\BB_0\text{ at }[\Sigma]\hookrightarrow \text{ germ }M_{[\Sigma]}\]
  as in the previous section. 
  

	For a given $D$-transversal spectral curve $\Sigma$ we denote as before by $\{p_\beta\}_{\beta\in Pol}$ the set $D\cap \Sigma$, by $(k_\beta\ge 1)_{\beta\in Pol}$ the corresponding multiplicities, and by $D_\Sigma$ the effective divisor $\sum_\beta k_\beta p_\beta $ on $\Sigma$.
	 {\bf We assume that 
	{${{Pol\ne \emptyset}}$}}. 

	 In order to define a formal germ of symplectic manifold $M_{[\Sigma]}$ and a  Lagrangian embedding of germ $B_{[\Sigma]}$ 
	 , we proceed as follows.
	 \begin{defn} \label{ML def}For $[\Sigma]\in \BB_0$ define a {\bf collection of movable loops}  as a tuple 
	 consisting of 
	 \begin{itemize}
	  \item[(1)] a {\it non-empty} finite collection of disjoint closed discs $(\D_i)_{i\in I}$ in $\Sigma-\{p_\beta\}_{\beta\in Pol}$ with real-analytic boundaries,
	  \item[(2)] a collection of germs of holomorphic 1-forms $\nu_i$ on $S$ defined in a neighborhood of  of $\D_i\subset S$ such that
	   $d\nu_i=\omega$ and $\nu_i$  is everywhere non-vanishing 1-form near $\partial \D_i$ with the kernel {\it transversal}
	   to $\partial \D_i$.
	   \end{itemize}
	 \end{defn}
	 
We denote by $ML(\Sigma)$ the set of such collections. We endow $ML(\Sigma)$ with the partial order:
 \[\sigma_1\prec \sigma_2\,,\]
  if $\sigma_2$ is obtained from $\sigma_1$ by adding  new discs and new germs of 1-forms satisfying the above conditions.

Obviously, the nerve of this poset is {\it contractible}. The reason for that is the simple fact that for any finite collection of closed discs (maybe intersecting each other) in the non-compact (since $Pol\ne \emptyset$) surface 
  $\Sigma-\{p_\beta\}_{\beta\in Pol}$, there exists another disc in  $\Sigma-\{p_\beta\}_{\beta\in Pol}$ disjoint from our collection.

	For every collection of movable loops $\sigma\in ML(\Sigma)$
	we define {\it a formal germ of infinite-dimensional manifold $G_\sigma$} by considering  formal deformations $\Sigma_0'$ of the  complex curve 
	\[\Sigma_0:=\Sigma-\sqcup_i int(\D_i)\subset P\]
	with real-analytic  boundaries, such that the boundary is locally a shift 
	 of $\partial \D_i$ along the foliation $\text{Ker}(\nu_i )$, and the integral of $\nu_i$ over the corresponding boundary component vanishes\footnote{The germ $G_\sigma$ can be thought of as a non-linear analog of the coisotropic subspace $G_\Sigma$ defined in \eqref{coisotropic space}.}. 
	  
	The germ $G_\sigma$ is naturally embedded into the  infinite-dimensional  vector space 
\[W_\sigma=\Big\{(\eta_i)_{i\in I}\,|\eta_i\text{ is }\C\text{-valued real analytic 1-form on }\partial\D_i,\,\text{and}\,\,\forall i\,\, \, \int_{\partial \D_i}\eta_i=0 \Big\}\]
endowed with the usual skew-symmetric pairing:
\[\langle (\eta_i^{(1)})_{i\in I},(\eta_i^{(2)})_{i\in I}\rangle:= \sum_i \int_{\partial \D_i} f_i^{(1)} \eta^{(2)}_i,\quad df_i^{(1)}=\eta_i^{(1)}\,.\]
  Therefore, $G_\sigma$ carries a {closed} 2-form $\omega_{\sigma}$, which is the pullback of the translation-invariant closed 2-form  on $W_\sigma$
  associated with the above pairing.
  
  Let us describe the tangent space to $G_\sigma$ and its skew-orthogonal.
The tangent space to $G_\sigma$ at any point $\Sigma_0'$ close to $\Sigma_0$ \footnote{Since we speak about formal germs, we mean
that  $\Sigma_0'$ is a deformation of $\Sigma_0$ over the spectrum of a local Artin algebra over $\C$.} 
is given by formula (compare with the proof of Proposition \ref{prop with discs})
\begin{equation}\label{formula big tangent}T_{\Sigma_0'}G_\sigma=
\Big\{\eta \in \Gamma({\Sigma'_0},\Omega^1_{{\Sigma'_0}}(D_{{\Sigma'_0}}))\,| \int_{\partial_i {\Sigma'_0}} \eta=0,\,\,\forall i\in I
\Big\}\,.\end{equation}
It is naturally embedded into $W_\sigma$, using foliations $Ker (\nu_i)$.
The  skew orthogonal to $T_{\Sigma_0'}G_\sigma$ is the subspace 
\begin{equation}\label{orthogonal}
(T_{\Sigma_0'}G_\sigma)^\perp= d\Gamma({\Sigma'_0},\mathcal{O}(-D_{{\Sigma'_0}}))\subset T_{\Sigma_0'}G_\sigma.
\end{equation}
The quotient space $T_{\Sigma_0'}G_\sigma/(T_{\Sigma_0'}G_\sigma)^\perp$ is a finite-dimensional symplectic vector space which we denote by $\H_{\Sigma'_0,D_{{\Sigma'_0}}}$. This space depends only on $\Sigma_0'$ (which is considered as an abstract complex curve with real-analytic boundary) and on an effective divisor $D_{{\Sigma'_0}}$ with the support in the interior of $\Sigma_0'$.

Suppose we have a germ of a  complex manifold endowed with a closed 2-form of constant rank. It is well-known that the subbundle of the tangent bundle given by the kernel of the form defines a foliation. The germ of the space of leaves of this foliation is naturally a germ of symplectic manifold. We apply these general considerations
 to the infinite-dimensional manifold  $(G_\sigma,\omega_\sigma)$ and obtain a well-defined finite-dimensional formal symplectic germ, which we denote by $M_\sigma$ (abusing the notation we will denote the symplectic form on $M_\sigma$ also by  $\omega_\sigma$ as  the 2-form  on $G_\sigma$ is its pullback).
 
 Let us fix $\Sigma$. Now we are going to identify {\it canonically} all symplectic germs $(M_\sigma,\omega_\sigma)$ for $\sigma \in ML(\Sigma)$.
 Recall the partial order $\prec$ on $ML(\Sigma)$.
 Obviously, if $\sigma_1\prec \sigma_2$,
	  we have a natural embedding
	  \[i_{\sigma_1,\sigma_2}:G_{\sigma_1}\hookrightarrow G_{\sigma_2}\]
	  and in this way we get  a functor from  the poset $ML(\Sigma)$ to the category of infinite-dimensional formal germs and embeddings.
	  One can see immediately that $i^*\omega_{\sigma_2}=\omega_{\sigma_1}$, and moreover $i_{\sigma_1,\sigma_2}$ induces a 
	  symplectomorphism of formal germs $M_{\sigma_1}\stackrel{\sim}{\to}M_{\sigma_2}$. The connectedness and 1-connectedness of the nerve
  of poset $ML(\Sigma)$ implies that we identify all formal germs $M_\sigma$ with each other.
    \begin{defn} For any point $[\Sigma]\in \BB_0$ define the associated symplectic formal germ $M_{[\Sigma]}$ as the result of the identifications of formal germs $M_\sigma$ for $\sigma\in ML(\Sigma)$ described above.
  \end{defn}

 It is obvious from our description that  formal germs $(M_{[\Sigma]})_{[\Sigma]\in \BB_0}$ form a {\it flat family} of formal germs, and symplectic structures on fibers are covariantly constant.
 Also it is straightforward to see that the flat family of germs $([\BB_\Sigma])_{[\Sigma]\in \BB_0}$ has natural covariantly constant Lagrangian embeddings into $(M_{[\Sigma]})_{[\Sigma]\in \BB_0}$.
 
 Finally let us describe a natural covariantly constant family of semi-affine structures on $M_{[\Sigma]}$.
 For given curve $\Sigma$ and $\sigma\in ML(\Sigma)$, let us consider a  subbundle of finite codimension of the tangent bundle to  $G_\sigma$. The subbundle is given by first order
 deformations such that the corresponding 1-form $\eta$ with poles at $D_{\Sigma_0'}$ on $\Sigma_0'$ (which is a small deformation of $\Sigma_0$ as a curve with boundary) is \emph{holomorphic}, i.e. has no poles. It is easy to see that this distribution is {\it integrable}, i.e. it defines a foliation $\FF_\sigma$ on $G_\sigma$.  This follows immediately 
 from the geometric interpretation of the above constraint on first order deformations. More precisely,
it means that intersection points
     $p_\beta$ of the moving curve with the divisor $D$ stay constant, and moreover $(k_\beta-1)$-jet of $\Sigma_0$ at $p_\beta$ stays constant.
      Also it is clear from this description that leaves of  $\FF_\sigma$ contain leaves of the foliation $Ker (\omega_\sigma)$.
       Hence $\FF_\sigma$ induces a foliation on $M_{[\Sigma]}$, and it is easy to see that it does not depend on $\sigma$ (since the maps $i_{\sigma_1,\sigma_2}$
        identify induced foliations).

        Let us denote by $\FF_{M_{[\Sigma]}}^{co}$ the resulting foliation on $M_{[\Sigma]}$. Pick any $\sigma\in ML(\Sigma)$. Since  the restriction of $\omega_\sigma$ to $T_{\FF_\sigma}$ vanishes, it follows that
 $\FF_{M_{[\Sigma]}}^{co}$ is coisotropic. The codimension of leaves of $\FF_{M_{[\Sigma]}}^{co}$ coincides with the codimension of leaves of
  $\FF_\sigma$, and is equal to $\sum_\beta k_\beta-1$ (the term $-1$ comes from the condition that the sum of residues of 1-forms $\eta$ from \eqref{formula big tangent} at points $p_\beta'\in \Sigma_0^{\prime}$  vanish).
  
Notice that there is a 1-form $\xi$ on $\FF_{M_{[\Sigma]}}^{co}$. In order to see this, let us
replace in \eqref{orthogonal} the space of functions vanishing at $p_\beta, \beta\in Pol$ by the space of all functions. Then the quotient of the tangent space to $\FF_{M_{[\Sigma]}}^{co}$ by this bigger subspace is isomorphic to $H^1(\Sigma, \C)$. Natural (surjective) projection to the quotient gives rise to an $H^1(\Sigma,\C)$-valued 1-form $\xi$ along leaves of $\FF_{M_{[\Sigma]}}^{co}$.

\begin{lmm} The form $\xi$ is closed.
  \end{lmm}
  {\it Proof.}  We will check vanishing of $d\xi$ using the lift to $G_\sigma$ for any given $\sigma\in ML(\Sigma)$. In order to do this, we will construct an $H^1(\Sigma,\C)$-valued functional on pairs of points of $G_\sigma$ lying in the same  leaf of  $\FF_\sigma$, and satisfying the cocycle property. The construction is similar to what we had before. We interpret $H^1(\Sigma,\C)$
   as the space of $\C$-valued linear functionals on $H_1(\Sigma_0',\Z)$ vanishing on classes of boundary loops $\partial_i\Sigma_0'$ (here as before,
   $\Sigma_0'$ is a small deformation of $\Sigma_0$). With any two deformed surfaces $\Sigma_0',\Sigma_0''$ and closed 1-chains $\rho',\rho''$ on these surfaces representing  homology classes identified by the Gauss-Manin connection, we associate a complex number by integrating $\omega=\omega_S$ along 2-chain in $S$ with the boundary equal to $\rho''-\rho'$. This assignment depends only on homology classes, vanishes on boundary loops, satisfies the cocycle property from Remark 5.3.2,
   and gives the pullback of $\xi$ as the first derivative at the diagonal.
  $\blacksquare$
	
Thus we have a foliation on each $\FF_\sigma^{co}$ with fibers which are fibers of the map $\xi: T_{\FF_\sigma}^{co}\to H^1(\Sigma,\C)$. Tangent space to a fiber is isomorphic to 
$(T_{\FF_{\sigma}} ^{co})^\perp$, which is the skew-orthogonal with respect to the symplectic form $\omega_\sigma$.
Thus the skew-orthogonal subbundle $( T_{   \FF_{M_{[\Sigma]    } }^{co}})^\perp\subset T M_{[\Sigma]}$ 
  gives an isotropic subfoliation. The quotient bundle
  \[T_{   \FF_{M_{[\Sigma]    }}^{co} }/( T_{   \FF_{M_{[\Sigma]    } }^{co}})^\perp\]
  is naturally isomorphic to a flat symplectic bundle of rank $2g$ (where $g$ is the genus of $\Sigma$) with the fiber at the base point equal to  $H^1(\Sigma,\C)$. 
  
   This finishes the construction of the canonical flat family of Lagrangian embeddings of $\BB_0$ into a family of semi-affine symplectic germs. 
  
 Also, we see immediately  that in the case when the foliation $\FF$ on $P$  exists, the semi-affine structure coming from general considerations is automatically weaker (see Definition 8.1.6) than the semi-affine structure coming from $\FF$ (which is in fact a symplectic affine structure on $M_\Sigma$). According to the discussion after the Definition 8.1.6 in that case we obtain isomorphic flat bundles of completed Weyl algebras. It is naturally to inquire whether the corresponding quantizations depend on chosen foliations.
Hence we ask the following question.

\begin{que} \label{Question:  symplectic invariance} Is it true that the flat family of quantizations of locally Lagrangian embeddings $\BB_{[\Sigma]}\hookrightarrow M_\Sigma$
does not depend on the choice of foliation $\FF$? 
 \end{que}
Presumably the answer is positive. We expect that it is related to the deep result about symplectic invariance in TR proven by Eynard and Orantin in [EynOr2].
Also, the positive answer will imply (via  finite jet approximations of the foliation) that the canonical quantization can be constructed even in the case when there is no foliation $\FF$ at all. 

\section{Two more speculations}

\subsection{Wave function for an individual spectral curve}

One can generalize considerations of this paper to the case of immersed (rather than embedded) spectral curves which can 
have singularities. The moduli space $\BB$ of immersed singular spectral curves contains an open subspace $\BB_0$, the moduli space of   
immersed smooth spectral curves. 

Having in mind the above remark let us make a comment about a quantization of an individual spectral curve $\Sigma$,  assuming that the open symplectic leaf $S$ is isomorphic to $T^{\ast}C$, where $C$ is a compact curve.
If the answer to the Question \ref{Question:  symplectic invariance} is positive, then one can study the behavior of the corresponding wave function $\psi_{B_{[\Sigma]}}$ near the complement $\BB-\BB_0$  (discriminant locus)   in the bigger space $\BB$.

Let $T_x^{\ast}, x\in C$ denote the cotangent fiber at  $x\in C$. For generic collection of points $x_1,\dots,x_n\in C$  let us pick some pre-images ${y}_1,\dots,{y}_n\in \Sigma$. Then the abstract curve  
$\Sigma_{{y}_1,...,{y}_n}$ is 
defined as the disjoint union $T_{x_1}^{\ast}\sqcup...\sqcup T_{x_n}^{\ast}  \sqcup \Sigma$ with identified pairs of points of $\Sigma$ and $T_{x_i}, 1\le i\le n$ corresponding to the same point from the collection ${y}_1,\dots,{y}_n\in \Sigma$.  This singular curve is naturally immersed into the Poisson surface $P$, and hence gives a point in the discriminant $\BB-\BB_0$. W can try to deform it to an immersed smooth curve.
\begin{rmk}
The reason for picking {\it more} than 1 point in the collection $(x_i)_{1\le i\le n}$ is the following. If we want to deform an immersion of the singular curve $\Sigma_{{y}_1,...,{y}_n}$, we need (at first order deformation)  a meromorphic 1-form on $\Sigma$ with first order poles at ${y}_1,\dots,{y}_n\in \Sigma$. In the case $n=1$ such a form is automatically holomorphic, because the sum of residues vanishes. Hence there is no smoothening of the singular curve in the case $n=1$.
\end{rmk}

It is plausible that with the Lagrangian subvarieties $\Sigma, {\Sigma_{y_1,...,y_n}}\subset S=T^{\ast}C$ one can associate cyclic $DQ$-modules with cyclic vectors  $\psi_\Sigma$ and $\psi_{\Sigma_{y_1,...,y_n}}$ respectively.

\begin{que}
How to relate  $\psi_\Sigma$ and $\psi_{\Sigma_{y_1,...,y_n}}$ to $\psi_{\BB_{[\Sigma]}}$?

\end{que}

Probably the possibility to extract $\psi_\Sigma$ from $\psi_{\BB}$ is encoded into some constraints which allow us to extend differential forms
from $\Omega^1(\Sigma-\cup_{\alpha}\{r_\alpha\})$ to the whole curve $\Sigma$.
Relations of $\psi_\Sigma$ to the topological recursion were discussed in [DuMu].

\subsection{dg-Airy structure associated with compact Calabi-Yau $3$-fold}

In this subsection we discuss a differential-graded version  of an Airy structure
associated with the moduli space of complex Calabi-Yau $3$-folds. It is well-known that the analog of the above-discussed moduli space $\BB$ is a complex Lagrangian cone in the middle cohomology. The corresponding wave function contains all the information about Gromov-Witten invariants. Moreover, it satisfies the HAE.

Let $Y$ be a compact complex $3$-dimensional Calabi-Yau manifold (CY 3-fold for short).
Deformation theory of complex structures on $Y$ is controlled by the differential-graded Lie algebra (DGLA for short)
\[\g_Y=\Gamma(Y, \oplus_{i\ge 0}\wedge^iT^{1,0}_Y\otimes \Omega^{0,\ast}_Y)[1]\,.\]
We split $\g_Y$ into the direct sum of $\Z$-graded subspaces:
\[\g_Y=\oplus_{i\ge 0}\g^i\,,\quad \g_i=\g_Y^i=\Gamma(Y,\wedge^iT^{1,0}_Y\otimes \Omega^{0,\ast}_Y)[1]\,.\]
 The differential is $d_\g=\overline{\partial}+div$, where $div$
is the divergence operator corresponding to the holomorphic volume for $\Omega^{3,0}$ (see e.g. [BarKo]).
The corresponding homological vector field on $\g_Y[1]$ is given 
$$\dot{\gamma}=(\overline{\partial}+div)(\gamma)+{[\gamma,\gamma]\over{2}}$$
and the Maurer-Cartan equation (MC equation in short) is the equation for zeroes of this field. The space of solutions of the MC equation
 carries a foliation with the tangent space
 at $\gamma$ consisting of 
 \[\delta(\gamma)=(\overline{\partial}+div)(\epsilon)+[\epsilon,\gamma],\qquad \epsilon\in \g_Y\,.\]

The {\it extended moduli space of complex structures}  
is the quotient of the space of solutions to 
the MC-equation by this foliation.
The extended moduli space is a smooth finite-dimensional $\Z$-graded analytic supermanifold by a generalization of Tian-Todorov theorem (see loc.cit.).

Let us make the change of variables $\mu=e^{\gamma}$. Then the above homological vector field becomes
\begin{equation}\label{linear homfield}\dot{\mu}=(\overline{\partial}+div)(\mu)\end{equation}
and the foliation becomes also linear:
\[\delta(\mu)=(\overline{\partial}+div)(\epsilon),\qquad \epsilon\in \g_Y\,.\]

In the new variable $\mu$ the image of an element $\gamma_0+\gamma_1$ belonging to the graded Lie sublagebra $\g^0\oplus \g^1\subset \g_Y$
has the form $\exp(\gamma_0+\gamma_1)=\mu_0+\mu_1+\mu_2+\mu_3$ where:
$$\mu_0=\exp(\gamma_0),\, \mu_1=\exp(\gamma_0)\cdot \gamma_1,\, \mu_2=\exp(\gamma_0)\cdot(\gamma_1^2/2),\, \mu_3=\exp(\gamma_0)(\gamma_1^3/6).$$

 Using contraction with volume element $\Omega^{3,0}$ we identify (as a $\Z/2$-graded supermanifold) the
  vector space $\g_Y[1]$ endowed with the linear homological  vector field \eqref{linear homfield} with the space of forms
  \[W:=\Gamma(Y,\oplus_{p.q\ge 0}\Omega^{p,q}_Y)[3]\]
  endowed with the  linear homological vector field $\dot{\alpha}=d\alpha$.
  The space $W$ is an infinite-dimensional symplectic vector space endowed with the symplectic form
  $\langle \omega_1,\omega_2\rangle=\int_X\omega_1\wedge \omega_2$.

\begin{prp} Let us define ${\mathcal L}=\exp(\g^0\oplus \g^1)$ considered as a subvariety of $W$. Then:

1) ${\mathcal L}$ is a quadratic Lagrangian subvariety in $W$.

2) The vector space $Ker(d)\subset W$ is a coisotropic vector subspace and $Im(d)$ is symplectic orthogonal to $Ker(d)$. 

3) The natural projection (symplectic reduction) $pr: {\mathcal L}\cap Im(d)\to H^{\bullet}(Y)\to H^3(Y,\C)$ 
maps ${\mathcal L}$
into the Lagrangian cone $L\subset H^3(Y,\C)$ which is the moduli space of pairs consisting of a complex structure on $X$ and a choice
of holomorphic volume form $vol_Y$.

\end{prp}

{\it Proof.} Parts 2) and 3) follow from the deformation theory of the Calabi-Yau structure on $X$, and can be derived in many ways 
(see e.g. [BarKo]).

The fact that $W$ is symplectic is general. It holds for any real compact manifold $Y$ such that $dim_{\R}Y\in 2+4\Z_{\ge 0}$ (e.g. for any 
odd-dimensional compact complex manifold).

Notice  that  the Calabi-Yau
pairing  implies the following isomorphisms (properly understood in the infinite-dimensional context):
$\g^0\simeq (\g^3)^{\ast}, \g^1\simeq (\g^2)^{\ast}$. 
Quadratic equations for ${\mathcal L}$ require passing to the coordinates $\mu_i, i=0,1,2,3$, in which
${\mathcal L}$ can be described by the equations
$$\mu_0\mu_2={\mu_1^2\over{2}},$$
$$\mu_0\mu_3={\mu_1\mu_2\over{3}}.$$

The fact that ${\mathcal L}$ is Lagrangian can be checked directly.
Namely, one checks first that the dimension ${\mathcal L}$ understood in the sense of infinite-dimensional smooth topology on $\g$ is
``one half'' of the dimension of $\g$.
Then one shows by a direct computation, that the algebra of 
functions vanishing on ${\mathcal L}$ is closed under the Poisson bracket. Such functions have form 
$$f_{a_0,a_1}(\mu_0,...,\mu_3)=
Tr\left(a_0(\mu_0\mu_2-\mu_1^2/2)+a_1(\mu_0\mu_3-\mu_1\mu_2/3\right)),$$
where $a_i\in \g^i, i=0,1$ and $Tr:\g^3\to \C$ is the functional given by the holomorphic volume form $vol_X$.
Then one computes $\{f_{a_0,a_1},f_{b_0,b_1}\}$ and sees that it gives a solvable Lie algebra structure on functions vanishing on ${\mathcal L}$.

This completes the proof.
$\blacksquare$

We should warn the reader that despite the existence of quadratic Lagrangian, there is no Airy structure in the compact case. The corresponding tensors $A,B,C$ belong to the {\it completed} tensor products. As a result the series which should give the wave function are divergent and do not give a meaningful answer.

\vspace{3mm}

{\bf References}

\vspace{2mm}

[ABCO] J. Andersen, G. Borot, L. Chekhov, N. Orantin, The ABCD of topological recursion, in progress.

\vspace{2mm}

[ACNP] J.E. Andersen, L.O. Chekhov, P. Norbury, R.C. Penner, Topological
recursion for Gaussian means and cohomological field theories, Theor.
Math. Phys. 185 (2015), 1685-1717.

\vspace{2mm}

[AgKas]  A. D'Agnolo, M. Kashiwara, On quantization of complex symplectic manifolds, arXiv:1008.5273.

\vspace{2mm}

[BarKo] S. Barannikov, M. Kontsevich, Frobenius Manifolds and Formality of Lie Algebras of Polyvector Fields,  arXiv:alg-geom/9710032.

\vspace{2mm}

[BCOV] M. Bershadsky, S. Cecotti, H. Ooguri and C. Vafa, Kodaira-Spencer theory of gravity and
exact results for quantum string amplitudes, Comm. Math. Phys. 165(2), 311–427 (1994).

\vspace{2mm}

[Cos] K. Costello, The Gromov-Witten potential associated to a TCFT, arXiv:0509264.

\vspace{2mm}

[DuMu] O. Dumitrescu, M.Mulase, Quantization of spectral curves for meromorphic Higgs bundles through topological recursion,
arXiv:1411.1023.

\vspace{2mm}

[EynOr1]  B. Eynard, N. Orantin,   Invariants of algebraic curves and topological expansion, arXiv:math-ph/0702045.
   
 \vspace{2mm}

[EynOr2]  B. Eynard, N. Orantin, About the $x-y$ symmetry of the $F_g$ algebraic invariants, arXiv:1311.4993.

\vspace{2mm}
[GelKaz] I. Gelfand, D. Kazhdan, Some problems of differential geometry and the calculation of the cohomologies
of Lie algebras of vector fields, Soviet Math. Doklady 12, 1367-1370, 1971.

\vspace{2mm}

[KoSo] M. Kontsevich, Y. Soibelman, Wall-crossing structures in Donaldson-Thomas invariants, integrable systems and Mirror Symmetry,  arXiv:1303.3253.

\vspace{2mm} [T] J. Tate, Residues of differentials on curves, Ann. Esci. ENS, ser. 4, 1:1,  p. 149-159, 1968.

\vspace{2mm}

\vspace{2mm} [Wit] E. Witten, Quantum background independence in String Theory, hep-th 9306122.

\vspace{10mm}

M.K. IHES, 35 route de Chartres, Bures-sur-Yvette, F-91440, France, {email: maxim@ihes.fr}

Y.S. Department of Mathematics, Kansas State University, Manhattan, KS 66506, USA, {email: soibel@math.ksu.edu}
\end{document}